\documentclass[final]{siamltex}
\usepackage{latexsym, graphicx, epsfig, amsmath, amsfonts,amssymb,subfigure}
\usepackage{multirow}
\usepackage{color}

\newtheorem{remark}{Remark}
\usepackage[percent]{overpic}
\usepackage{booktabs,epstopdf}
\usepackage{algorithm}
\usepackage{hyperref}
\hypersetup{hypertex=true}
\usepackage{algorithmicx}
\usepackage{algpseudocode}
\floatname{algorithm}{Algorithm}
\def\mb{\mathbf}
\def\newx{\mathbf{x}}

\title{Failure-informed adaptive sampling for PINNs, Part II: combining with re-sampling and subset simulation}

\author
{
Zhiwei Gao\thanks{School of Mathematics, Southeast University, Nanjing 210096, China.}
\and Tao Tang\thanks{Division of Science and Technology, BNU-HKBU United International College, Zhuhai, China. Email: tangt@sustech.edu.cn. }
\and Liang Yan\thanks{School of Mathematics, Southeast University, Nanjing 210096; Nanjing Center for Applied Mathematics, Nanjing 211135, China. Email: yanliang@seu.edu.cn. LY's work was supported by the NSF of China (No.12171085).}
\and Tao Zhou\thanks{Institute of Computational Mathematics, Academy of Mathematics and Systems Science, Chinese Academy of Sciences, Beijing 100190, China. Email: tzhou@lsec.cc.ac.cn. This work is supported by the National Key R\&D Program of China (2020YFA0712000), NSF of China (under grant numbers 12288201), the Strategic Priority Research Program of Chinese Academy of Sciences (Grant No. XDA25010404), and youth innovation promotion association (CAS). }
}

\begin{document}
\graphicspath{figures/}
\maketitle

\begin{abstract}
This is the second part of our series works on failure-informed adaptive sampling for physic-informed neural networks (FI-PINNs). In our previous work \cite{gao2022failure}, we have presented an adaptive sampling framework by using the failure probability as the posterior error indicator, where the truncated Gaussian model has been adopted for estimating the indicator. In this work, we present two novel extensions to FI-PINNs. The first extension consist in combining with a re-sampling technique, so that the new algorithm can maintain a constant training size.  This is achieved through a cosine-annealing, which gradually transforms the sampling of collocation points from uniform to adaptive via training progress. The second extension is to present the subset simulation algorithm as the posterior model (instead of the truncated Gaussian model) for estimating the error indicator, which can more effectively estimate the failure probability and generate new effective training points in the failure region. We investigate the performance of the new approach using several challenging problems, and numerical experiments demonstrate a significant improvement over the original algorithm.
    \end{abstract}

\pagestyle{myheadings}
\thispagestyle{plain}
\markboth{Z. Gao, T. Tang, L. Yan and T. Zhou}
{AFI-PINNS}

\section{Introduction}

Scientific machine learning (Sci-ML) \cite{E_CICP,Edeep18,raissi2019physics,sirignano2018dgm,zang2020weak} has gained popularity in recent years for solving partial differential equations (PDE). The Sci-ML methods are different from conventional ML methods in that they can take into account physical information from the system. Additionally, similar to general ML methods, Sci-ML methods can be data-driven by providing a large amount of training data. Although it is possible to use both physical and data information, it still presents a significant challenge that must be carefully exploited.

Physic-informed neural networks(PINNs) \cite{lu2021deepxde,raissi2019physics} have emerged as a potential solution to the aforementioned problem. The key approach of PINNs is to embed the physical information, i.e., the PDE residual function induced by the neural network, into the loss function as a soft penalty. As a result, the PDE solution is transformed into an unconstrained optimization problem by minimizing the self-supervision loss function. However, training PINNs can be difficult in some cases. This is due in part to the unbalance of different loss terms, several of which may contain high order derivatives and discontinuities inherited from the PDE \cite{krishnapriyan2021characterizing,wang2021understanding,wang2022and}.

To address this issue, independent solutions such as sequence to sequence learning \cite{krishnapriyan2021characterizing}, training weights for different loss terms \cite{bischof2021multi,mcclenny2020self,xiang2022self}, and adaptive sampling strategies \cite{daw2022rethinking,lu2021deepxde,peng2022rang,subramanian2022adaptive,tang2021deep,wu2023comprehensive}have been proposed. It is well noticed that the performance of PINNs  greatly depends on the distribution and location of the \textit{collocation} points. Adaptive sampling techniques are the most effective methods for selecting collocation points in practice. Initially, \cite{raissi2019physics} proposed a uniform selection of collocation points from the system domain. In most cases, this approach falls short of capturing the most crucial details of the real underlying problem. The residual-based refinement method (RAR) proposed in \cite{lu2021deepxde} is a better option because it selects adaptive samples based on the residual function values induced by the network. However, this method still has the drawback of the effective sample size not being large enough to precisely locate the area with larger residual errors. In our recent work \cite{gao2022failure}, we presented a novel mathematical framework, namely failure-informed PINNs (FI-PINNs), which can integrate various adaptive sampling strategies for PINNs. FI-PINNs is inspired by adaptive finite element methods for solving PDEs. The basic idea is to use  a performance function, e.g. the PDE residual, to determine a region known as {\it failure region}, where the learning accuracy is relatively low.  Then we can define the {\it failure probability} as an error indicator to describe the reliability of the PINNs based on this failure region. In particular, we have proposed to use the truncated Gaussian as a simple posterior model to estimate the failure probability and to generate new training points within the failure region. Empirical studies on the benchmark  problems showed that FI-PINNs can effectively capture the solution structure.

In the current work, we present several novel extensions of FI-PINNs. More precisely:
 \begin{itemize}
 \item We describe an efficient annealing FI-PINNs that re-sampling the collocation dataset in an annealing manner while maintaining the same dataset size during the training phase to reduce computational costs. This is achieved through a cosine-annealing, which gradually transforms the sampling of collocation points from uniform to adaptive via training progress.

\item We present \textit{subset simulation} as the posterior model, instead of the truncated Gaussian model, to estimate the failure probability and to generate new training points. It is shown that the new model is more flexible for complex solutions, e.g., solutions with several peaks.

\item We investigate the role of various performance function $\mathcal{Q}$ in the training. In particular, besides the residue function itself, we also investigate the use of its gradient.

\item We present several numerical experiments to demonstrate that the new annealing FI-PINNs can further accelerate training convergence. Furthermore, we show that the new method outperforms the original algorithm.
\end{itemize}

The remainder of this paper is structured as follows.   The fundamental principles of PINNs and FI-PINNs are introduced in Section \ref{background}.   In Section \ref{annealing_adaptive_sampling}, we present the annealing failure-informed PINNs and the subset simulation strategy.  We provide numerical experiments to verify the effectiveness of our sampling strategy in Section \ref{numerical_experiments}, and this is followed by some concluding remarks in Section \ref{conclusion}.

\section{Background}
\label{background}
In this section, we first provide a brief overview of the PINNs for solving PDEs.  The failure-informed PINNs(FI-PINNs) framework will then be introduced.
\subsection{Physics-informed neural networks}
\label{PINNs}
Consider the following form of partial differential equations(PDE)
\begin{equation}
\begin{split}
\mathcal{N}(\newx,u) &= 0,\quad \newx\in \Omega,\\
\mathcal{B}(\newx,u)&=0,\quad\newx\in \partial \Omega,
\end{split}
\end{equation}
where $\mathcal{N}$ is a linear or nonlinear partial differential operator that operates on the variable $\newx$. $\mathcal{B}$ is the operator that constrain the boundary condition to have a well-posed solution. The key approach of PINNs is to use a deep neural network (DNN) $u(\newx;\theta)$ parameterized by $\theta$ to approximate the PDE solution. Thus, solving the PDE is equivalent to solving the following unconstrained optimization problem:
\begin{equation}
\label{loss_function}
\theta^{*} = \arg\min_{\theta \in \Theta}\mathcal{J}(\theta),
\end{equation}
where $\Theta$ is a fixed parameter space.
In general, the loss function $\mathcal{J}(\theta)$ consists of two parts, the PDE loss $\mathcal{J}_c(\theta)$ and the boundary loss $\mathcal{J}_b(\theta)$, which are defined as
\begin{equation}
   \begin{split}
   \mathcal{J}_{c}(\theta) &= \|r(\newx;\theta)\|^2_{\Omega}:=\int_{\Omega}|r(\newx;\theta)|^{2}d\newx ,\\
   \mathcal{J}_{b}(\theta) &= \|\mathcal{B}(\newx,u(\newx;\theta))\|^2_{\partial \Omega}:=\int_{\partial \Omega}|\mathcal{B}(\newx,u(\newx;\theta))|^{2}d\newx ,
   \end{split}
\end{equation}
where
\begin{equation}
\label{residual_function}
r(\newx;\theta) := \mathcal{N}(\newx,u(\newx;\theta))
\end{equation} denotes the PDE residual function induced by the neural network. It is noted that the $\mathcal{J}_{c}(\theta)$ incorporates the physical information, while $\mathcal{J}_{b}(\theta)$ makes use of the domain information of the boundary parts. This is consistent with our previous explanation. In practical, the loss function is often discretized by given a training dataset $\mathcal{D}$, consisting of collocation dataset $\{\newx_i^c\}_{i=1}^{N_c} \subset \mathcal{D}_{c}$, boundary dataset $\{\newx_i^b\}_{i=1}^{N_b}  \subset\mathcal{D}_{b}$. Therefore, the loss function can be discretized as
\begin{equation}
\begin{split}
\hat{\mathcal{J}}_{c}(\theta) &= \frac{1}{N_{c}}\sum_{i=1}^{N_{c}}|r(\newx^c_{i};\theta)|^{2},\\
\hat{\mathcal{J}}_{b}(\theta) &= \frac{1}{N_{b}}\sum_{i=1}^{N_{b}}|\mathcal{B}(\newx^b_{i},u(\newx^b_{i};\theta))|^{2},
\end{split}
\end{equation}
where $N_{c}, N_{b}$ are the number of points for each dataset.
Then, the optimal parameters can be obtained by minimizing the following discrete loss function
\begin{equation}
\label{loss_term}
\hat{\theta}^{*} = \arg\min_{\theta\in \Theta}(\hat{\mathcal{J}}_{c}(\theta)+ \lambda_{b}\hat{\mathcal{J}}_{b}(\theta)),
\end{equation}
where $\lambda_{b}$ is the weight of the boundary loss term.  It should be noted that utilizing various weights could produce diverse outcomes. As stated in \cite{bischof2021multi}, weights need to be carefully adjusted for a variety of problems in order to perform better.

\subsection{Failure-informed PINNs}
\label{Fi-PINNs}
It is known that conventional PINNs may fail to converge due to the existence of several failure modes\cite{krishnapriyan2021characterizing,wang2021understanding,wang2022and}, including the curse of dimensionality, the low regularity, the discontinuities etc.  Some adaptive sampling techniques have been proposed to solve these probelms, see e.g., \cite{daw2022rethinking,lu2021deepxde,subramanian2022adaptive,tang2021deep,wu2023comprehensive}.  In \cite{gao2022failure}, we present a novel mathematical framework called failure-informed PINNs (FI-PINNs), which can combine various adaptive sampling strategies for PINNs. The basic idea is to define a non-negative performance function $\mathcal{Q}(\newx)$ to quantify the system's reliability from the viewpoint of reliability analysis. We define the system as unreliable at a point when $\mathcal{Q}(\newx)$ exceeds a predetermined tolerance $\epsilon_{r}$. The failure region $\Omega_{\mathcal{F}}$ is made up of these points collectively, and it is defined as
\begin{equation}
\Omega_{\mathcal{F}} :=\{\newx:\mathcal{Q}(\newx)>\epsilon_{r}\}.
\end{equation}
By creating adaptive samples from the failure region and updating the collocation dataset with these points, the performance of the PINNs can be enhanced after the network has been retrained. As a result, the failure region will shrink as the system becomes more stable and reliable.
We can define the failure probability $P_{\mathcal{F}}$ over the failure region $\Omega_{\mathcal{F}}$ as
\begin{equation}\label{faiP}
P_{\mathcal{F}} = \int_{\Omega}\omega(\newx)\mathcal{I}_{\Omega_{\mathcal{F}}}(\newx)d\newx  = \int_{\Omega_{\mathcal{F}}}\omega(\newx)d\newx,
\end{equation}
where $\omega(\newx)$  is the prior distribution of $\newx$, and $\mathcal{I}_{\Omega_{\mathcal{F}}}$ is the indicator function with
\begin{equation}
\label{indicator_function}
\mathcal{I}_{\Omega_{\mathcal{F}}}(\newx) =  \begin{cases}
1, \quad \newx\in \Omega_{\mathcal{F}},\\
0,\quad otherwise.
\end{cases}
\end{equation}
Because the failure probability typically decreases as the failure region shrinks, it can be used to automatically halt training when it is smaller than a tolerance $\epsilon_p$.  In some ways, the failure probability functions  $P_{\mathcal{F}}$ as an error indicator for training set refinement, similar to the classical adaptive FEM. Additionally, as demonstrated in \cite{gao2022failure}, these two tolerances, $\epsilon_r$ and $\epsilon_p$, can  ensure the convergence of the FI-PINNs framework after the training is finished.

One benefit of the FI-PINNs algorithm is that the candidate points for updating can be obtained simultaneously when the approximate method, particularly the sampling method, is used to estimate the failure probability.  Many adaptive sampling techniques can therefore be unified within this framework. The residual-based refinement method  proposed in \cite{lu2021deepxde}, for instance, can be viewed as approximating the failure probability using Monte Carlo method. For the sake of efficiency,  in \cite{gao2022failure} we choose  a truncated Gaussian as a simple model from the importance sampling respective in order to approximatively determine the failure probability and to produce new training points. It should be noted that selecting the appropriate $\mathcal{Q}(\newx)$ also has a significant impact on performance. Different proposals have been made at the moment, such as the residual function $r(\newx;\theta)$ and its gradient function  defined as
\begin{equation}
\label{gradient_residual_function}
r_{\newx}(\newx;\theta) = \frac{\partial r(\newx;\theta)}{\partial \newx}.
\end{equation}
In this paper, we decide to use the aforementioned functions as the  $\mathcal{Q}$ performance function.

\begin{remark}
   Note that the residual function may be very large in some areas due to the low regularity of the solution. Thus, the residual function is normalized with respect to the first update $\|r_{1}(\newx;\theta)\|_{\Omega}$, that is $\hat{r}_{k}(\newx;\theta) = \frac{r_{k}(\newx;\theta)}{\|r_{1}(\newx;\theta)\|_{\Omega}}$, where $r_{k}(\newx;\theta)$ is the residual function at the $k$th update in the framework FI-PINNs. This technique is also used with the gradient function.
\end{remark}

\begin{remark}
   While FI-PINNs has been shown to be more effective over the original PINNs \cite{gao2022failure}, there are, however, many possible ways to further improve its performance. For instance, the truncated Gaussian model can be changed into other models so that it can handle more complex solutions structures. Moreover, in the original FI-PINNs, the adaptive updating rule is to add new collocation points into the previous training dataset, meaning that the training set will become larger and larger, as the adaptivity proceed. This can also be improved by re-sampling techniques so that we can control the size of training dataset. These will be discussed in detail in the next section.
\end{remark}

\section{FI-PINNs: combining with re-sampling and  subset simulation}
\label{annealing_adaptive_sampling}

 We present in this section two main extensions to FI-PINNs, namely, the re-sampling technique and the subset simulation algorithm.

\subsection{Annealing failure-informed PINNs}

As mentioned in the last section, in the original FI-PINNs, the adaptive updating rule is to add new collocation points into the previous training dataset, meaning that the training set will become larger and larger, as the adaptivity proceed.  To overcome this limitation, we propose an efficient annealing-based adaptive sampling method to update the collocation dataset within the FI-PINNs framework,  motivated by \cite{daw2022rethinking,subramanian2022adaptive}.  The proposed method is henceforth refereed as annealing failure-informed PINNs.  The new algorithm consists of two key parts. First, the algorithm decides whether or not to update the collocation points while training, rather than waiting until training is finished.  The second important distinction is that only the components of the training set are changed; the size of the training set itself does not change during training.  We will demonstrate that, especially for problems with local features, this approach results in better localization of the collocation points.

\begin{figure}[htbp]
   \centering
   \includegraphics[width = 0.7\textwidth]{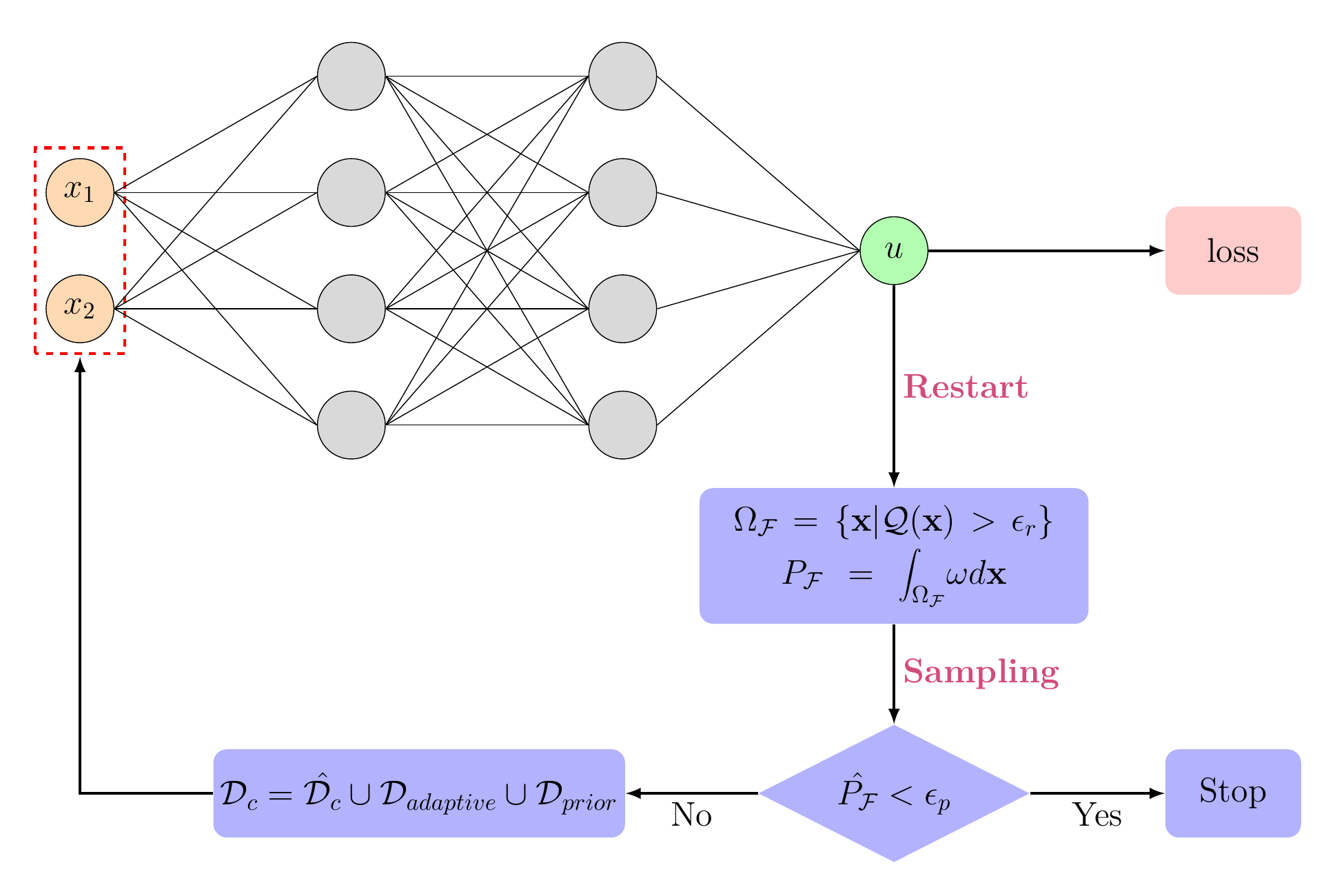}
   \caption{The framework of AFI-PINNs}
   \label{Annealing-FIPINNs}
   \end{figure}

We now go into more detail about how the annealing framework we described in the previous description was put into practice. The schematics of  annealing FI-PINNs  is shown in Fig. \ref{Annealing-FIPINNs}.

The annealing FI-PINNs framework is broken into two stages. The first  is to determine  when to update the training set, which is referred to as the {\it restart} phase. The second is  how to resample the point specifically. For the first phase, an easy solution is to swap points on a regular schedule, i.e., we restart every given number of epochs throughout training.  In this work, we will resample the collocation points when  the training loss is no longer decreasing for a few epochs.  We now going to  introduce our  adaptive sampling procedure using an annealing manner.  When a restart is required at the  epoch $T_{restart}$, we   estimate the  failure probability $\hat{P}_{\mathcal{F}}$  first using a sampling method, such as  the SAIS method suggested in \cite{gao2022failure} or the subset simulation we will use in this paper, and  simultaneously generate the candidate adaptive dataset $\mathcal{D}_{adaptive}$ of $N_{\mathcal{F}}$ collocation points from the failure region.   If $\hat{P}_{\mathcal{F}}$ is less than the  predefined tolerance $\epsilon_{p}$,  the entire training can be terminated according to our FI-PINNs framework.  Otherwise, we update the training points using an annealing manner.  For the specifics of the updating rule, the new  collocation dataset $\mathcal{D}_{c}$ will consist three components:  the adaptive collocation dataset $\mathcal{D}_{adaptive}$, the pre-selected collocation dataset $\hat{\mathcal{D}_{c}}$ with size of $N_{c}(1-\eta)$ generate from the initial dataset  $\mathcal{D}_{c}$,  and a new dataset $\mathcal{D}_{prior}$    with size $\eta N_{c} - N_{\mathcal{F}}$ samples from the prior distribution.  Here,  $0<\eta<1$ is the resample proportion.  Using this setting, the number of collocation points stays the same after updating.   The reason for including pre-selected points is that the network needs to retain the learned characteristics from previous training. While adaptive samples and prior samples help support the network in learning new characteristics from the failure region and other parts of the problem domain, respectively.  This type of combination will greatly improve training effectiveness while also reducing computational costs and speeding up convergence.

\begin{figure}[t]
\centering
\includegraphics[width = \textwidth]{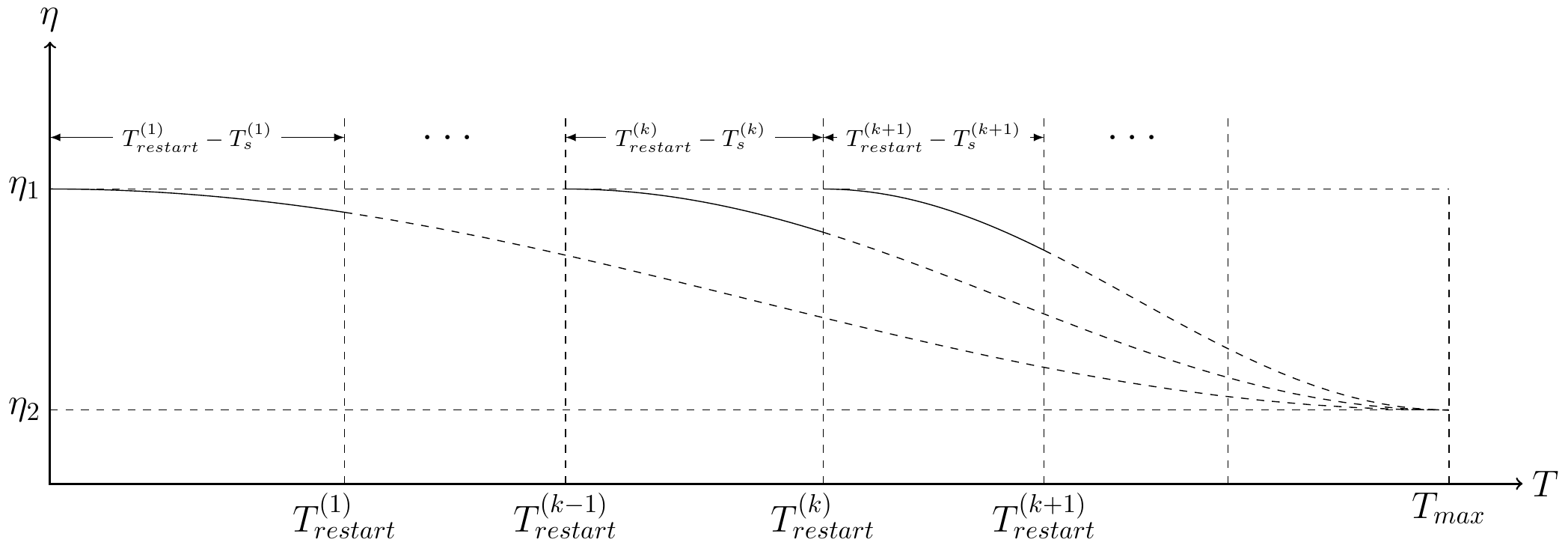}
\caption{The proportion of adding adaptive samples.  $T^{(k)}_{restart}$ is the epoch of the $k$-th restart,  $T^{(k)}_{s}=T^{(k-1)}_{restart}$ is the epoch of the last restart.}
\label{proportion}
\end{figure}

To determine the resample proportion $\eta$,  a cosine annealing manner \cite{subramanian2022adaptive} is used based on the restart  and maximum epochs. This trick aims to select the proportion adaptively in order to maintain  balance between different parts of new dataset.  In more detail, suppose the $k$-th restart epoch is $T_{restart}^{(k)}$, and the total training epoch is $T_{max}$. When the  training is restarted, the proportion can be calculated as
\begin{equation}
\label{Annealing}
\eta = a\left(1+b\cos\left(\frac{\pi (T_{restart}^{(k)} - T_{s}^{(k)})}{T_{max} - T_{s}^{(k)}}\right)\right),
\end{equation}
where $a,b$ are positive constants used to control the upper bound and lower bound of the proportion, $T_{s}^{(k)}$ is the epoch of last restart, i.e., $T_{s}^{(k)}=T_{restart}^{(k-1)}$. Fig.\ref{proportion} demonstrates the trend for the proportion during the training process, where $\eta_{1} = a(1+b)\in(0,1), \eta_{2} = a(1-b)\in (0,1)$.  In general, if $N_{\mathcal{F}} < N_{c}\eta:=N_s$, we can draw $N_{s} - N_{\mathcal{F}}$ samples from the prior, denoted as $\mathcal{D}_{prior}$. Otherwise, we can construct the $\mathcal{D}_{adaptive}$ by draw the  $N_{s}$ samples from the candidate adaptive dataset. In this case, $\mathcal{D}_{prior} = \emptyset$. This procedure will be repeated until the training epoch reaches the maximum epoch or the failure probability in some restart epoch is less than the tolerance.   The pseudo-code of our proposed annealing sampling strategy is shown in Algorithm \ref{Algorithm2}.

 The efficient annealing sampling strategy greatly accelerates training convergence. In detail, at each restart, the loss will typically increase suddenly and then decrease quickly with a smaller predicted error. This allows for the use of fewer samples to achieve the same or higher accuracy. The novel annealing FI-PINNs framework can also reduce computational costs by keeping the size of the collocation dataset invariant while terminating early based on the failure probability. This process differs from traditional PINNs, which train on the same dataset throughout. Instead, without increasing in number, the training points gradually transition through an annealing process from the initial uniform training pool to the final collocation points in failure regions.
This updating rule has superiority over the  FI-PINNs as it can achieve higher accuracy with a smaller training dataset, which will be verified in our experiments.

 \begin{algorithm}[t]
\caption{Annealing FI-PINNs}
\label{Algorithm2}
\begin{algorithmic}[1]
\Require  ~ A DNN network $u(\mb{x};\theta)$, collocation dataset $\mathcal{D}_{c}$, boundary dataset $\mathcal{D}_{b}$, number of collocation points $N_{c}$, maximum epochs $T_{\max}$, performance tolerance $\epsilon_{r}$, failure probability tolerance $\epsilon_{p}$.
\State $s\leftarrow 1, T_{s} \leftarrow 0$
\While  {$s \le T_{\max}$}
\State Train $u(\mb{x}; \theta)$ using the training dataset $\mathcal{D}_{b}, \mathcal{D}_{c}$
\If {\textbf{Restart}}
\State $T_{restart}\leftarrow s$
\State $\eta \leftarrow \text{cosine-schedule($T_{restart}$, $T_{s}$, $T_{max}$)}$ using Eq.\eqref{Annealing}
\State Estimate failure probability $\hat{P}_{\mathcal{F}}$ using sampling method, and get $N_{\mathcal{F}}$ candidate points $\{\newx_{adpative}:  \newx_{adpative}\in \Omega_{\mathcal{F}}\}$
\If{$\hat{P}_{\mathcal{F}}<\epsilon_{p}$}
\State \textbf{Early stop}
\Else
\State Set the adaptive dataset $\mathcal{D}_{adaptive} = \{\newx_{adpative}\}$
\State Uniformly generate $(1-\eta)N_{c}$  points $\{\newx_{pr}: \newx_{pr}\in \mathcal{D}_{c}\}$, set $\hat{\mathcal{D}}_{c}= \{\newx_{pr}\}$
\If {$N_{\mathcal{F}} < \eta N_{c}$}
\State Generate  $\eta N_{c} - N_{\mathcal{F}}$  points $\{\newx_{r}:\newx_{r} \sim \omega(\newx)\}$, set $\mathcal{D}_{prior}=\{\newx_{r}\}$
\Else
\State Set $\mathcal{D}_{prior} = \emptyset$
\EndIf
\State $\mathcal{D}_{c} \leftarrow  \mathcal{D}_{adaptive}  \cup \hat{\mathcal{D}}_{c} \cup \mathcal{D}_{prior}$
\EndIf
\State $T_{s}\leftarrow T_{restart}$
\EndIf
\State $ s \leftarrow s + 1$
\EndWhile
\end{algorithmic}
\end{algorithm}

\begin{remark}
Note that the different choice of restart criteria for the training of PINNs can be applied. For example, one can try restarting according to the test error. Different annealing methods can also be applied, such as sine annealing, exponential annealing etc.
\end{remark}

\subsection{Subset simulation}
\label{subset_simulation}

It is worthy noting that our novel framework can use a variety of sampling techniques to produce the adaptive dataset. In this paper, we use the subset simulation (SS)\cite{au2001estimation,zuev2015subset} based on MCMC method to accomplish this task.   SS is one of the most commonly used  variance reduction techniques for developing an efficient estimator for failure probability $\hat{P}_{\mathcal{F}}$ in Eq.(\ref{faiP}). It has been successfully used in wide range of contexts and  applications.  The main idea is to define a series of nested failure regions that include the target failure region $\Omega_{\mathcal{F}}$,
\begin{equation}
\Omega_{\mathcal{F}}  \subset \Omega_{\mathcal{F}_{m}}\subset \Omega_{\mathcal{F}_{m-1}} \subset \ldots \subset \Omega_{\mathcal{F}_{1}} \subset \Omega_{\mathcal{F}_{0}} = \Omega.
\end{equation}
Then, the failure probability can be divided into
\begin{equation}
\begin{split}
P_{\mathcal{F}} = \mathcal{P}(\Omega_{\mathcal{F}})
&= \mathcal{P}(\Omega_{\mathcal{F}_{1}}|\Omega_{\mathcal{F}_{0}})\mathcal{P}(\Omega_{\mathcal{F}_{2}}|\Omega_{\mathcal{F}_{1}})\cdots \mathcal{P}(\Omega_{\mathcal{F}_{m}}|\Omega_{\mathcal{F}_{m-1}})\mathcal{P}(\Omega_{\mathcal{F}}|\Omega_{\mathcal{F}_{m}}),
\end{split}
\end{equation}
where $\mathcal{P}(\Omega_{\mathcal{F}_{k+1}}|\Omega_{\mathcal{F}_{k}}) = \frac{\mathcal{P}(\Omega_{\mathcal{F}_{k+1}})}{\mathcal{P}(\Omega_{\mathcal{F}_{k}})}, 0\leq k\leq m-1$ denote the conditional probability given the intermediate failure region $\Omega_{\mathcal{F}_{k}}$, $\mathcal{P}(\Omega_{\mathcal{F}}|\Omega_{\mathcal{F}_{m}})$ is the final conditional probability. Once the necessary preparations have been made, the main challenge is to determine the intermediate failure regions and conditional probabilities.

In the procedure of subset simulation, the conditional probability is set to a constant intermediate failure probability $p\in(0,1)$.  It is assumed that  $N_{p} = N_{s}p$ and $p^{-1}$ are both positive integers. These are equal to the number of chains and the number of samples per chain at a specific  $p$, as will be seen in the following section.  In the SS literature, an advisable choice is $p=0.1$. The number of samples $N_s$ controls the statistical accuracy of results.  A simulation run begins with the unconditional probability, which generates $N_s$ i.i.d. samples of $\newx$  from $\omega(\newx)$ using direct Monte Carlo method. The corresponding $\mathcal{Q}$ values are computed and arranged in descending order, resulting in an ordered list denoted by $\mathcal{S}_{0} = \{\newx_{0}^{(i)}\}_{i=1}^{N_{s}}$.  Then the intermediate failure region $\Omega_{\mathcal{F}_{1}}:= \{\newx:\mathcal{Q}(\newx)>\epsilon^{(1)}_{r}\}$, where $\epsilon^{(1)}_{r}$ is determined as the $(N_p+1)$-th largest sample values of $\newx$ at level $0$, i.e.,
$ \epsilon^{(1)}_{r} = \newx_{0}^{(N_p+1)}.$

We can obtain the intermediate failure region $\Omega_{\mathcal{F}_{k+1}}$ by employing a similar method. Suppose the samples set $\mathcal{S}_{k} = \{\newx_{k}^{(i)}\}_{i=1}^{N_{s}}$ are from the conditional distribution $\rho(\cdot|\Omega_{k})$ and have been rearranged according to their $\mathcal{Q}$ values in descending order. Then the ($k+1$)th \textit{intermediate failure region} $\Omega_{\mathcal{F}_{k+1}}$ is defined as
\begin{equation}
   \label{middle_failure_region}
\Omega_{\mathcal{F}_{k+1}} := \{\newx:\mathcal{Q}(\newx)>\epsilon^{(k+1)}_{r}\},
\end{equation}
where $\epsilon^{(k+1)}_{r}$ is determined as the $(N_p+1)$-th largest sample values of $\newx$ at level $k$, i.e.,
$$ \epsilon^{(k+1)}_{r} = \newx_{k}^{(N_p+1)}.$$
Thus, the conditional probability $\mathcal{P}(\Omega_{\mathcal{F}_{k+1}}|\Omega_{\mathcal{F}_{k}})$ can be estimated as
\begin{equation}
\mathcal{P}(\Omega_{\mathcal{F}_{k+1}}|\Omega_{\mathcal{F}_{k}}) = \int_{\Omega_{\mathcal{F}_{k+1}}} \rho(\newx|\Omega_{\mathcal{F}_{k}})d\newx\approx \frac{1}{N_{s}}\sum_{i=1}^{N_{s}}\mathcal{I}_{\Omega_{\mathcal{F}_{k+1}}}(\newx_{k}^{(i)}) = p.
\end{equation}
It should be noted that the samples $\{\newx_{k}^{(i)}\}_{i=1}^{N_{p}}$ also follow the conditional density $\rho(\cdot|\Omega_{\mathcal{F}_{k+1}})$.  For this reason, in order to keep the size of the samples constant, the next task is to generate the final $(N_{s}-N_{p})$ samples from $\rho(\cdot|\Omega_{\mathcal{F}_{k+1}})$.  In practice, the \textit{modified Metropolis Hasting} algorithm (MMA) \cite{zuev2015subset} can be used to generate $N_{p}$ chains with initial \textit{seeds} $\{\newx_{k}^{(i)}\}_{i=1}^{N_{p}}$.  There will be $\frac{1}{p}-1$ new samples accepted in each chain during the sampling procedure. Therefore, the next samples set $\mathcal{S}_{k+1}$ can be formed as $\{\newx_{k+1}^{(i)}\}_{i=1}^{N_{s}}= \{\newx_{k+1}^{(i)}\}_{i=1}^{N_{s}- N_{p}}\cup \{\newx_{k}^{(i)}\}_{i=1}^{N_{p}}$, where $\{\newx_{k+1}^{(i)}\}_{i=1}^{N_{s} - N_{p}}$ are samples in the $N_{p}$ chains.  The aforementioned procedure of generating  MCMC samples and increasing the simulation level $k$ is repeated until the desired threshold level  is reached.  To this end,   at level $k$,  the number of failure samples is calculated as
\begin{equation}
\label{failure_samples_calculation}
N_{\mathcal{F}_{k}} = \sum_{i=1}^{N_{s}}\mathcal{I}_{\Omega_{\mathcal{F}}}(\newx_{k}^{(i)}).
\end{equation}
The simulation can be terminated if $N_{p}>N_{\mathcal{F}_{k}}$ as the intermediate failure region is close to the real failure region. If the simulation ends at $k = m$, then the conditional probability  $\mathcal{P}(\Omega_{\mathcal{F}}|\Omega_{\mathcal{F}_{m}})$ can be estimated as
\begin{equation}
\mathcal{P}(\Omega_{\mathcal{F}}|\Omega_{\mathcal{F}_{m}}) = \int_{\Omega_{\mathcal{F}}}\rho(\newx|\Omega_{\mathcal{F}_{m}})d\newx \approx \frac{1}{N_{s}}\sum_{i=1}^{N_{s}}\mathcal{I}_{\Omega_{\mathcal{F}}}(\newx_{m}^{(i)}) = \frac{N_{\mathcal{F}_{m}}}{N_{s}}:=q.
\end{equation}
Finally, the failure probability can be approximated as
\begin{equation}
\label{failure_probability}
\begin{split}
P_{\mathcal{F}} &=  \mathcal{P}(\Omega_{\mathcal{F}_{1}}|\Omega_{\mathcal{F}_{0}})\mathcal{P}(\Omega_{\mathcal{F}_{2}}|\Omega_{\mathcal{F}_{1}})\cdots \mathcal{P}(\Omega_{\mathcal{F}_{m}}|\Omega_{\mathcal{F}_{m-1}})\mathcal{P}(\Omega_{\mathcal{F}}|\Omega_{\mathcal{F}_{m}})\\
&\approx p^{m}q = \hat{P}_{\mathcal{F}}^{SS}.
\end{split}
\end{equation}

\begin{algorithm}[t]
\caption{Subset simulation}
\label{Algorithm1}
\begin{algorithmic}[1]
\Require   The performance function $\mathcal{Q}$, the prior distribution $\omega$, the updating proportion $\eta$, the number of collocation points $N_{c}$, level probability $p$.
\State $N_{s}\leftarrow N_{c}\eta, N_{p}\leftarrow N_{s}p$
\State Generate $\mathcal{S}_{0}= \{\newx_{0}^{(i)}\}_{i=1}^{N_{s}}$ from $\omega$ and rearrange $\mathcal{S}_{0}$ according to $\mathcal{Q}$ values in descending order
\State Calculate $N_{\mathcal{F}_{k}}$ using Eq.\eqref{failure_samples_calculation}
\State $k\leftarrow 0$
\While  {$N_{\mathcal{F}_{k}} \le N_{p}$}
\State Calculate the intermediate failure region $\Omega_{\mathcal{F}_{k+1}}$ using Eq.\eqref{middle_failure_region}
\State Generate $\mathcal{S}_{k+1}=\{\newx_{k+1}^{(i)}\}^{N_s}_{i=1}$ using MMA method with initial seeds $\{\newx_{k}^{(i)}\}_{i=1}^{N_{p}}$
\State Rearrange $\mathcal{S}_{k+1}$ according to $\mathcal{Q}$ values in descending order
\State Calculate $N_{\mathcal{F}_{k+1}}$ using Eq.\eqref{failure_samples_calculation}
\State $k \leftarrow k + 1$
\EndWhile
\State Get the estimated failure probability $\hat{P}_{\mathcal{F}}^{SS}$ using Eq.\eqref{failure_probability}\\
\Return The failure probability $\hat{P}_{\mathcal{F}}^{SS}$,  failure samples set $\{\newx_{k}^{(i)}\}_{i=1}^{N_{\mathcal{F}}}$.
\end{algorithmic}
\end{algorithm}

It is worth noting that the $\mathcal{D}_{adaptive}$ represents the failure samples $\{\newx_{m}^{(i)}\}_{i=1}^{N_{{\mathcal{F}}}}$ generated in the final approximation. The number of failure samples is generally $N_{\mathcal{F}}\leq N_{s}$, which is consistent with our annealing sampling procedure. Furthermore, because the subset simulation employs the MCMC method, it can handle multi-modal failure regions. This property is useful for dealing with PDEs with local features. The whole procedure of SS can be summarized in Algorithm \ref{Algorithm1}.

\section{Numerical experiments}
\label{numerical_experiments}

In this section, we shall illustrate the advantages of annealing FI-PINNs with subset simulation using four numerical tests, including the two-dimensional problem with two and four peaks, a time-dependent wave equation, and a ten-dimensional Poisson equation. To better present the results, we shall perform the following two types of approaches:

\begin{itemize}
\item The residual function defined in Eq.\eqref{residual_function} is selected as the performance function for  annealing FI-PINNs, denoted as ``R-FIPINN".
\item The gradient function of the residual defined in Eq.\eqref{gradient_residual_function} is select as the performance function  for  annealing FI-PINNs, denoted as ``G-FIPINN".
\end{itemize}
We also test the annealing adaptive PINNs for comparison, which only restart and do not decide whether to stop early based on the failure probability. We generate the adaptive dataset at each restart using the residual-based refinement method. We designated this scheme as ``MC-FIPINN".

 To measure the performance of an approximation, we will use the following  relative $L_{2}$ error, defined as
\begin{equation}
err_{L_{2}} = \frac{\sqrt{\sum_{i=1}^{N}|u(\newx_{i}) - \hat{u}(\newx_{i})|^{2}}}{\sqrt{\sum_{i=1}^{N}|u(\newx_{i})|^{2}}}.
\end{equation}
where $\{u(\newx_{i})\}_{i=1}^{N}$ and $\{\hat{u}(\newx_{i})\}_{i=1}^{N}$ are the exact  and  predicted solution,  respectively.

Unless otherwise stated, the network used in our experiments is a fully connected network with a \textit{Tanh}  activation function using 5 hidden layers, each with 64 neurons. The network was trained using the \textit{LBFGS} optimizers  with a learning rate of 0.3 and a maximum of 50,000 iterations.  Additionally, the two parameters used in Eq.(\ref{Annealing}) are $a = 0.5, b  =1$ and the weights in Eq.\eqref{loss_term} are set to 1. The performance tolerance $\epsilon_{r}$ is set to 0.1, and the failure probability tolerance is set to 0.01.  The number of boundary points is set to 800. We select the size of the collocation dataset to be 500, 1000, 1500, and 2000 respectively for each experiment to test whether fewer collocation samples can achieve high accuracy.

\subsection{Two-dimensional problem with multiple peaks}

\begin{figure}[t]
\centering
(a)\includegraphics[width = 0.4\textwidth]{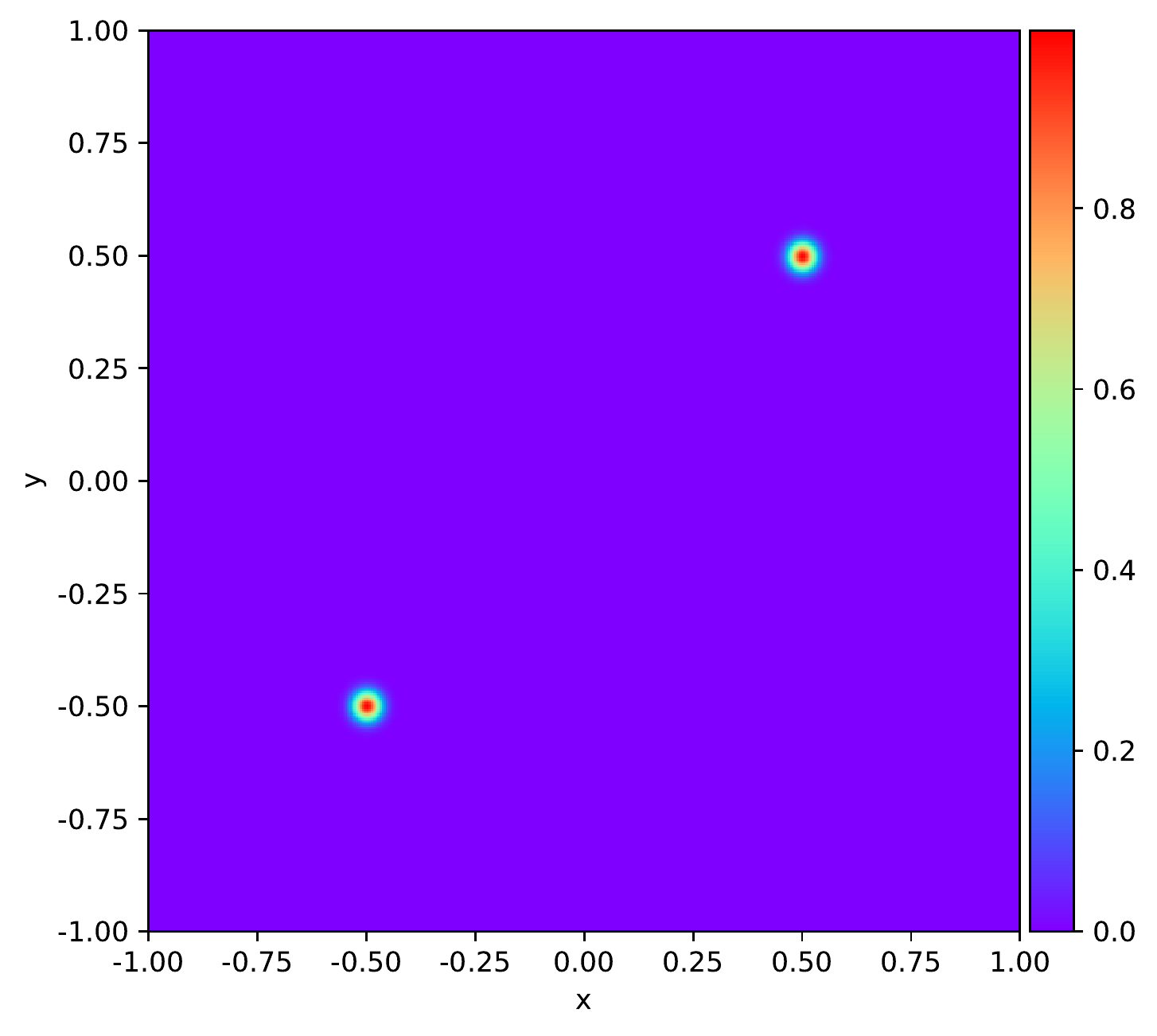}
(b)\includegraphics[width = 0.4\textwidth]{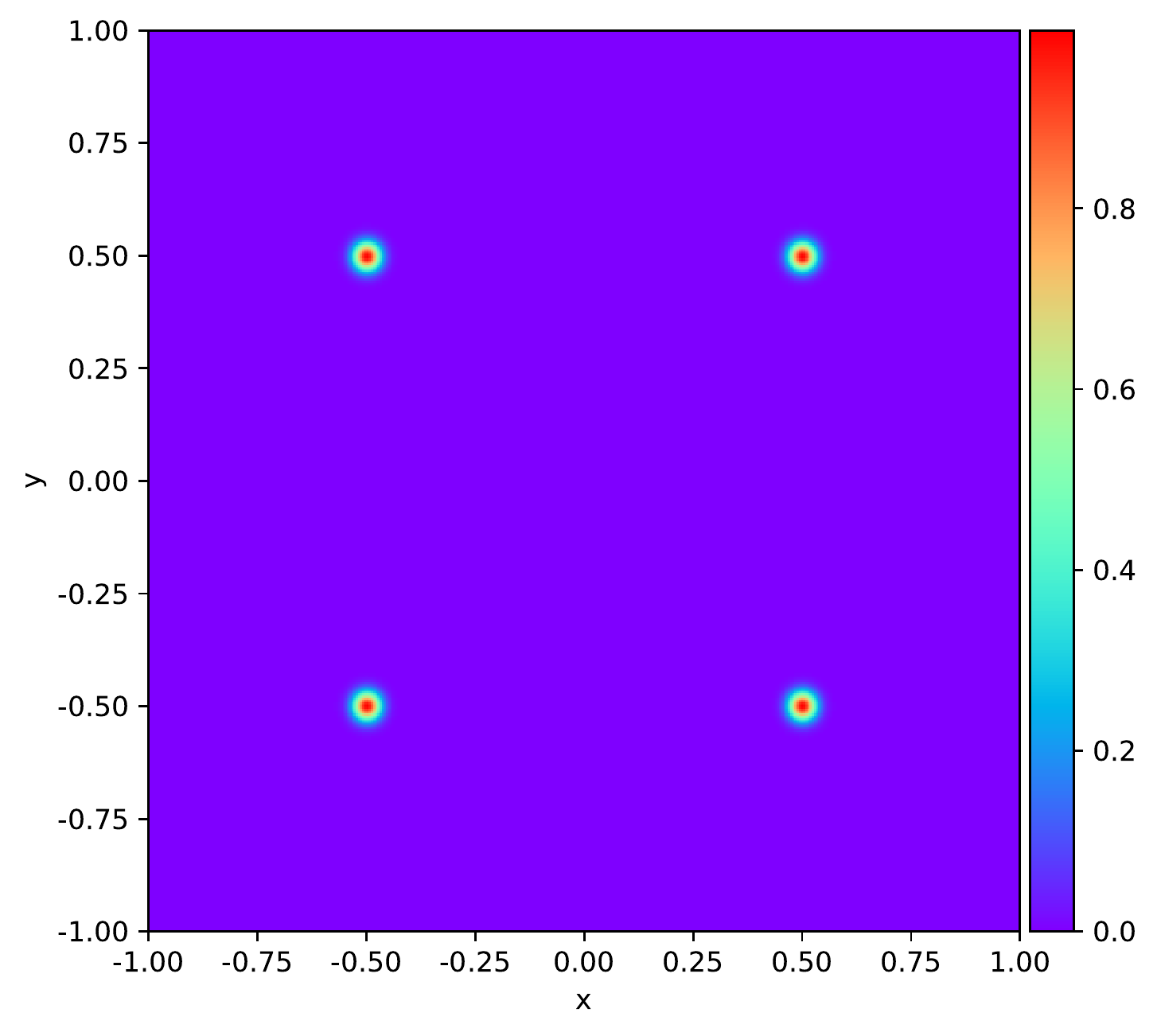}
\vspace{-0.4cm}
\caption{The contour plot of exact solution. (a) two peaks case; (b) four peaks case.}
\label{two_peak_true_solution}
\end{figure}

Here we consider the  following two-dimensional PDE with multiple peaks,
\begin{equation}
\begin{split}
-\bigtriangledown\cdot [u(x,y)\bigtriangledown(x^{2}+y^{2})] &+ \bigtriangledown^{2}u(x,y) = f(x,y),\quad (x,y)\in \Omega, \\
u(x,y) &= b(x,y),\quad (x,y)\in \partial \Omega,
\end{split}
\end{equation}
where $\Omega = [-1,1]\times[-1,1]$. The true solution is
\begin{equation}
u(x,y) = \sum_{i=1}^{k}e^{\left(-1000((x - x_{i})^{2} + (y-y_{i})^{2})\right)},
\end{equation}
which has several peaks at $(x_{i},y_{i}),i=0,\ldots,k$ and will decay to zero exponentially in other places. This problem has low regularity and is chosen here to test whether our algorithm can handle problems with multiple failure regions.

\subsubsection{Two peaks case}

\begin{figure}[htbp]
   \begin{center}
   \begin{overpic}[width = 0.3\textwidth]{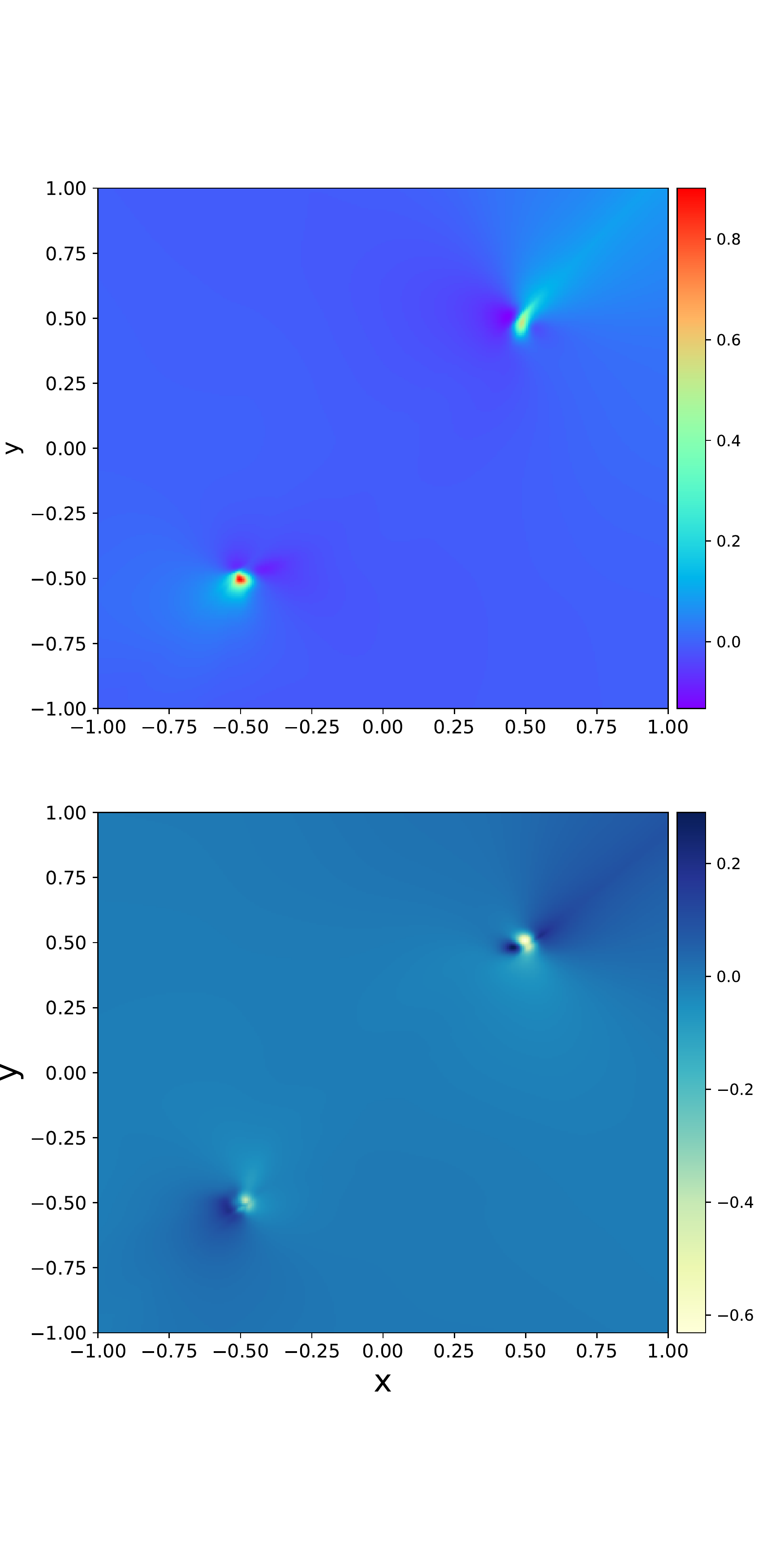}
      \put(15,90){\small MC-FIPINN}
   \end{overpic}
   \begin{overpic}[width = 0.3\textwidth]{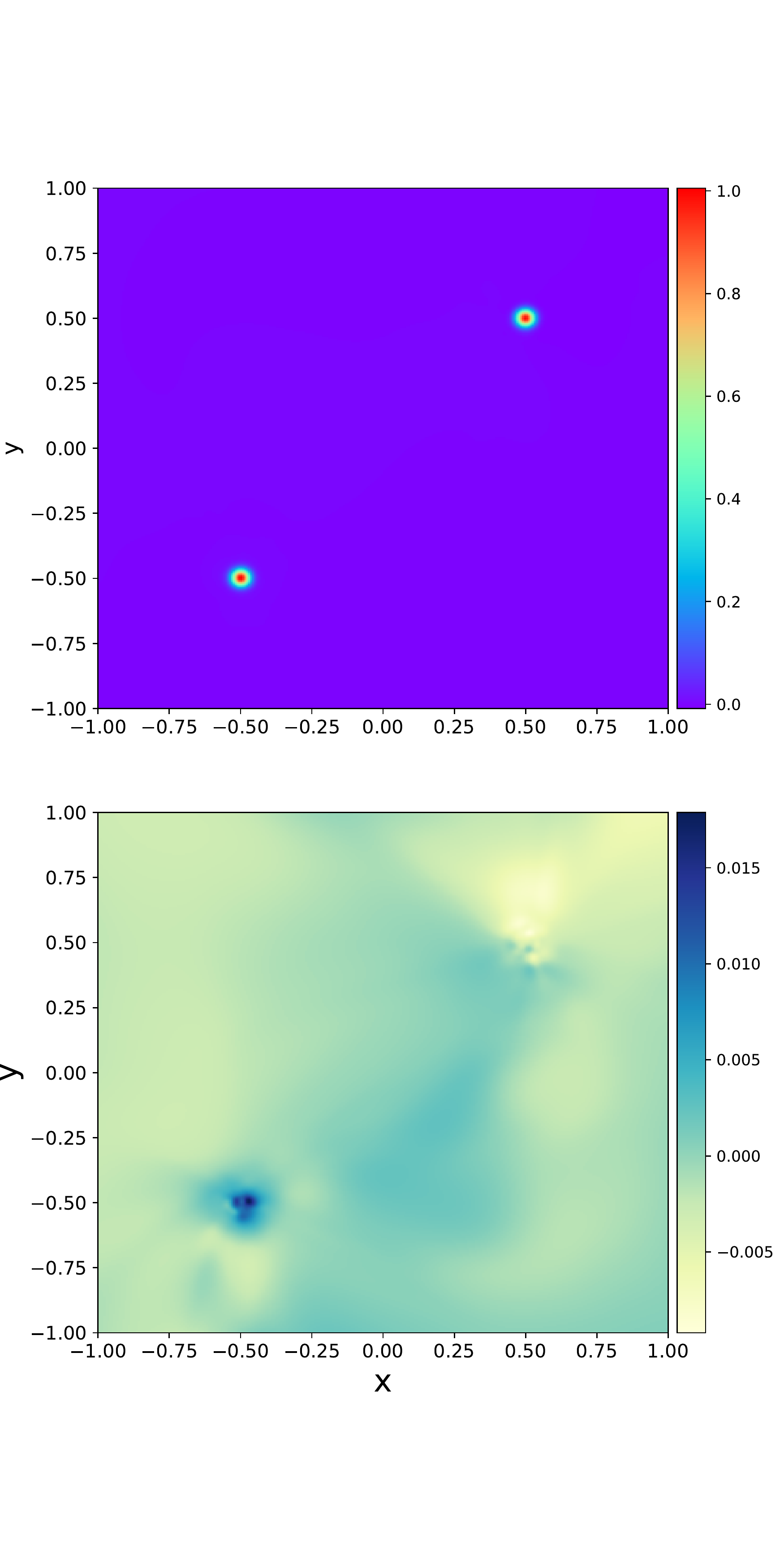}
      \put(15,90){\small R-FIPINN}
   \end{overpic}
   \begin{overpic}[width = 0.3\textwidth]{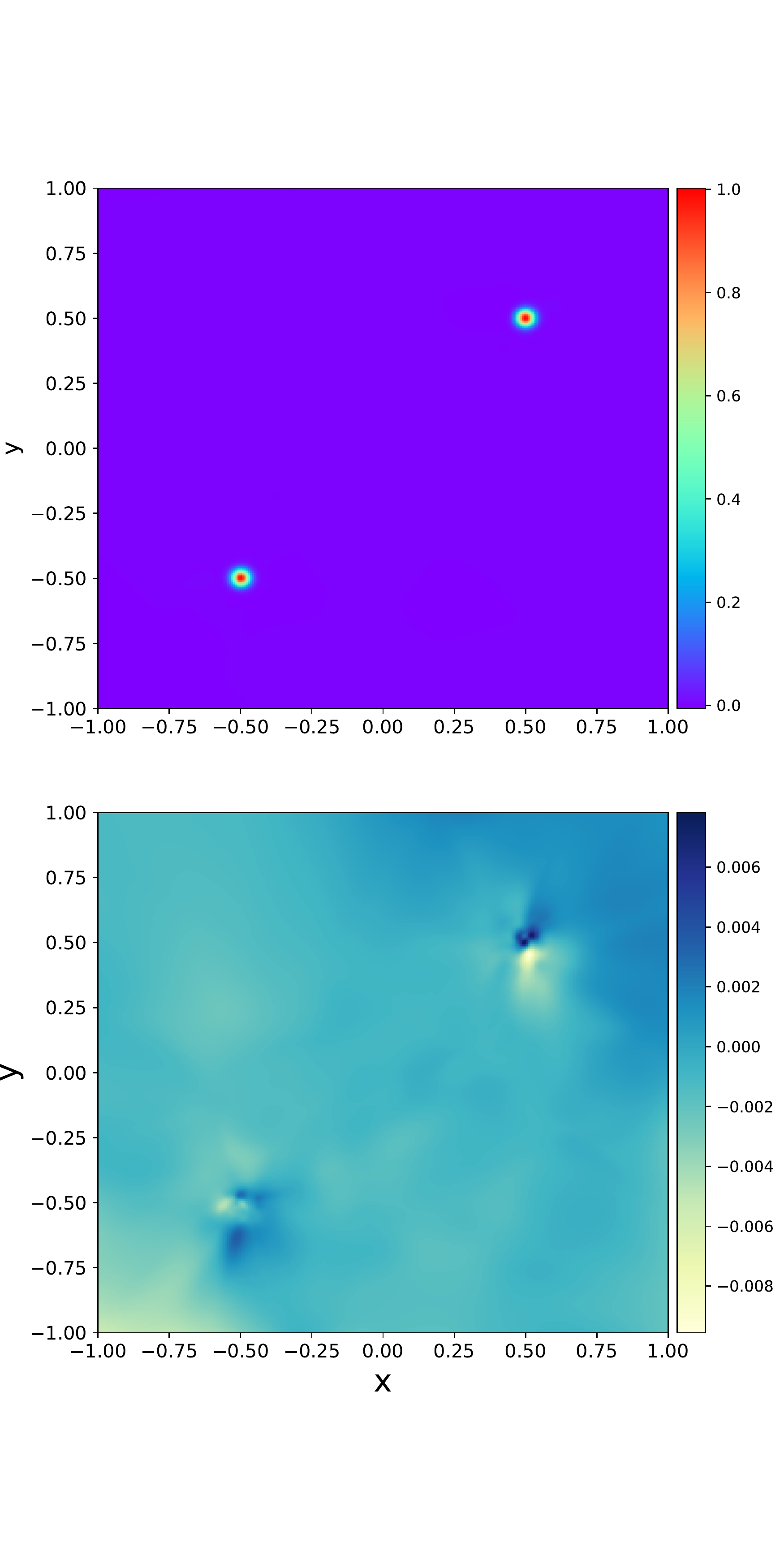}
      \put(15,90){\small G-FIPINN}
   \end{overpic}
   \end{center}
   \vspace{-1cm}
   \caption{The predicted solution (top) and the corresponding absolute error (bottom) obtained by different methods for two peaks case ($N_{c} = 1000$).}
   \label{two_peak_predictd_solution}
   \end{figure}

\begin{figure}[t]
\begin{center}
\begin{overpic}[width = 0.3\textwidth]{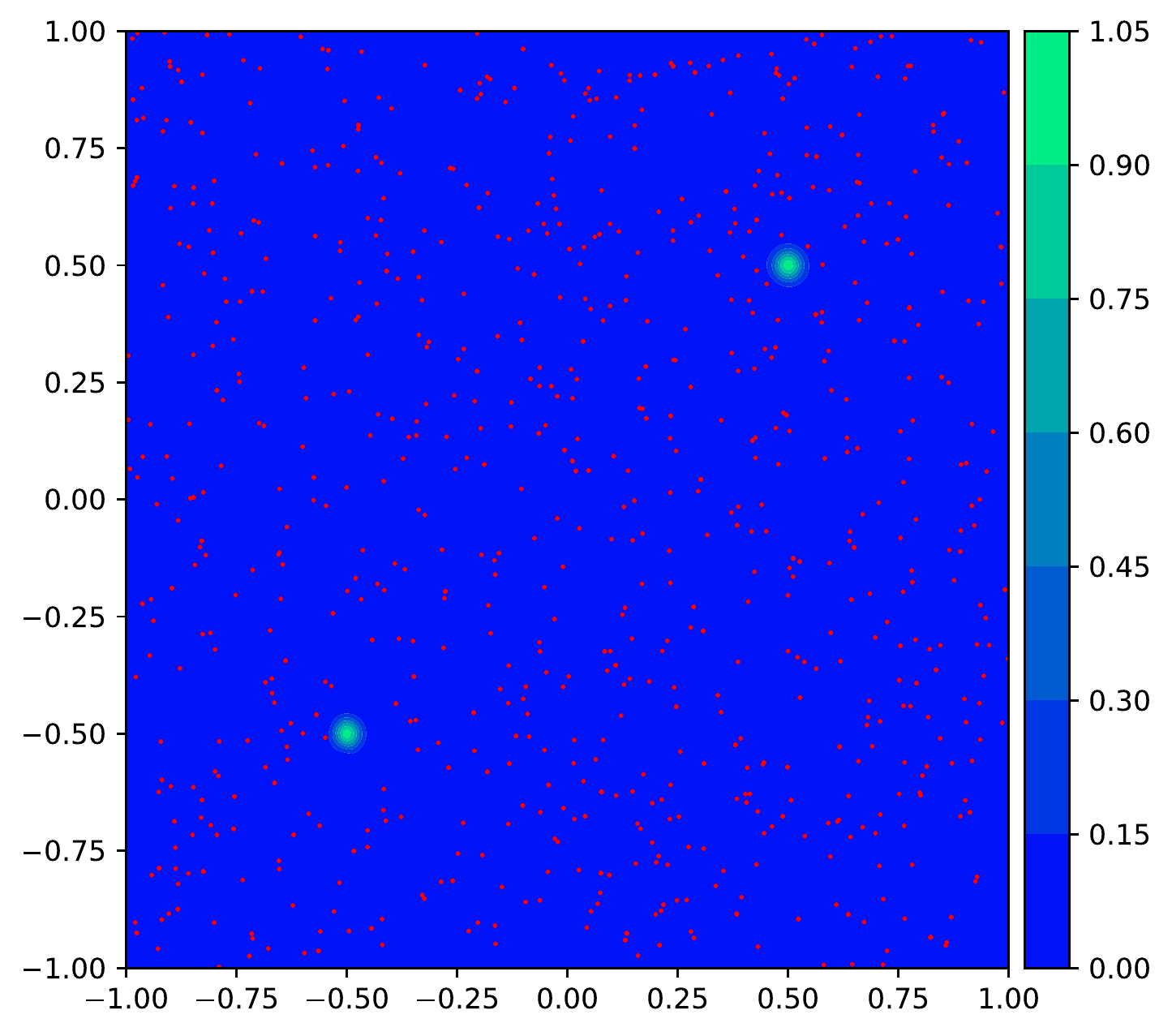}
   \put(35,90){\small MC-FIPINN}
\end{overpic}
\begin{overpic}[width = 0.3\textwidth]{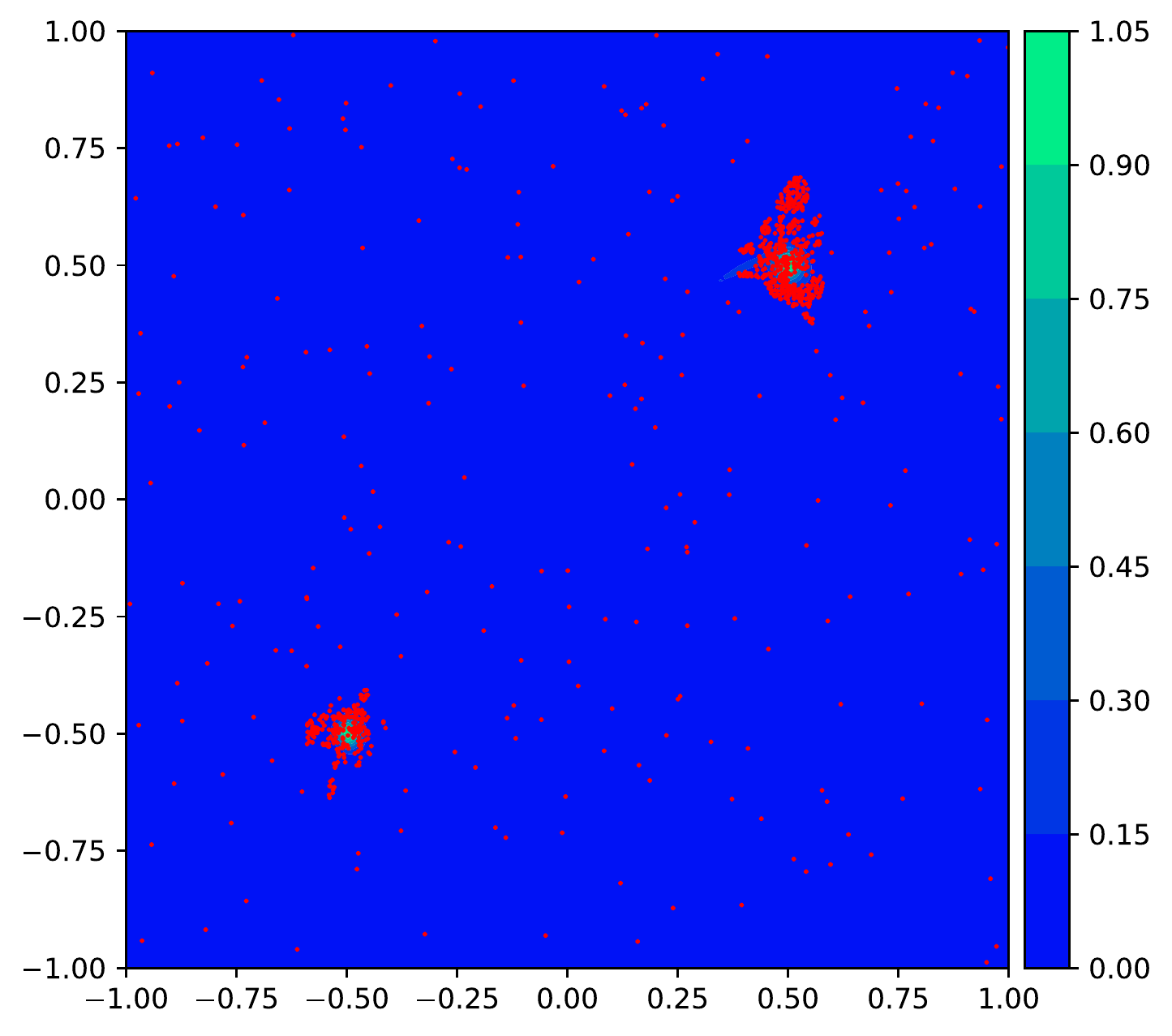}
   \put(35,90){\small R-FIPINN}
\end{overpic}
\begin{overpic}[width = 0.3\textwidth]{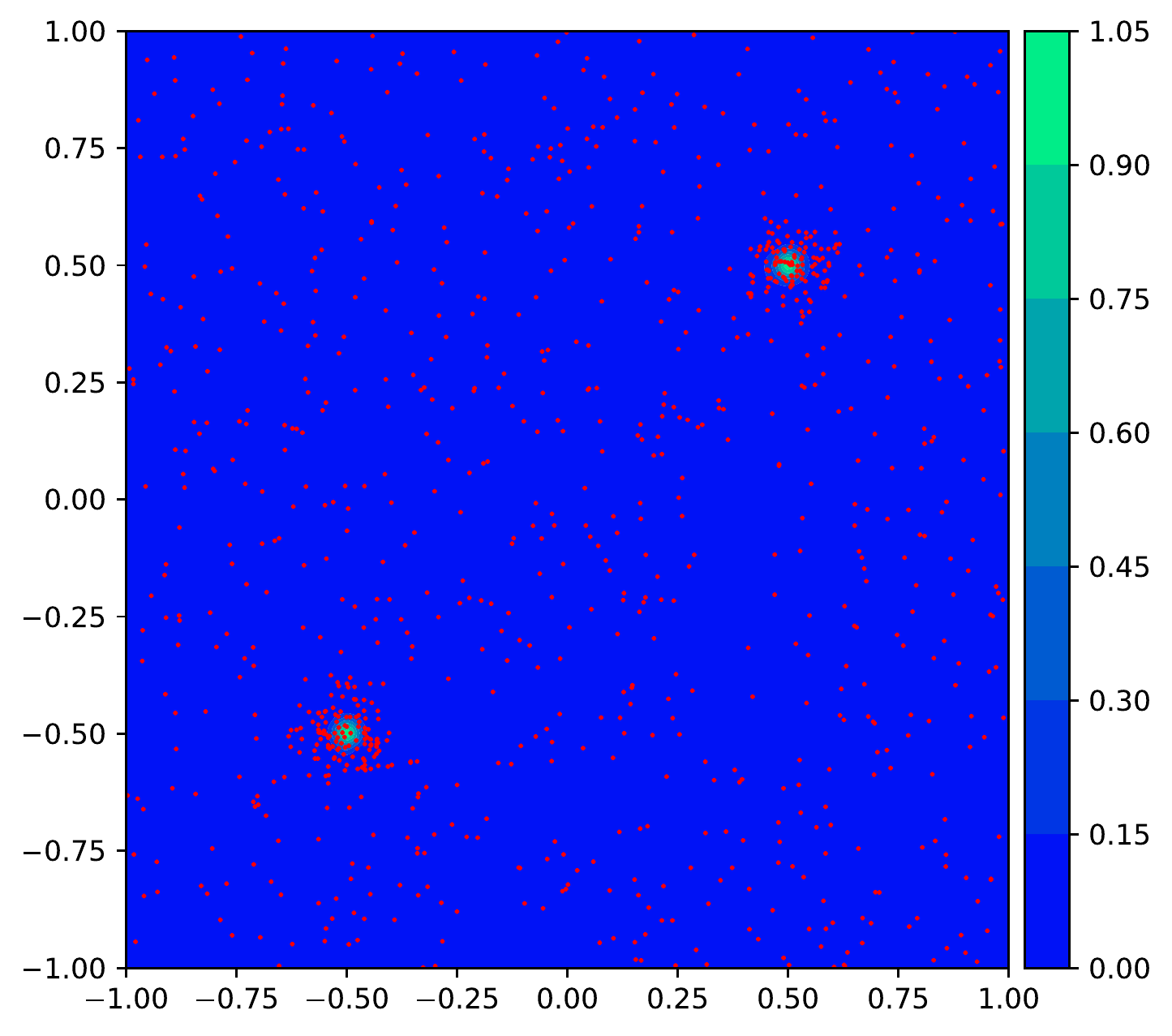}
   \put(35,90){\small G-FIPINN}
\end{overpic}
\end{center}
\vspace{-0.2cm}
\caption{Final distribution of collocation points obtained by using different methods ($N_{c}=1000$).}
\label{two_peak_samples}
\end{figure}

\begin{figure}[htbp]
\begin{center}
   \begin{overpic}[width = 0.3\textwidth]{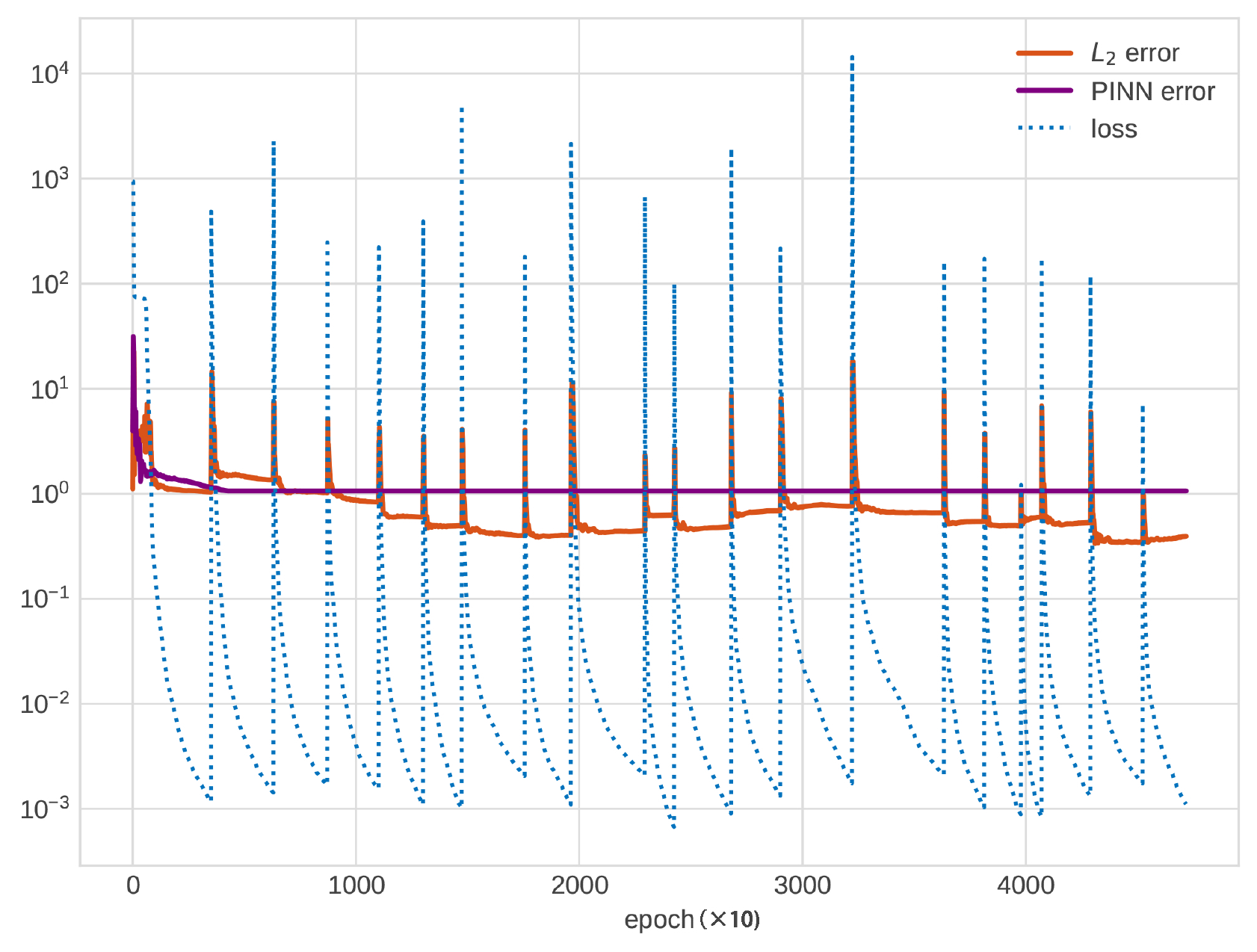}
      \put(35,78){\small MC-FIPINN}
   \end{overpic}
   \begin{overpic}[width = 0.3\textwidth]{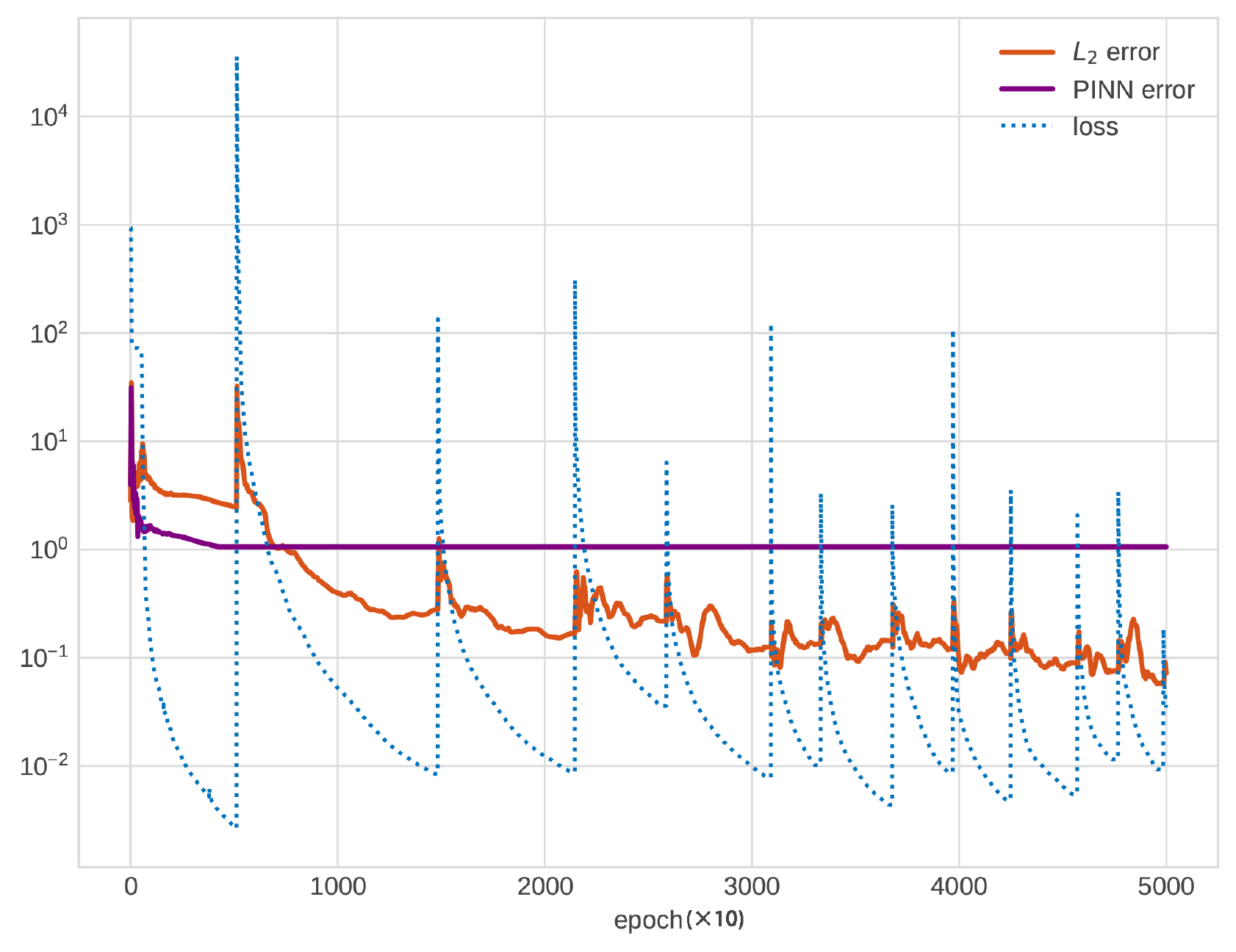}
      \put(35,78){\small R-FIPINN}
   \end{overpic}
   \begin{overpic}[width = 0.3\textwidth]{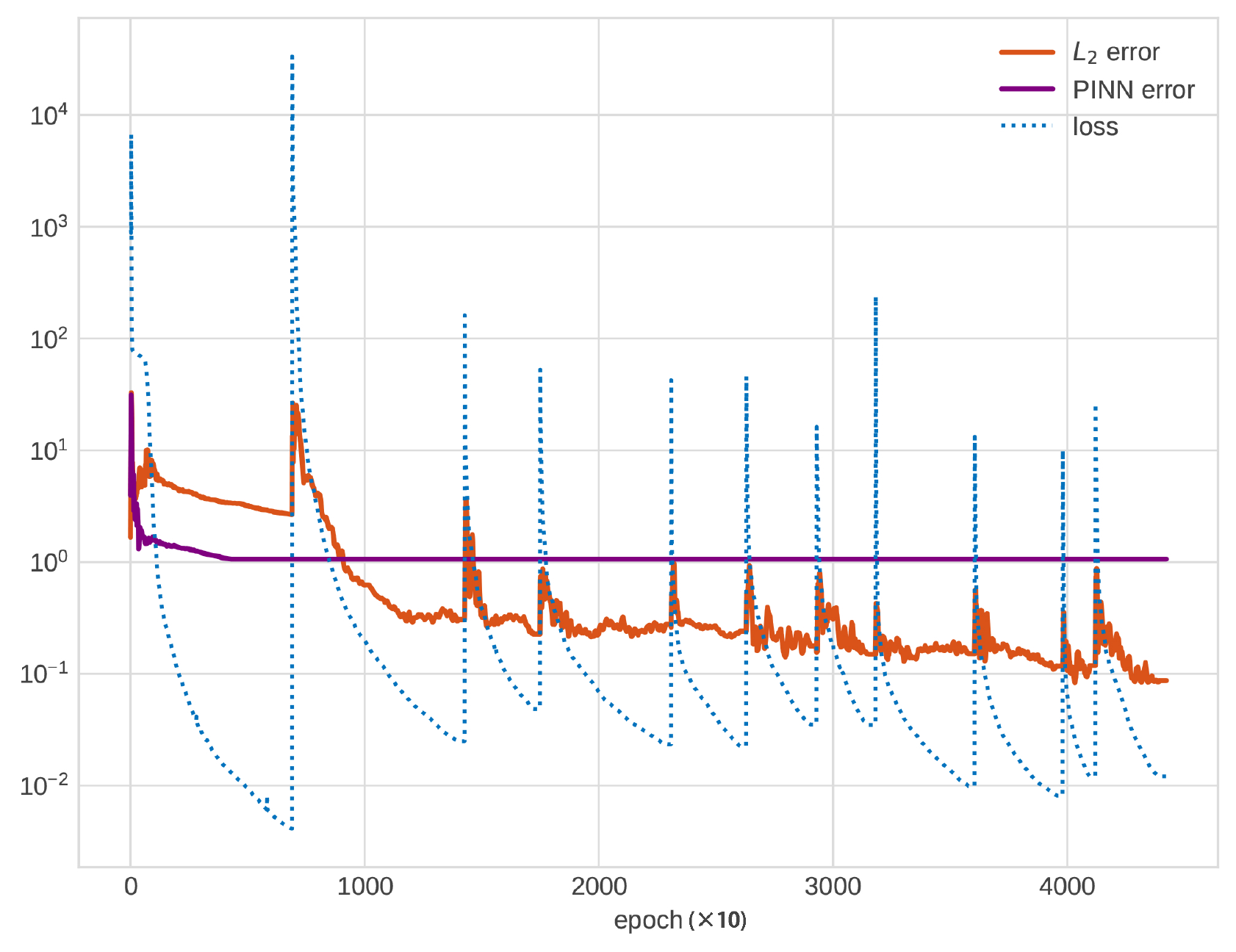}
      \put(35,78){\small G-FIPINN}
   \end{overpic}
\end{center}
\vspace{-0.2cm}
\caption{Relative error and the training loss compared to vanilla PINN  during the training process for two peaks problem ($N_{c} = 1000$).}
\label{two_peak_compared_error}
\end{figure}

\begin{figure}[htbp]
   \begin{center}
      \begin{overpic}[width = 0.45\textwidth]{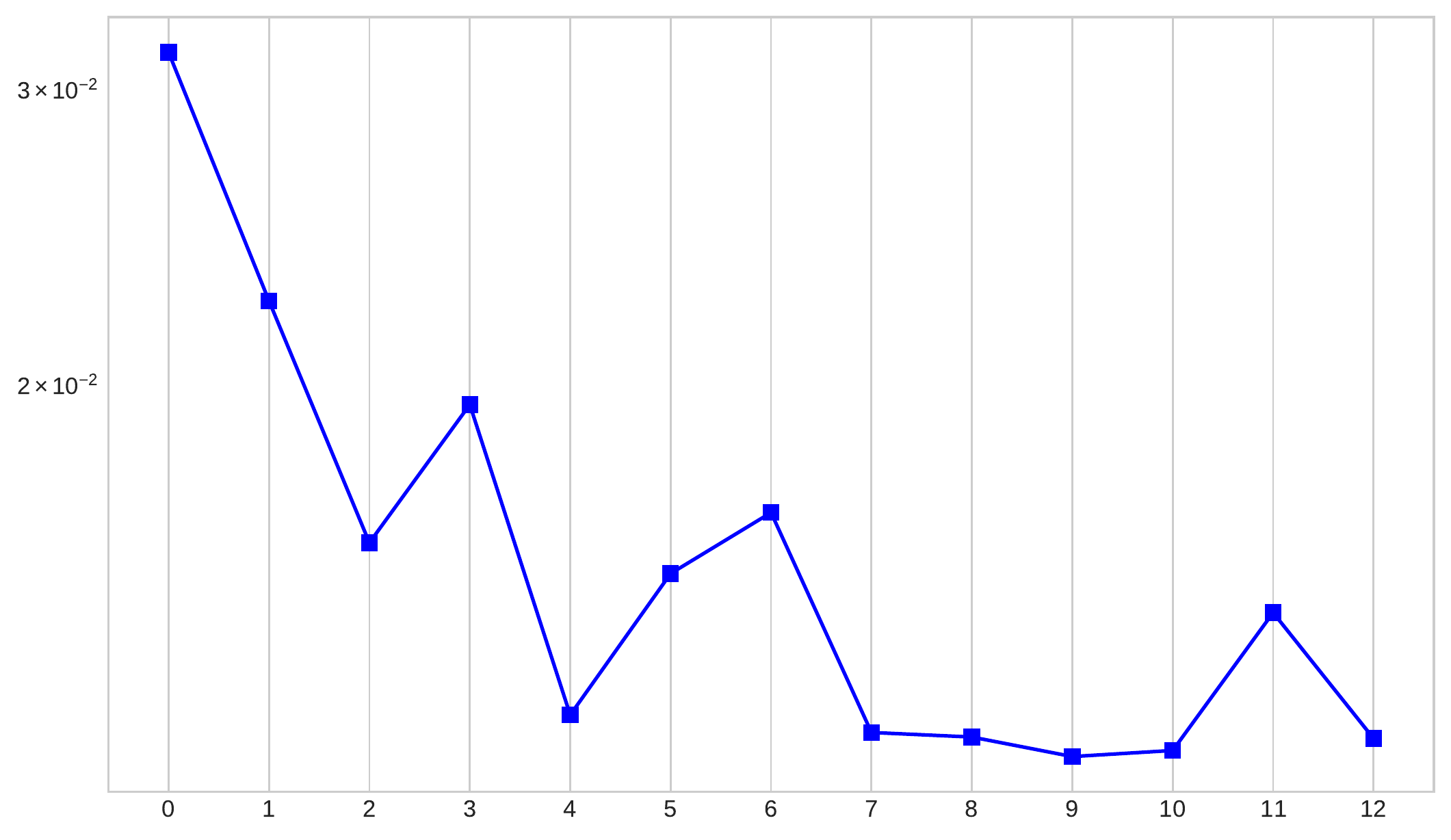}
         \put(35,60){\small R-FIPINN}
      \end{overpic}
      \begin{overpic}[width = 0.44\textwidth]{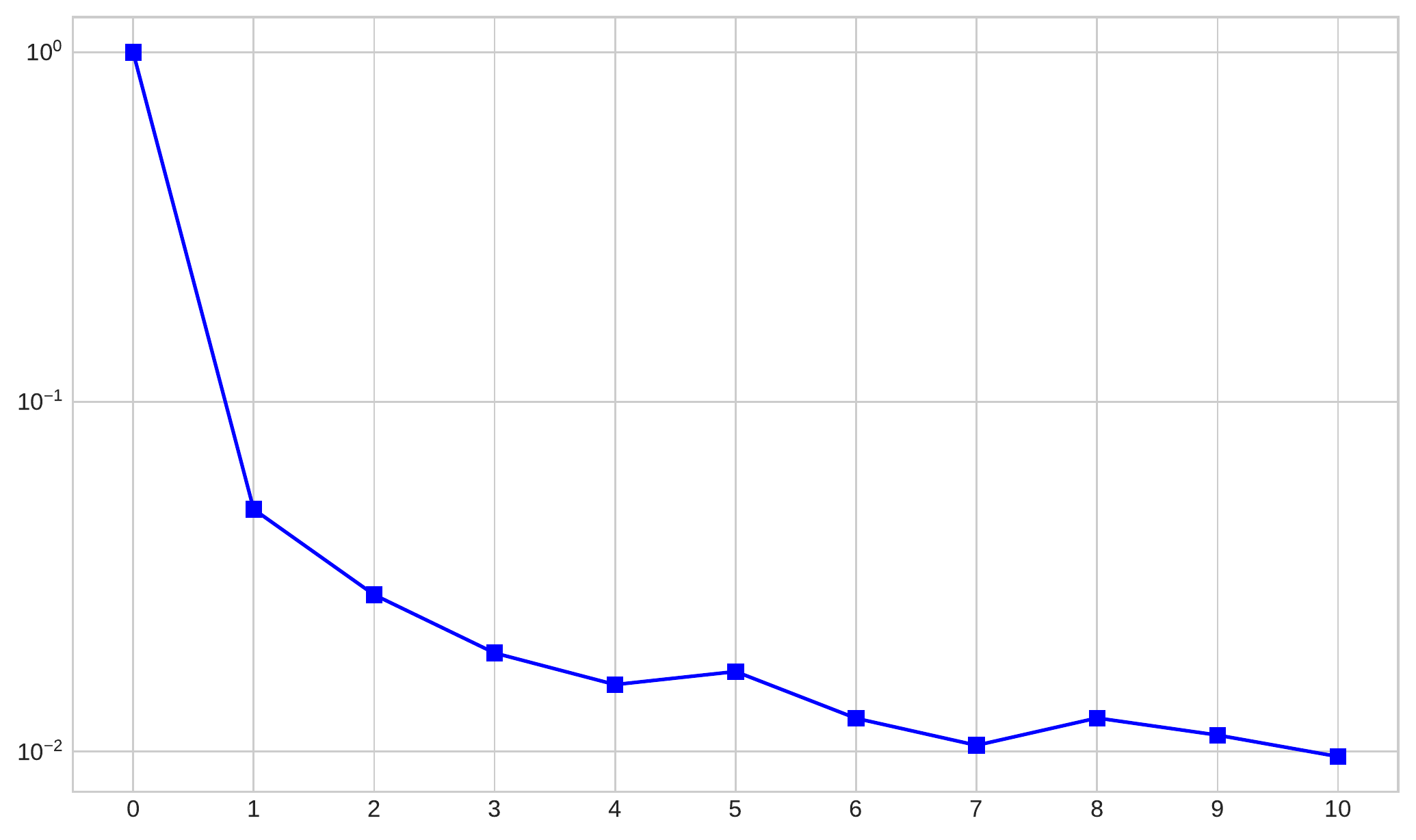}
         \put(35,60){\small G-FIPINN}
      \end{overpic}
   \end{center}
   \vspace{-0.2cm}
   \caption{Estimated failure probabilities over restarts  for two peaks problem ($N_{c} = 1000$).}
   \label{two_peak_failure_probability}
   \end{figure}

   \begin{figure}[htbp]
\centering
(a)\includegraphics[width = 0.4\textwidth]{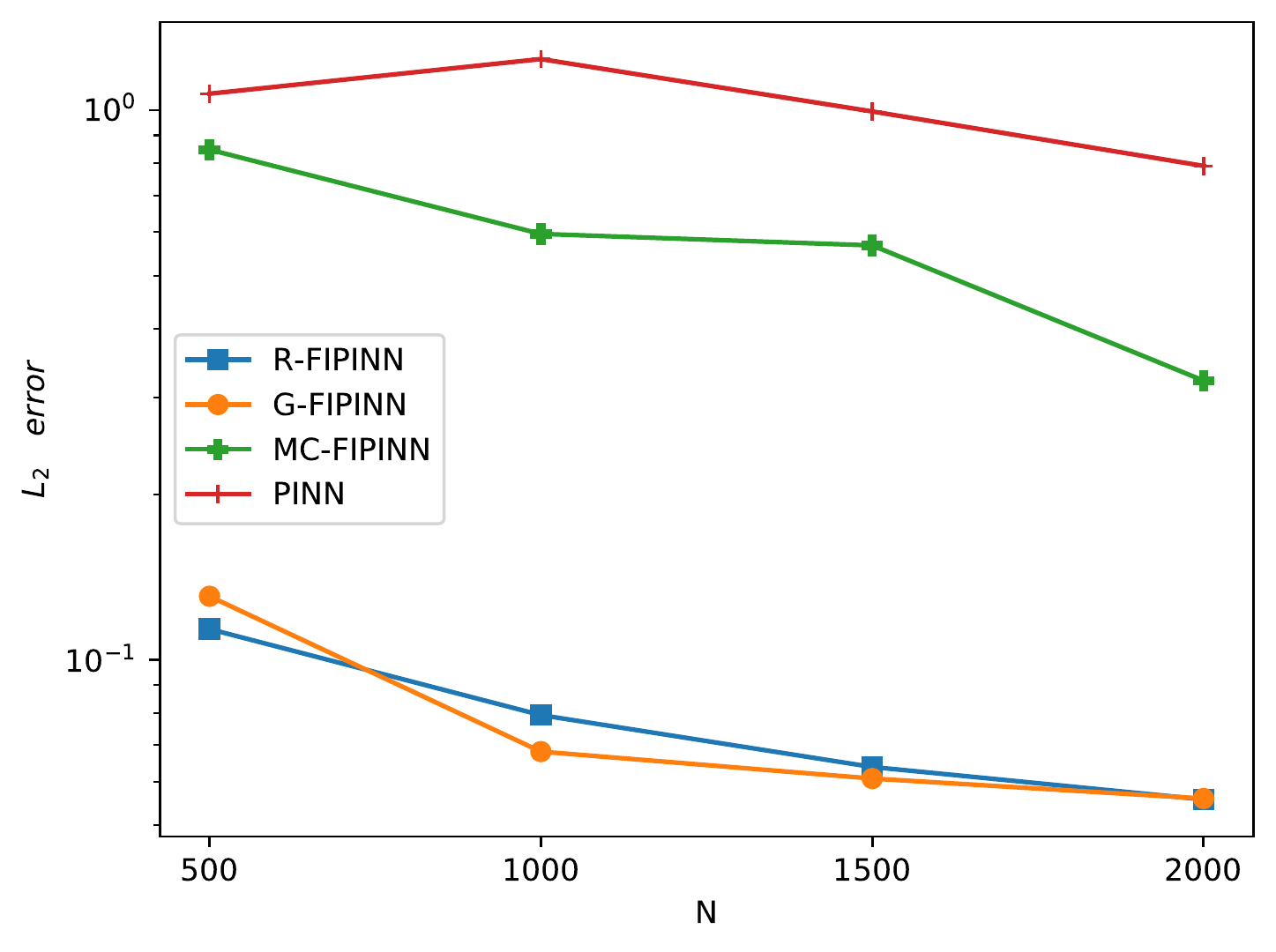}
(b)\includegraphics[width = 0.4\textwidth]{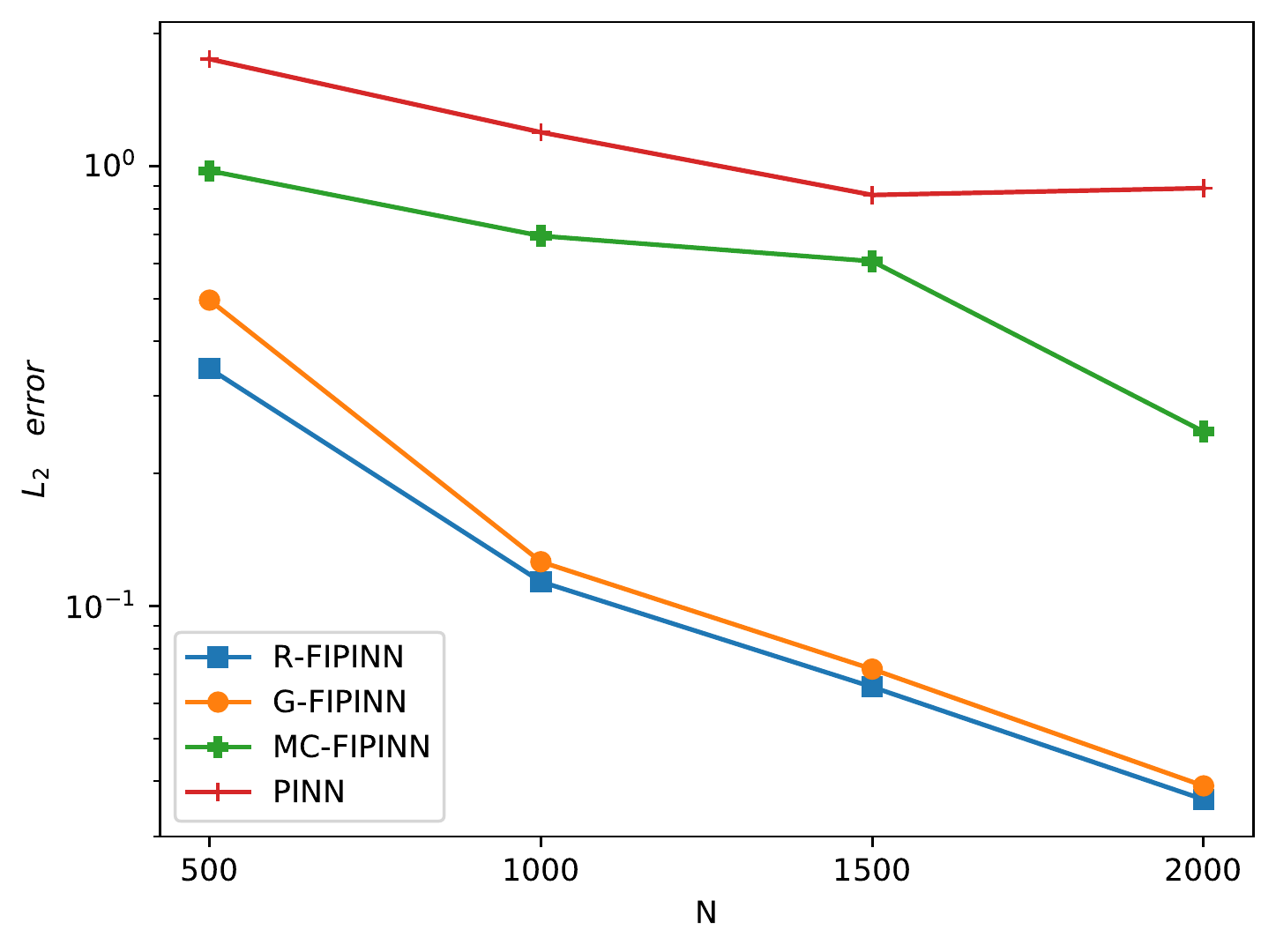}
\vspace{-0.2cm}
\caption{The $L_{2}$ error obtained by different methods when $N_{c}$ varies from 500 to 2000. (a)  two peaks case; (b) four peaks case.}
\label{two_peak_error}
\end{figure}

We consider the problem with two peaks, $(0.5,0.5)$ and  $(-0.5,-0.5)$, which is demonstrated in Fig. \ref{two_peak_true_solution}(a).

We first investigate the performance of the AFI-PINNs method using the relatively small number of collocation pints of $N_c =1000$.   In  Fig.\ref{two_peak_predictd_solution}, the predicted solution and the corresponding absolute error obtained using various methods are displayed. It is not surprising that even using this small number of collocation points,  the predicted solution obtained by our method differs indistinguishable from the true solution. This is due to the adaptive samples produced by the subset simulation focusing on the failure region during restarts, as shown in Fig.\ref{two_peak_samples}.  We can see that the distribution of the collocation points obtained by AFI-PINNs are much more concentrated around the two peaks $(0.5,0.5), (-0.5, -0.5)$. As a result, updating the collocation dataset improves the system reliability after the network has been retrained.

To better understand how the annealing restart works, we plot the predicted error and the training loss in Fig.\ref{two_peak_compared_error} to make a comparison with the vanilla PINN predicted error. We can see that the predicted error achieved by our method decreases much faster than the vanilla PINN, despite the fact that there is a sudden upward jump with each restart.  While the predicted error achieved by the MC-FIPINN decreases more slowly than our method. Also, as shown in Fig.\ref{two_peak_compared_error},  annealing restart can not only accelerate the convergence but also stabilize the training. This is because the new collocation dataset still contains new uniform samples generated from the prior, which can help balance the different areas of the domain.  As a result, the loss will suddenly increase and then decrease quickly with a smaller  error with each restart. Even after many updates to the training dataset, the  error obtained by MC-FIPINN sampling barely moves.  The estimated failure probability in Fig.\ref{two_peak_failure_probability}. It is clear that the failure probability gradually decrease during restarts and our training can stop early when it is smaller than the tolerance. This phenomenon verifies the effectiveness of our novel framework.

We then repeat the experiment, but this time we change the number of collocation points, $N_c$, from 500 to 2000, and present the relative error in Fig. \ref{two_peak_error}(a). We notice a continuous pattern in which AFI-PINNs consistently achieve up to an order of magnitude better error than vanilla PINNs and MC-FIPINN.  Also  note that the results obtained by the R-FIPINN and G-FIPINN are nearly identical in  this case. Furthermore, we can see that our technique can still produce a respectable performance with fewer collocation points.

\subsubsection{Four peaks case}

\begin{figure}[htbp]
      \begin{center}
         \begin{overpic}[width = 0.3\textwidth]{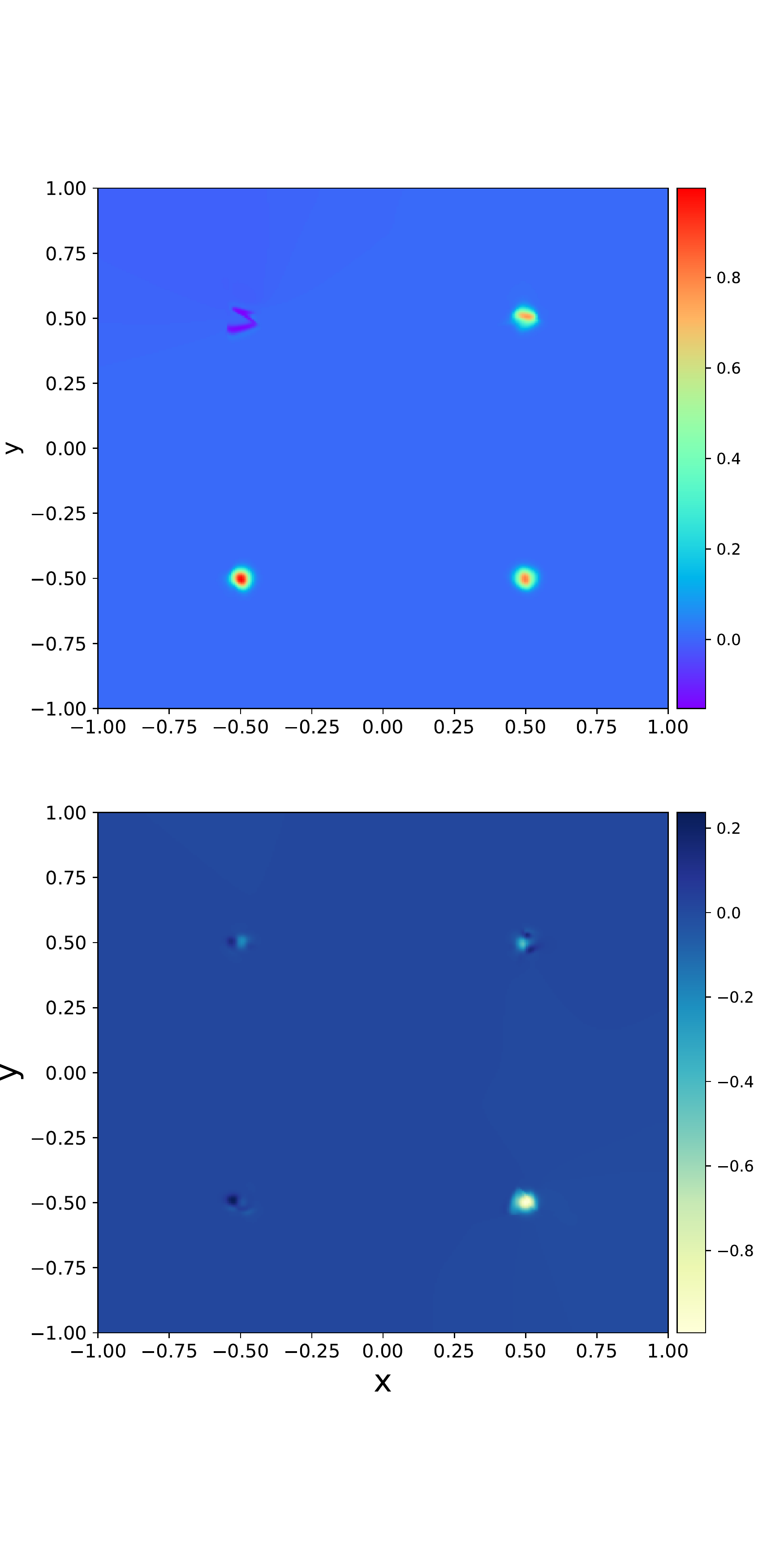}
            \put(15,90){\small MC-FIPINN}
         \end{overpic}
         \begin{overpic}[width = 0.3\textwidth]{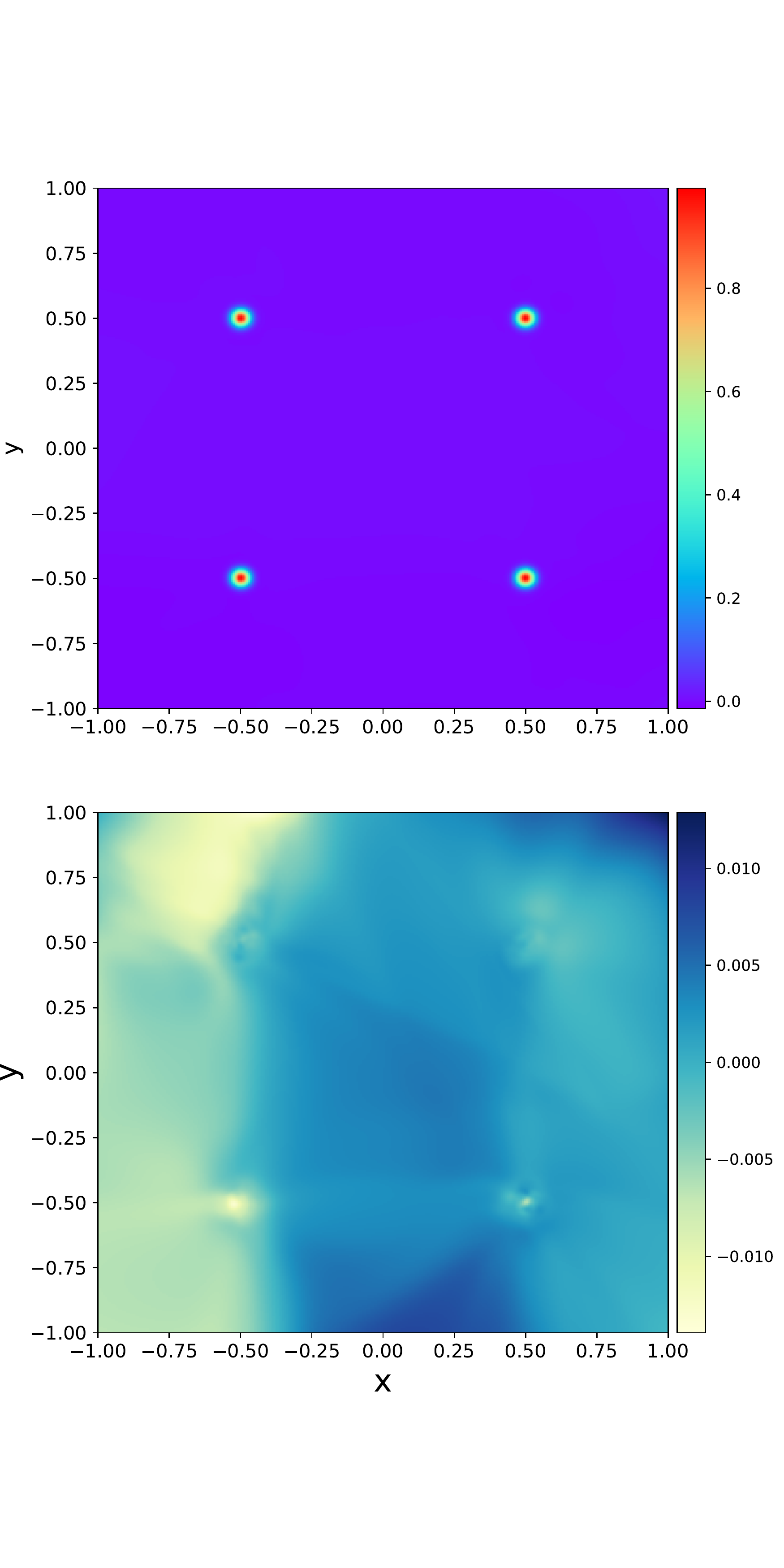}
            \put(15,90){\small R-FIPINN}
         \end{overpic}
         \begin{overpic}[width = 0.3\textwidth]{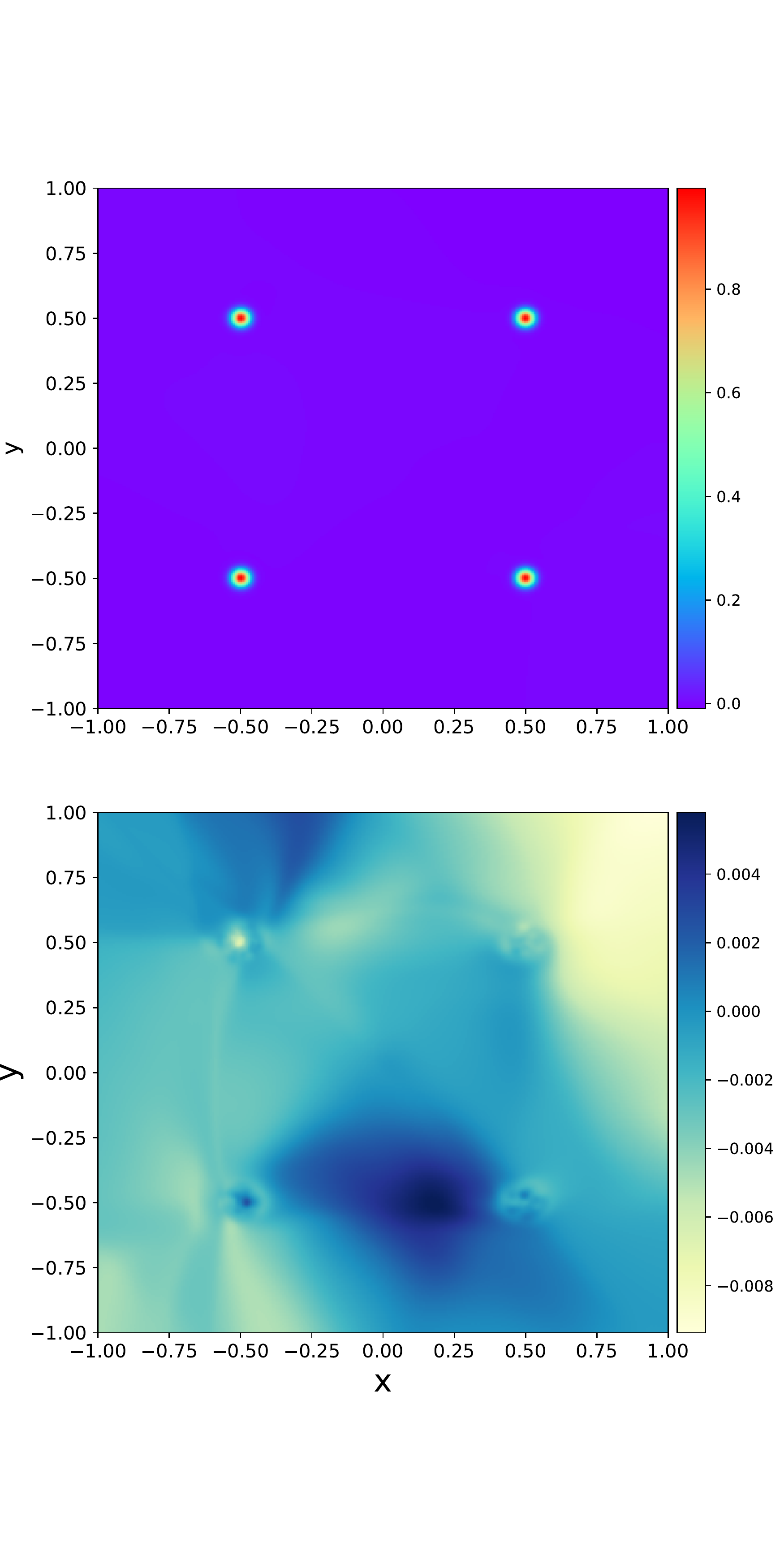}
            \put(15,90){\small G-FIPINN}
         \end{overpic}
      \end{center}
      \vspace{-.2cm}
      \caption{The predicted solution (top) and the corresponding absolute error (bottom) obtained by different methods for four peaks case ($N_{c} = 2000$).}
      \label{four_peak_predictd_solution}
      \end{figure}

\begin{figure}[htbp]
      \begin{center}
         \begin{overpic}[width = 0.3\textwidth]{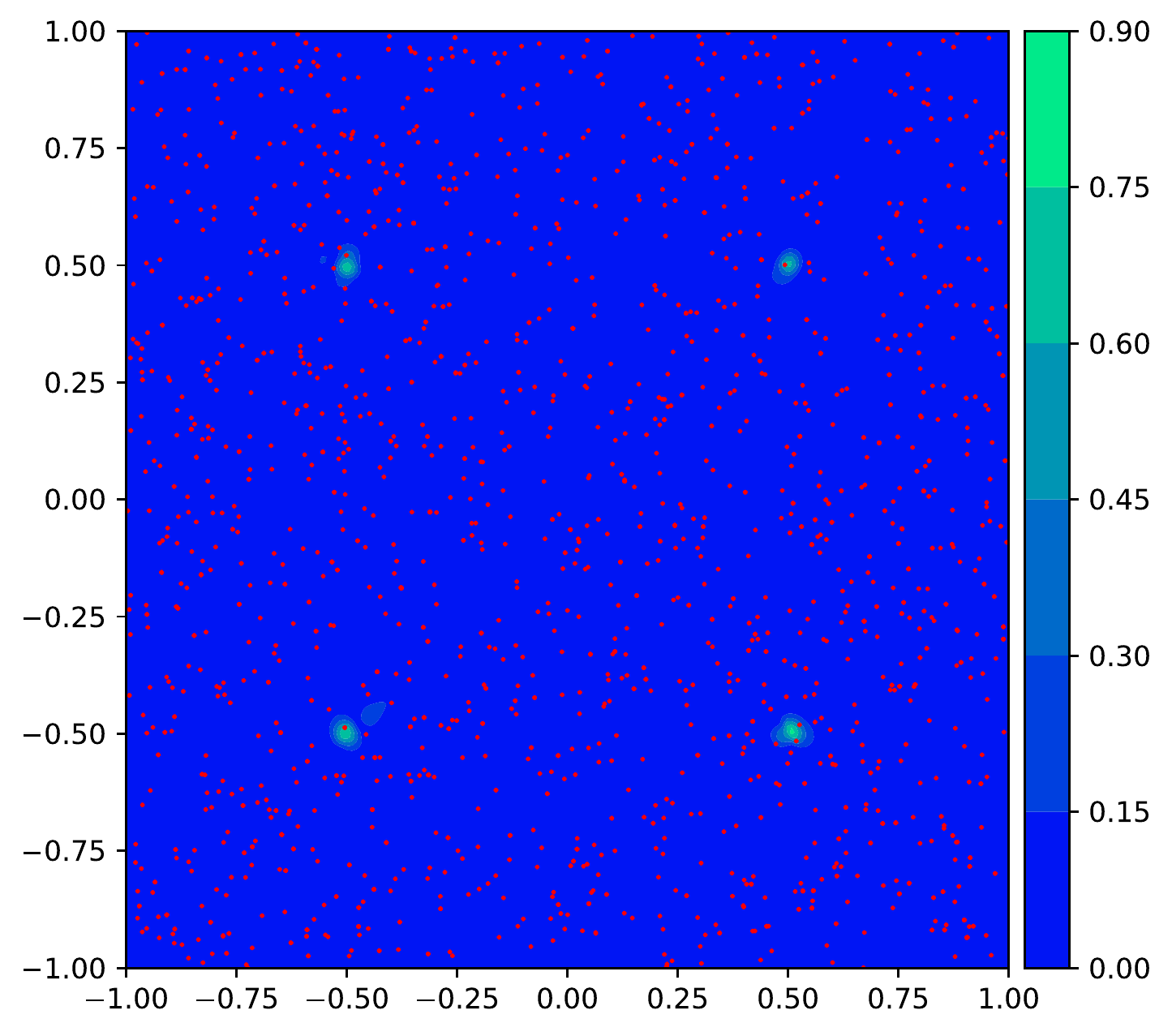}
            \put(35,90){\small MC-FIPINN}
         \end{overpic}
         \begin{overpic}[width = 0.3\textwidth]{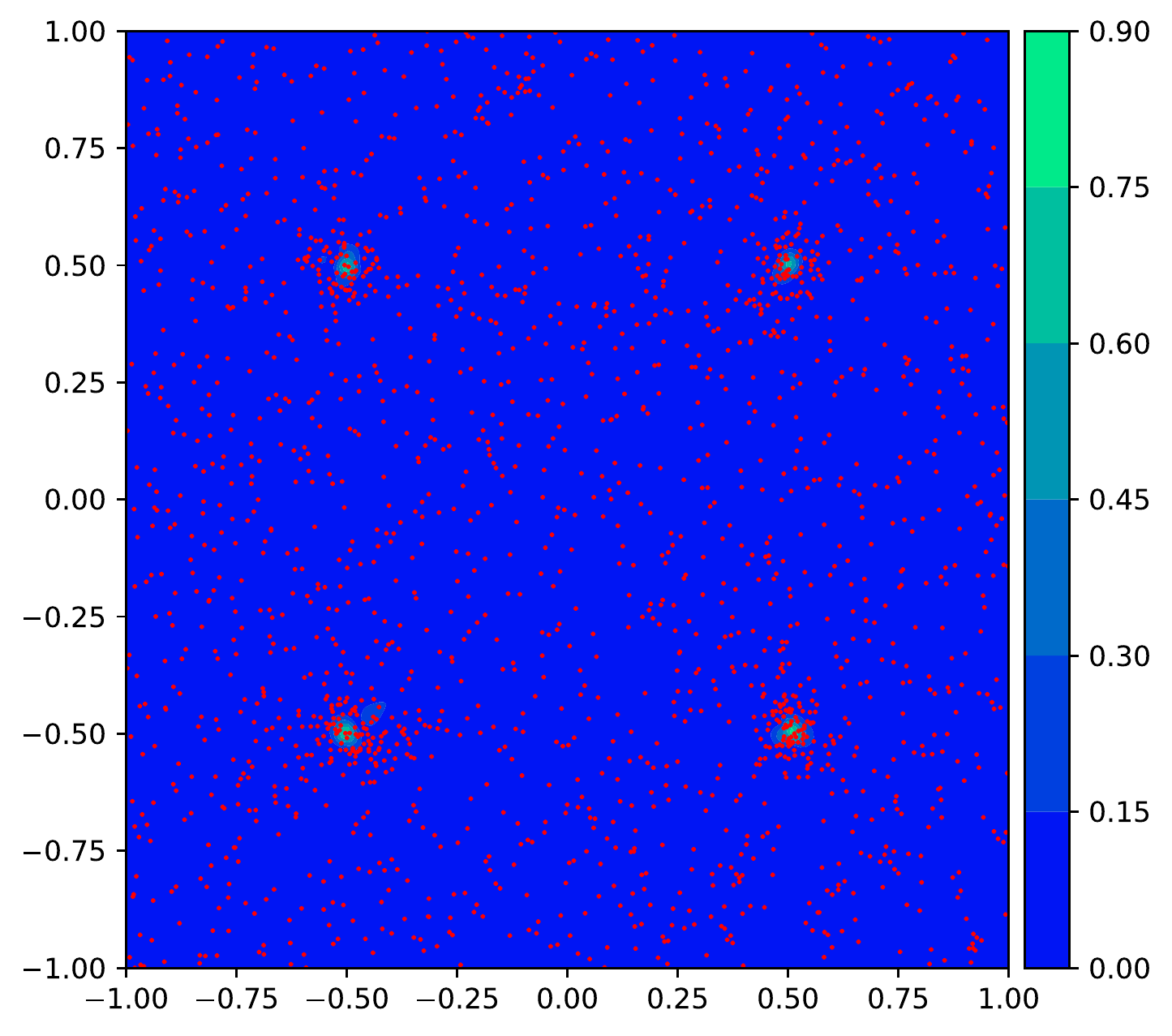}
            \put(35,90){\small R-FIPINN}
         \end{overpic}
         \begin{overpic}[width = 0.3\textwidth]{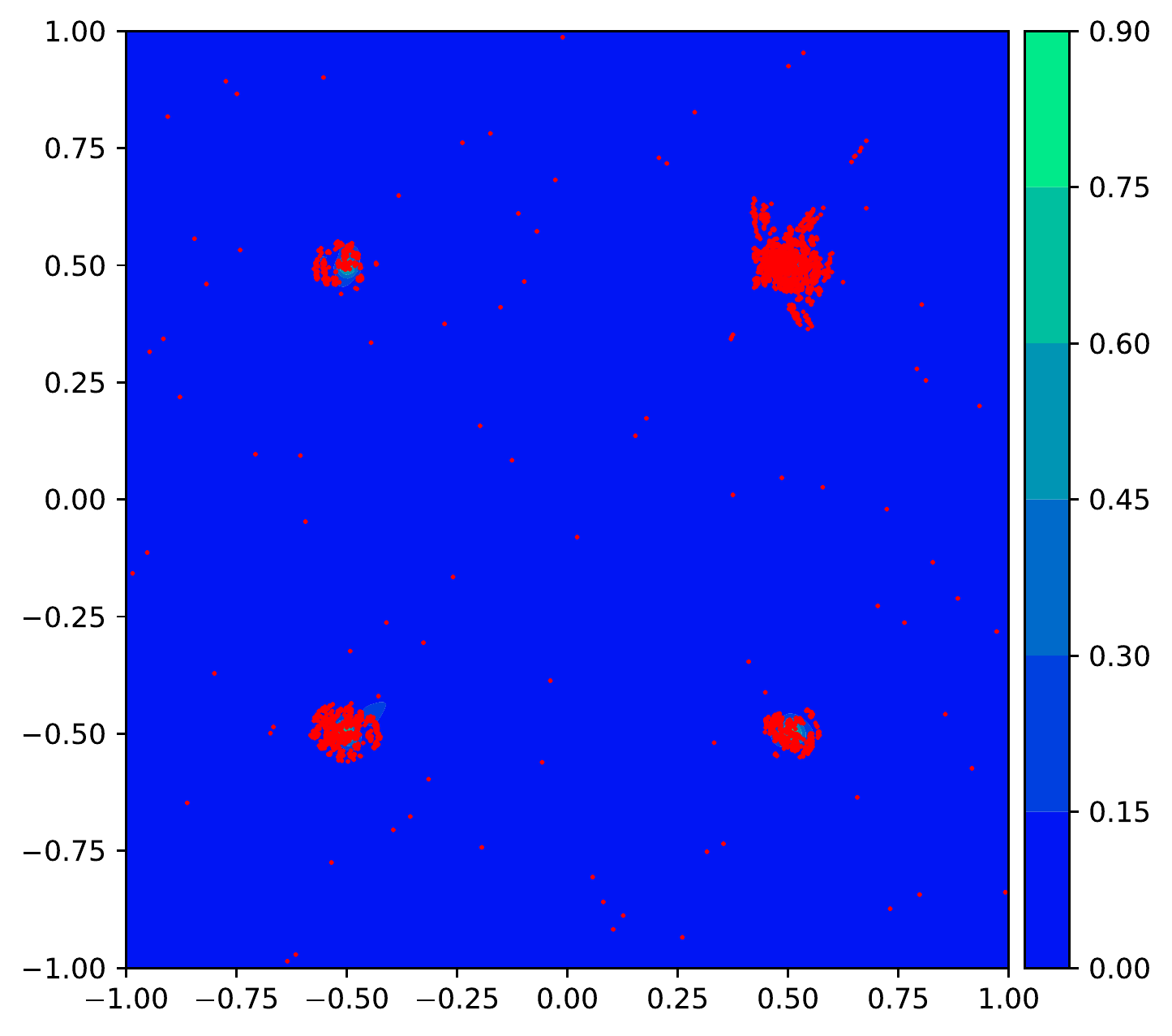}
            \put(35,90){\small G-FIPINN}
         \end{overpic}
      \end{center}
      \vspace{-0.2cm}
      \caption{Final distribution of collocation points obtained by using different methods ($N_{c}=2000$).}
      \label{four_peak_samples}
      \end{figure}

\begin{figure}[htbp]
         \begin{center}
            \begin{overpic}[width = 0.305\textwidth]{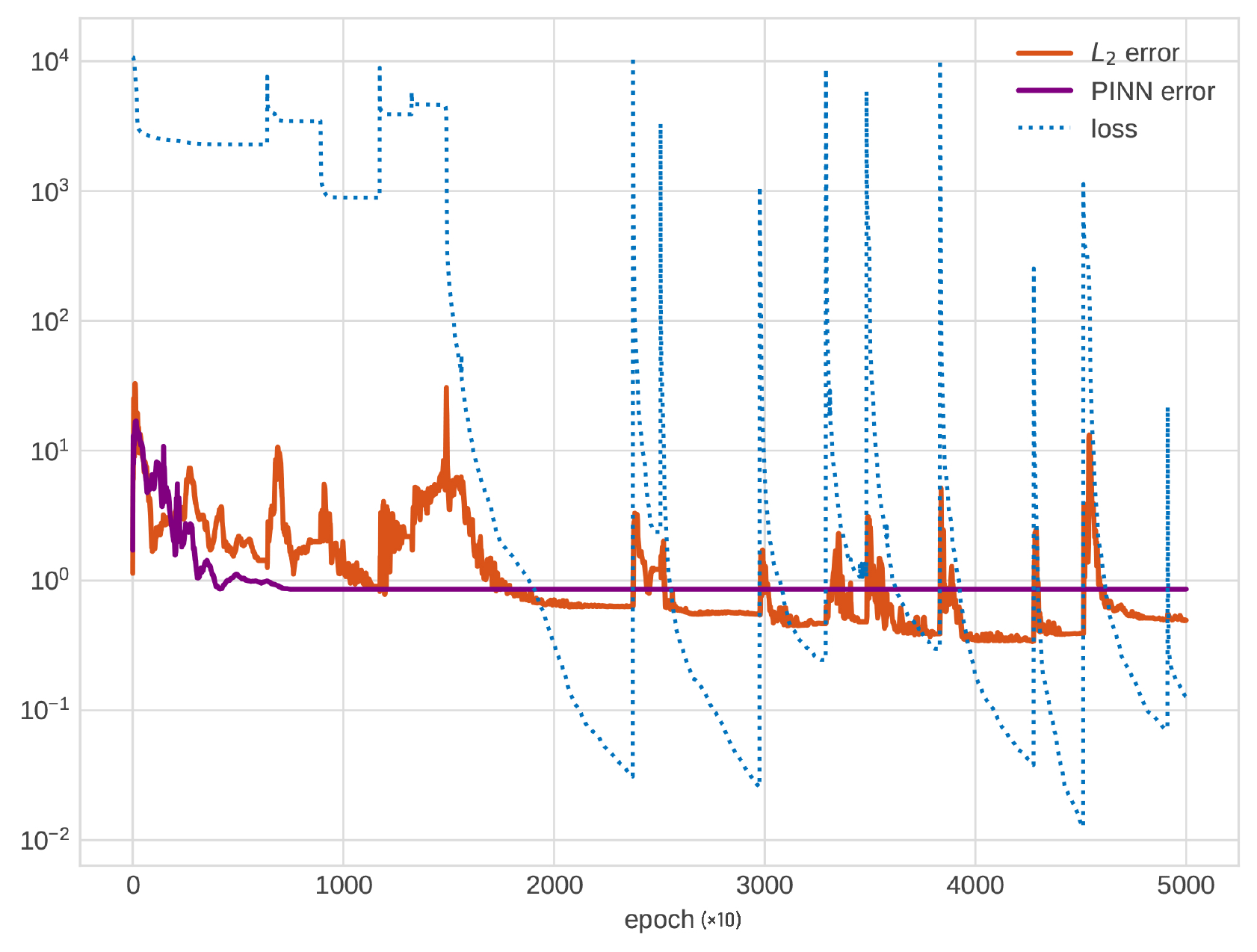}
               \put(35,78){\small MC-FIPINN}
            \end{overpic}
            \begin{overpic}[width = 0.3\textwidth]{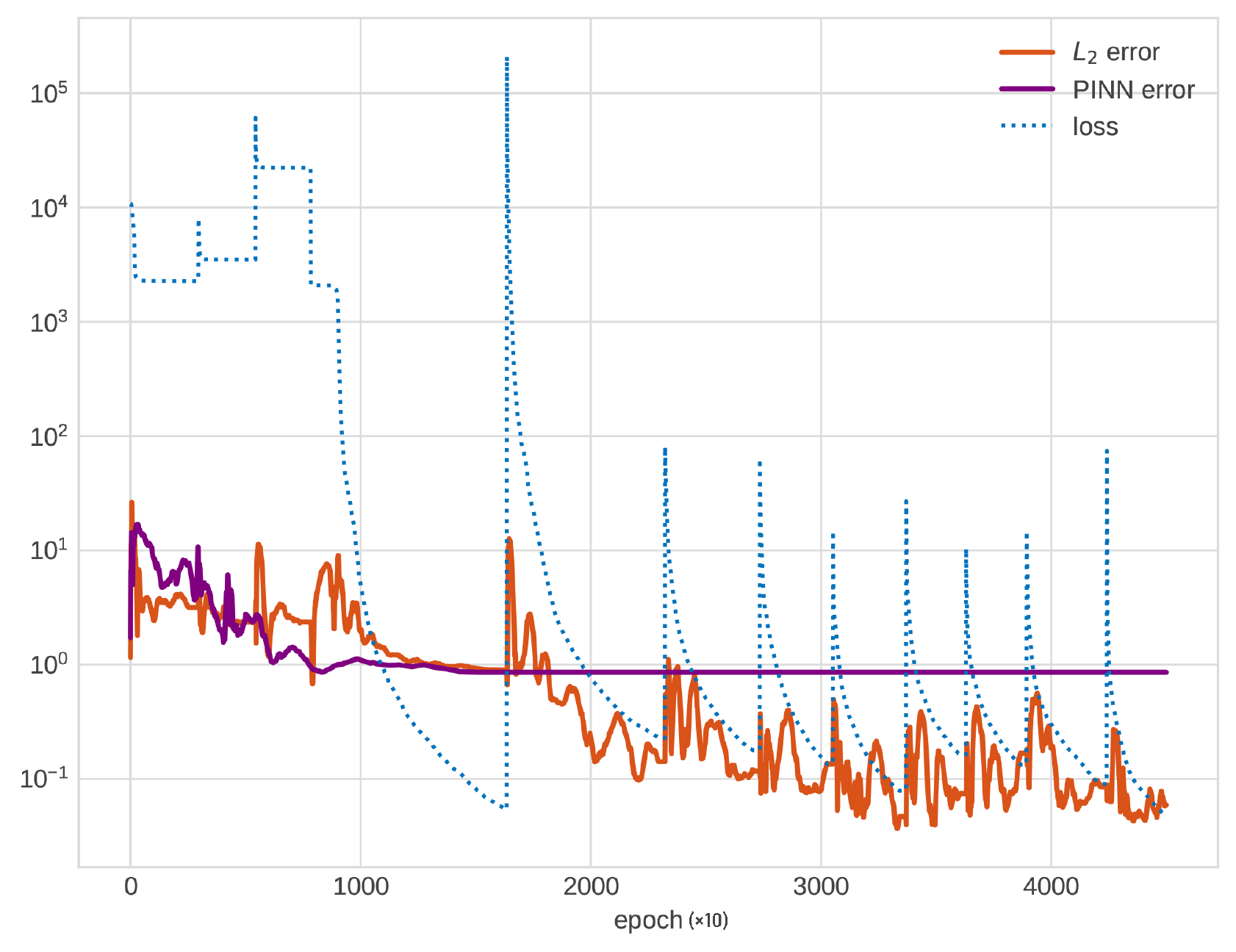}
               \put(35,78){\small R-FIPINN}
            \end{overpic}
            \begin{overpic}[width = 0.3\textwidth]{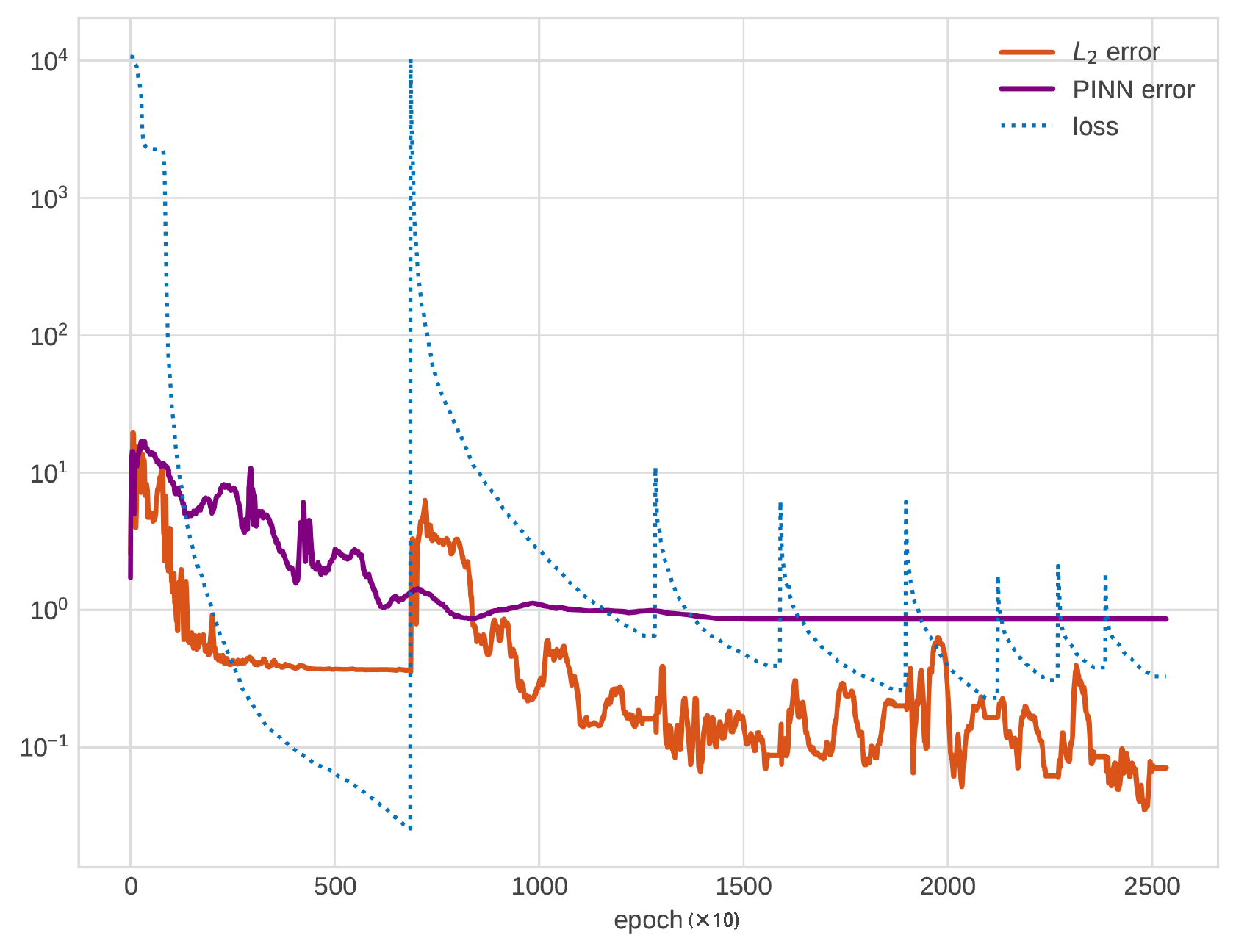}
               \put(35,78){\small G-FIPINN}
            \end{overpic}
         \end{center}
         \vspace{-0.2cm}
         \caption{Relative error and the training loss during the training process for four peaks problem ($N_{c} = 2000$).}
         \label{four_peak_compared_error}
         \end{figure}

\begin{figure}[htbp]
         \begin{center}
            \begin{overpic}[width = 0.40\textwidth]{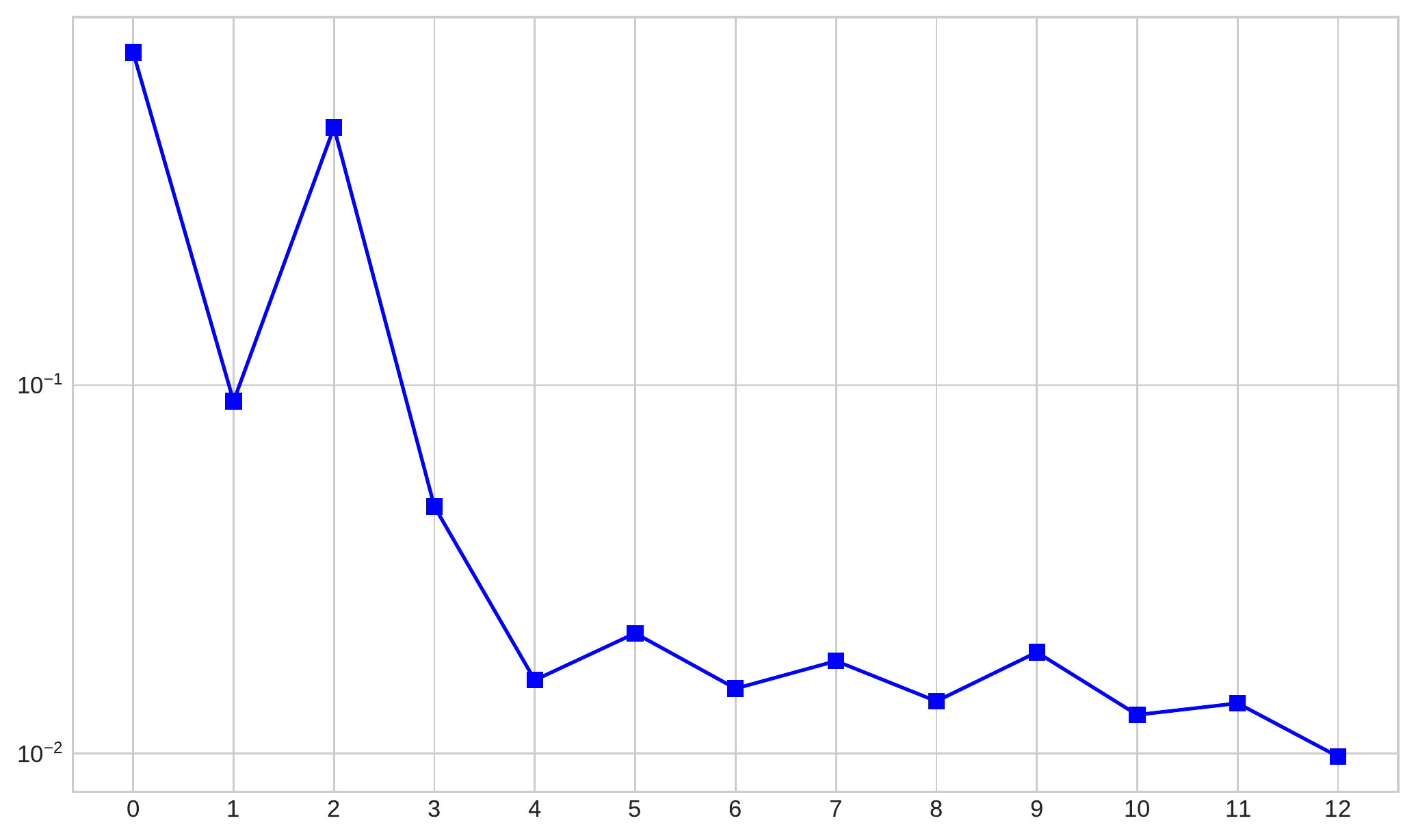}
               \put(35,60){\small R-FIPINN}
            \end{overpic}
            \begin{overpic}[width = 0.41\textwidth]{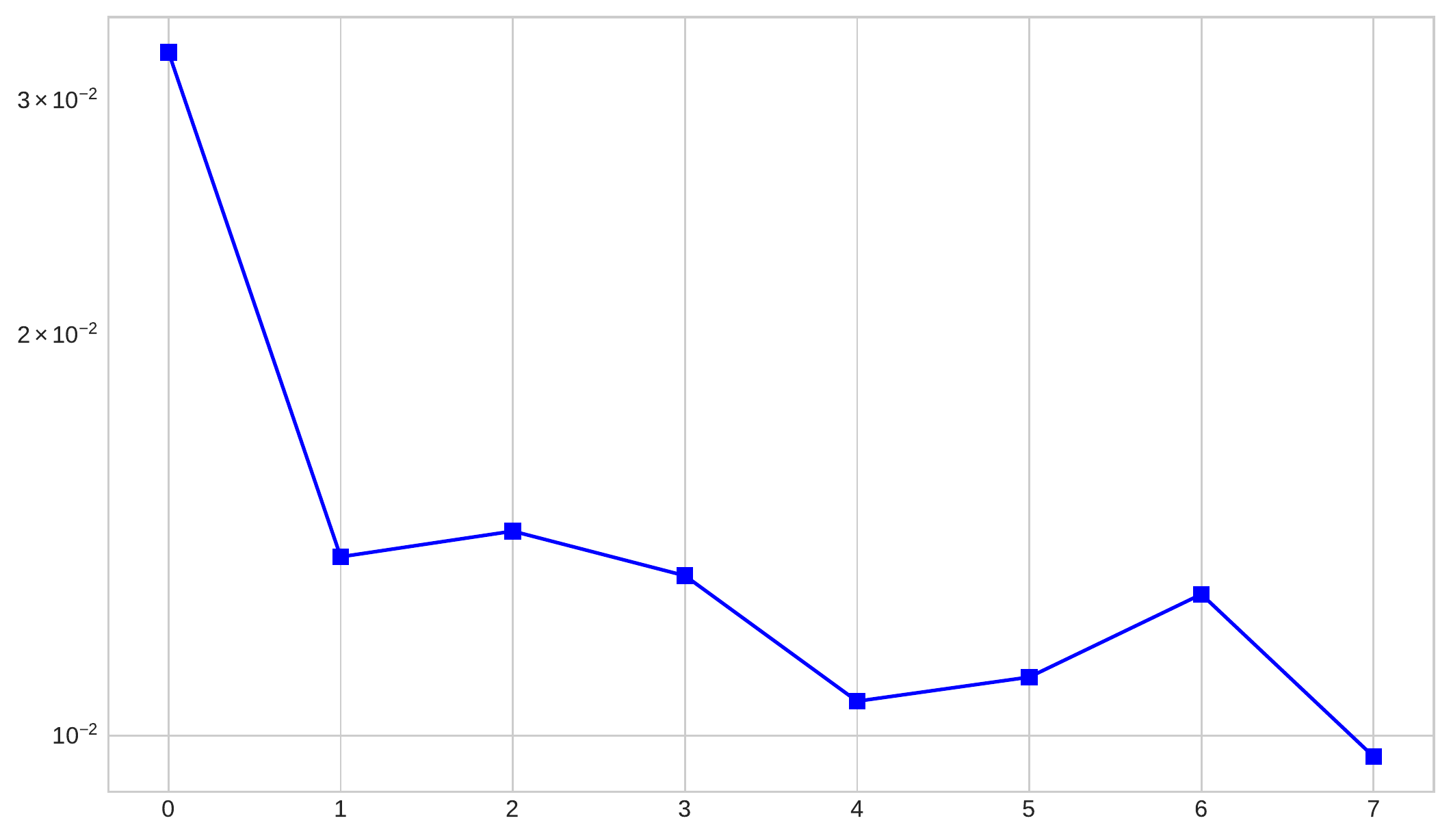}
               \put(35,60){\small G-FIPINN}
            \end{overpic}
         \end{center}
         \vspace{-0.2cm}
         \caption{Estimated failure probabilities over restarts  for four peaks problem ($N_{c} = 2000$).}
         \label{four_peak_failure_probability}
         \end{figure}

We  now consider the problem with four peaks, $(\pm 0.5, \pm 0.5)$, which has worse regularity and is more difficult to train. The true solution is depicted in Fig. \ref{two_peak_true_solution}(b).   To improve the expressivity of the neural network, we increase the number of hidden neurons to 128 while leaving the other parts unchanged in this case.  We notice a very similar behavior for the predict errors as in the earlier experiments.     In the same manner as before, we begin with $N_c=2000$ collocation points, and we provide the numerical results in Figs. \ref{four_peak_predictd_solution}-\ref{four_peak_failure_probability}.

We can clearly see from Fig.\ref{four_peak_predictd_solution} that the error is much smaller than the results obtained by the MC-FIPINN, particularly at the peaks. The collocation distribution  shown in Fig.\ref{four_peak_samples} can explain this phenomenon. These points can provide more useful information for  training and force the network to pay more attention in failure  region.   It is obvious that  during restarts, the failure probability will gradually decrease during restarts, as observed in Fig.\ref{four_peak_failure_probability}.  In conjunction with Fig. \ref{four_peak_compared_error}, we can see that after 7 restarts, the failure probability obtained by G-FIPINN algorithm is less than the tolerance and we can stop earlier. Only 25,000 training steps are needed  to train the network in this instance, which can reduce the computation cost of the training period.  The relative errors with various numbers of collocation points $N_c$ are shown in Fig. \ref{two_peak_error}(b). Similar to the  case of two peaks, larger values of $N_c$ perform better performance, but the AFI-PINNs produces superior outcomes.

\subsection{The time-dependent wave equation}

\begin{figure}[htbp]
\centering
\includegraphics[width = 0.4\textwidth]{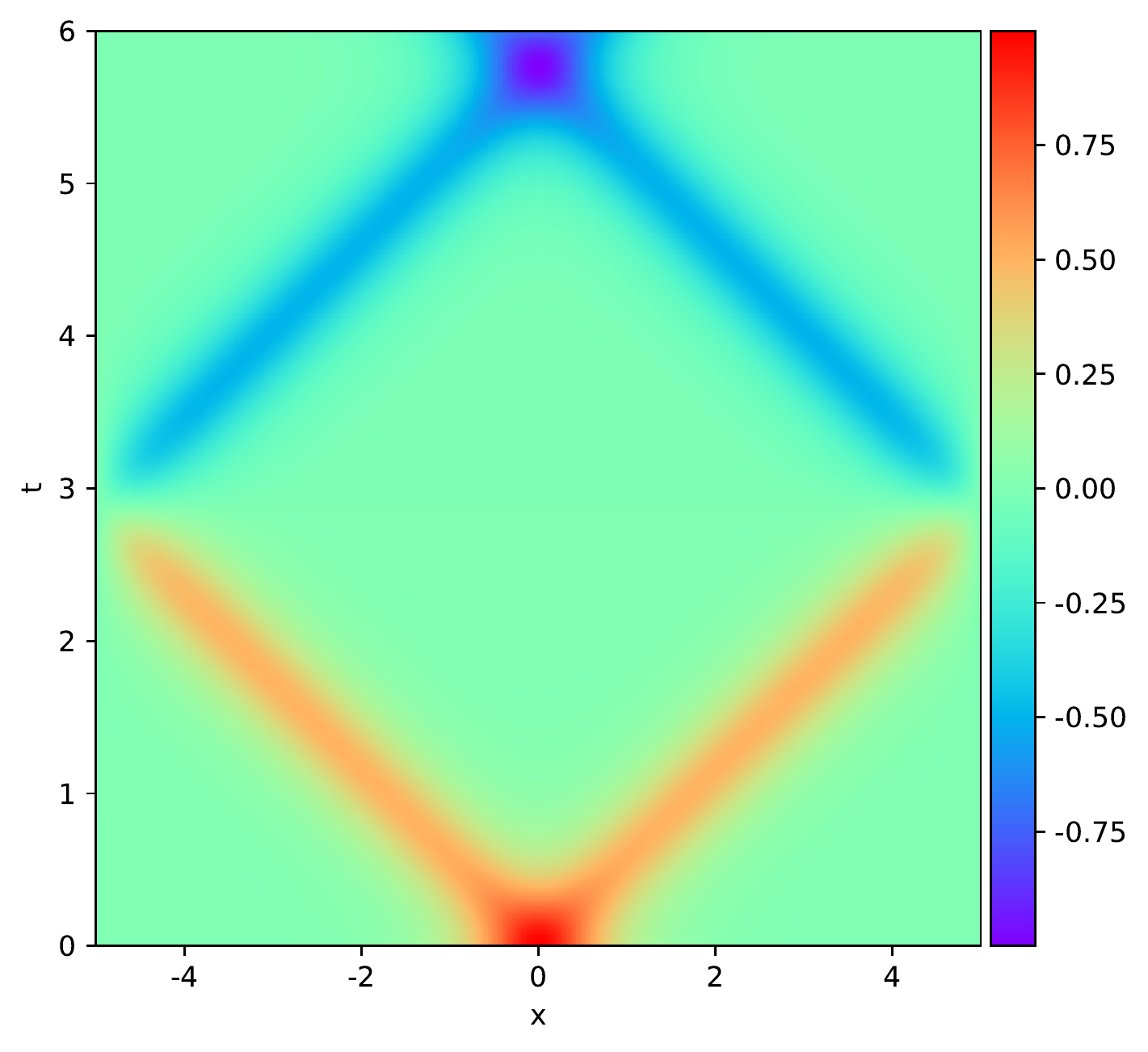}
\caption{The contour plot of the exact solution.}
\label{wave_true_solution}
\end{figure}

Here we consider the following one-dimensional time-dependent wave equation
\begin{equation}
\begin{split}
&\frac{\partial^{2}u}{\partial t^{2}} - 3\frac{\partial^{2}u}{\partial x}^{2} = 0, \quad (t, x) \in [0,6] \times [-5,5],\\
&u(0,x) = \frac{1}{\cosh(2x)} - \frac{0.5}{\cosh(2(x - 10))} - \frac{0.5}{\cosh(2(x+10))}, \\
&\frac{\partial u}{\partial t}(0,x) = 0, ,\\
&u(t,-5) = u(t,5) = 0,
\end{split}
\end{equation}
where the true solution is
\begin{equation}
\begin{split}
u(t,x) &= \frac{0.5}{\cosh(2(x - \sqrt{3}t))} - \frac{0.5}{\cosh(2(x - 10 + \sqrt{3}t))} \\
&+ \frac{0.5}{\cosh(2(x + \sqrt{3}t))} - \frac{0.5}{\cosh(2(x + 10 - \sqrt{3}t))}.
\end{split}
\end{equation}
We can clearly see that the true solution is approximately not connected, as demonstrated in Fig.\ref{wave_true_solution}. This problem is used to test effectiveness of our algorithm for the time dependent and multi-modal cases.  Here we set the failure probability tolerance to 0.0001.

\begin{figure}[htbp]
\centering
\includegraphics[width = 0.45\textwidth]{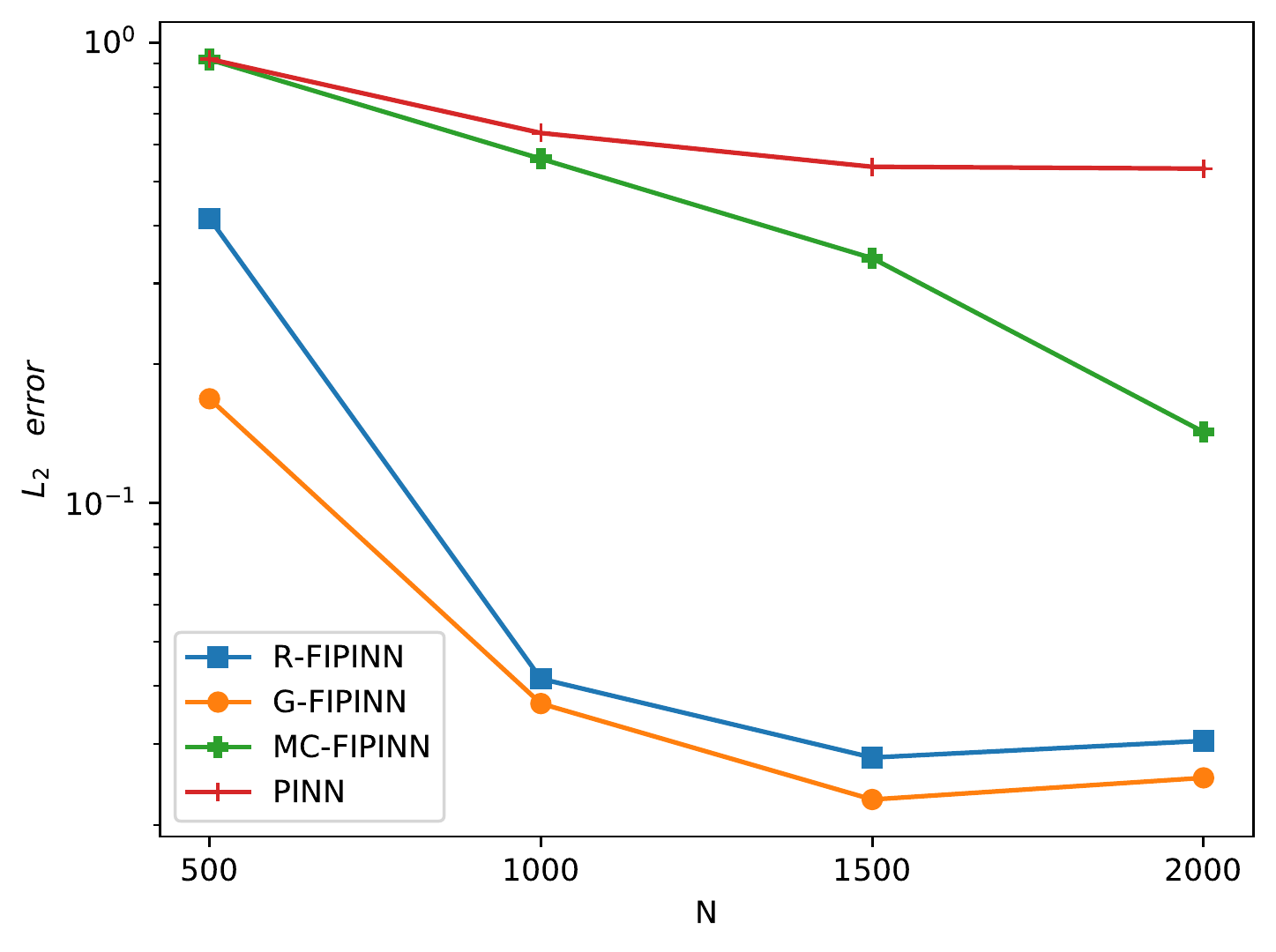}
\caption{ L2 error obtained by different methods when $N_{c}$ varies from 500 to 2000 for the wave problem.}
\label{wave_error}
\end{figure}

\begin{figure}[htbp]
\begin{center}
\begin{overpic}[width = 0.3\textwidth]{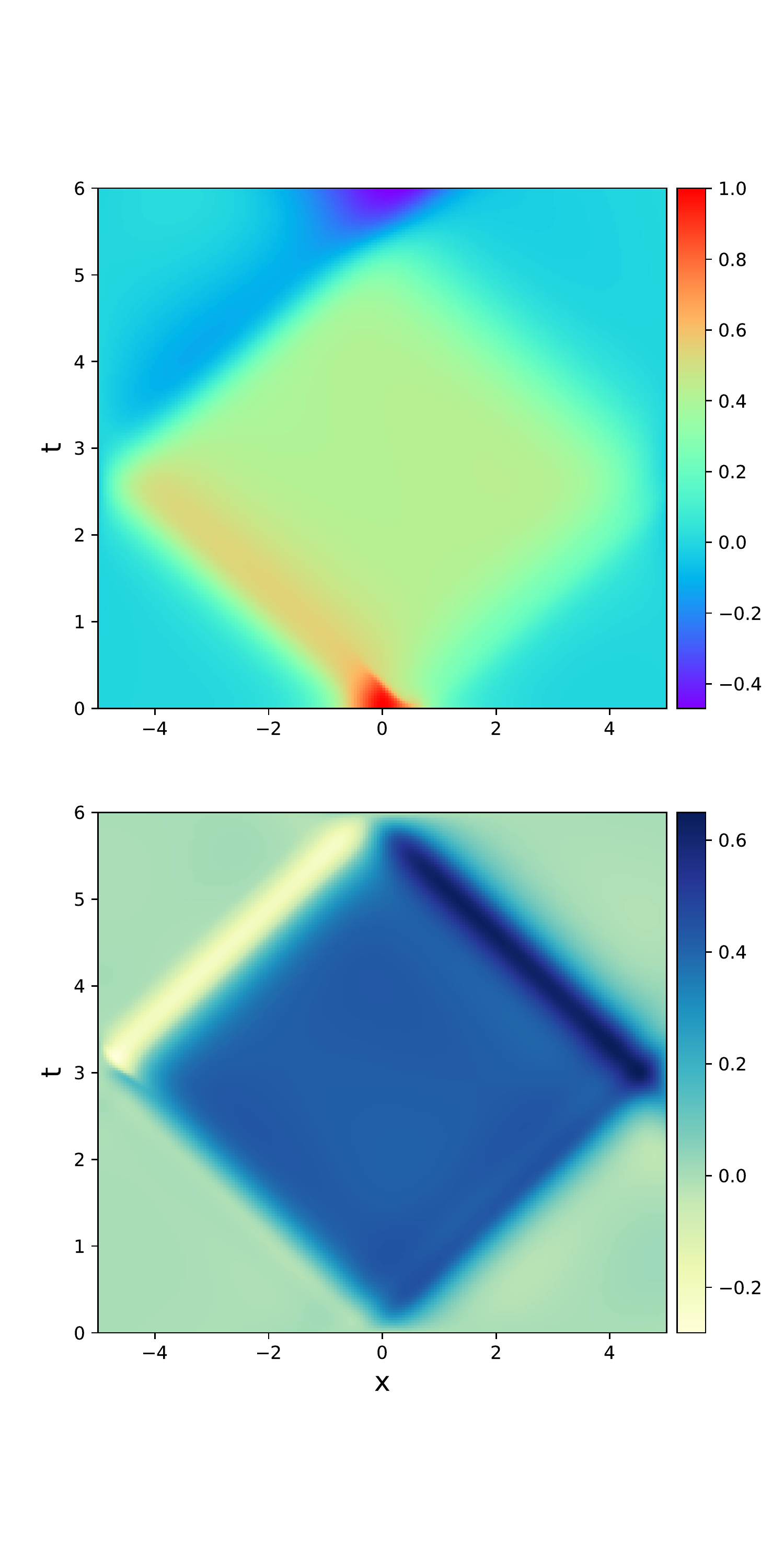}
   \put(15,90){\small MC-FIPINN}
\end{overpic}
\begin{overpic}[width = 0.3\textwidth]{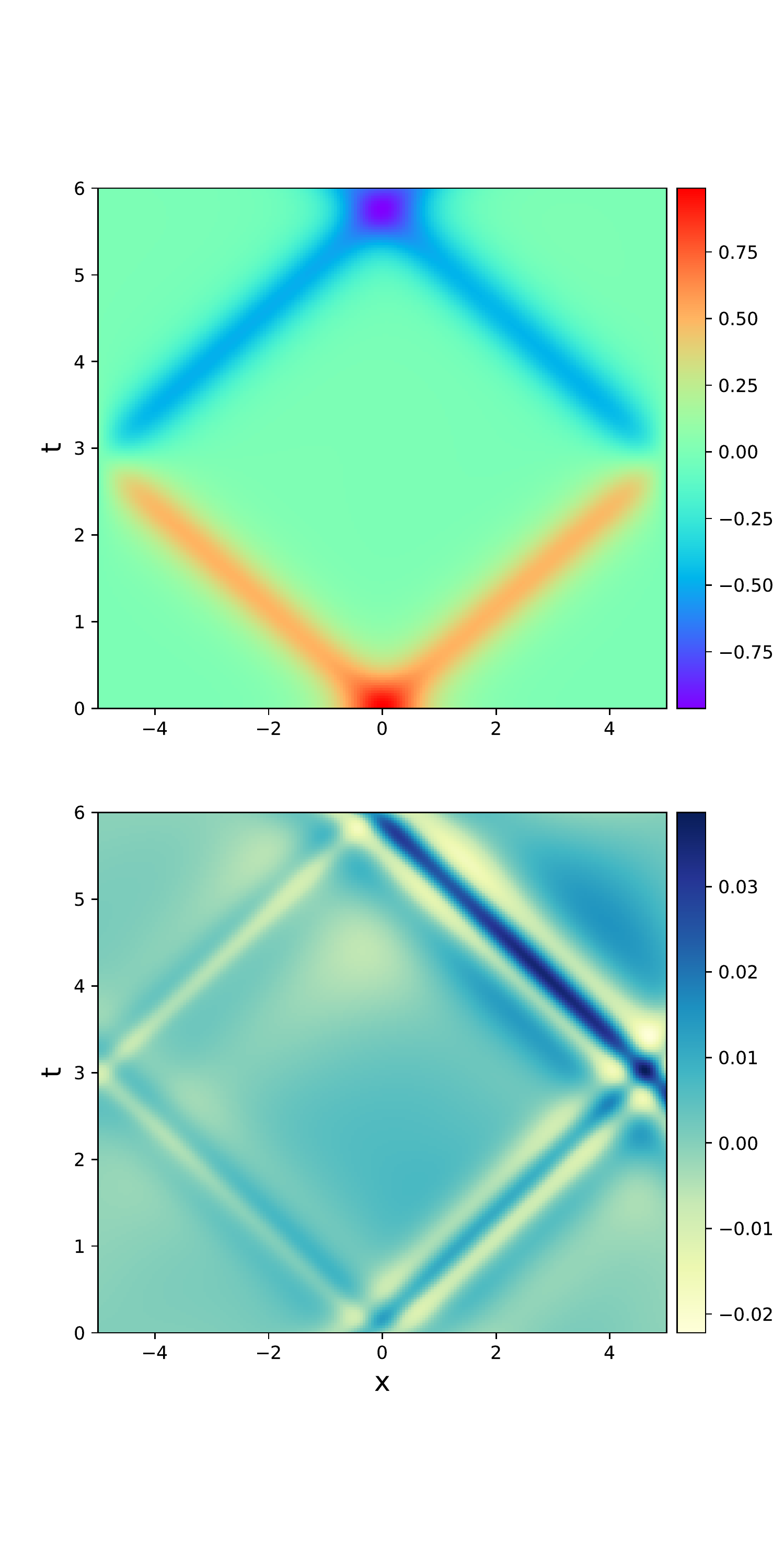}
   \put(15,90){\small R-FIPINN}
\end{overpic}
\begin{overpic}[width = 0.3\textwidth]{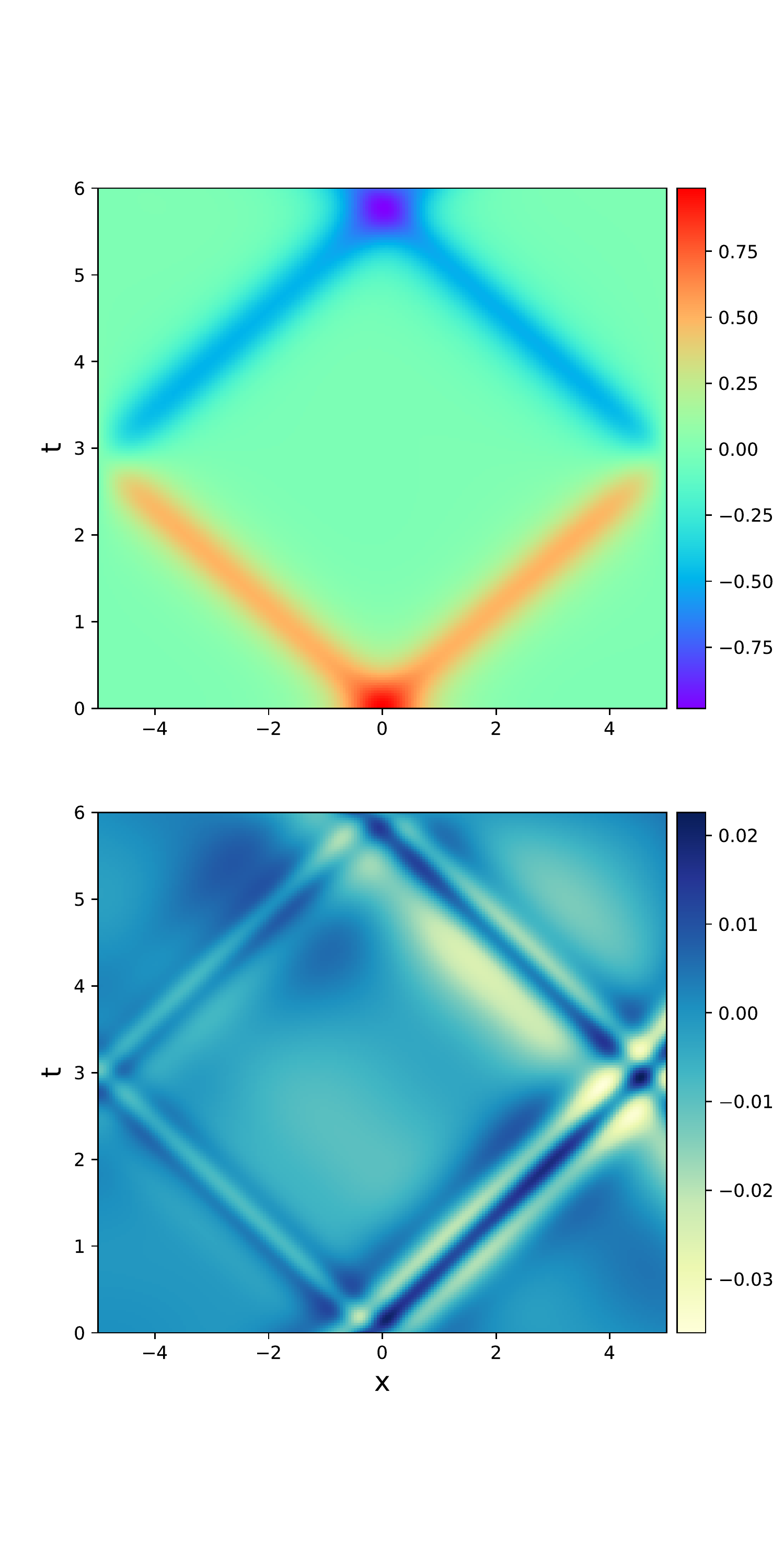}
   \put(15,90){\small G-FIPINN}
\end{overpic}
\end{center}
\vspace{-1cm}
\caption{The predicted solution (top) and the corresponding absolute error (bottom) obtained by different methods for  the wave problem ($N_{c} = 1000$).}
\label{wave_predictd_solution}
\end{figure}

\begin{figure}[htbp]
\begin{center}
   \begin{overpic}[width = 0.29\textwidth]{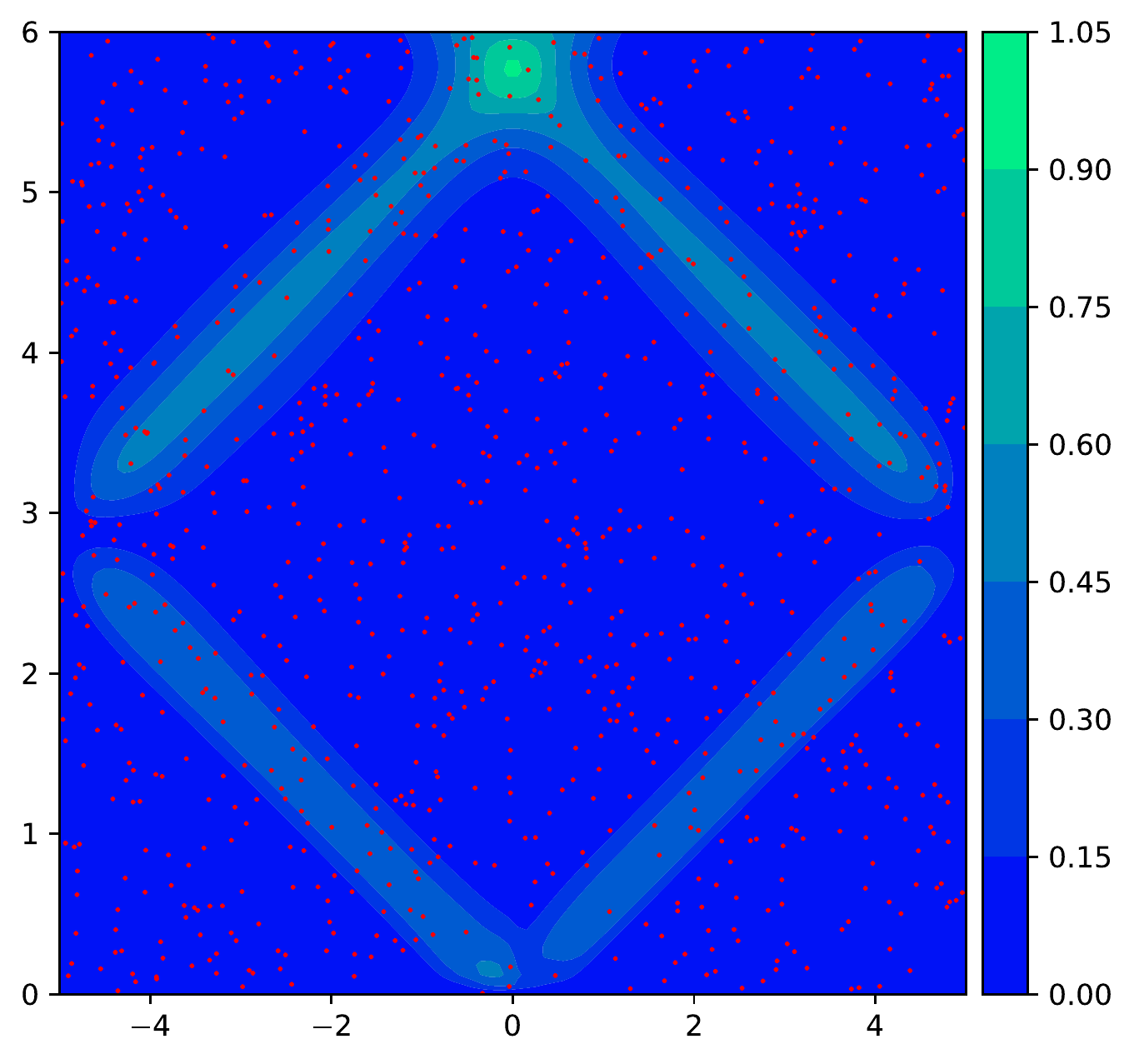}
      \put(35,95){\small MC-FIPINN}
   \end{overpic}
   \begin{overpic}[width = 0.29\textwidth]{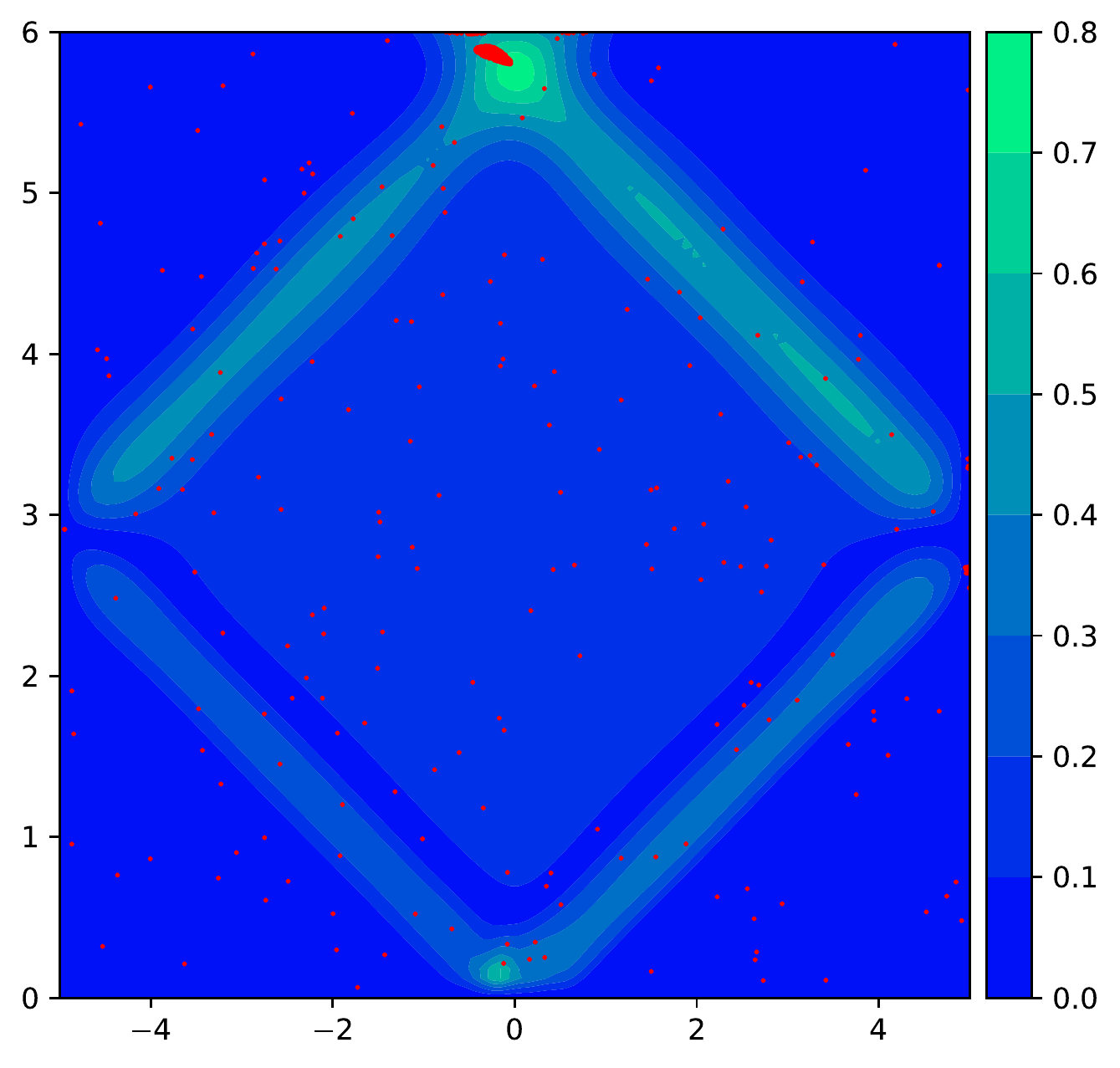}
      \put(35,95){\small R-FIPINN}
   \end{overpic}
   \begin{overpic}[width = 0.29\textwidth]{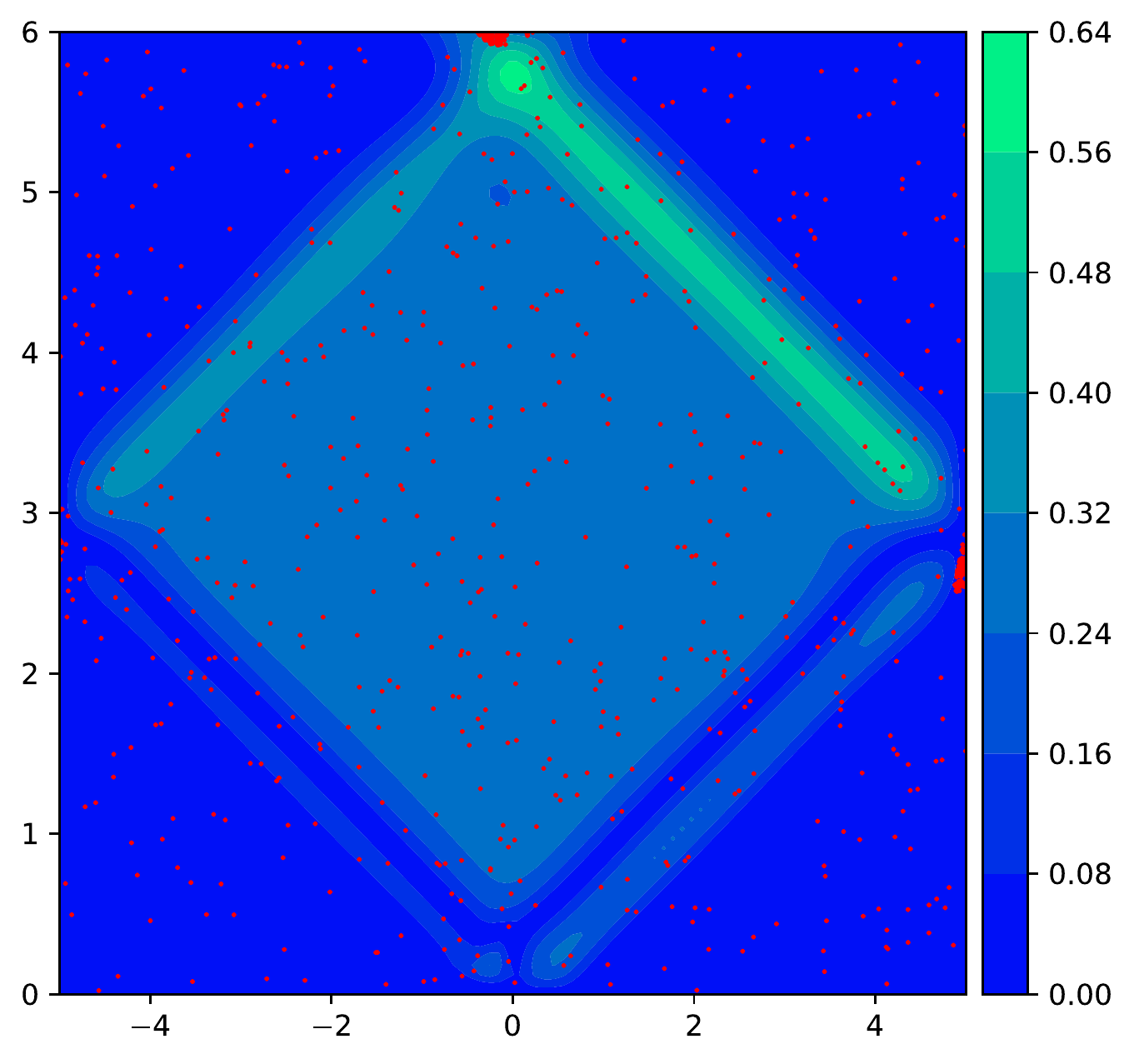}
      \put(35,95){\small G-FIPINN}
   \end{overpic}
\end{center}
\vspace{-0.2cm}
\caption{Final distribution of collocation points  for the wave problem ($N_{c}=1000$). }
\label{wave_samples}
\end{figure}

  We show the relative errors for the baseline and AFI-PINNs schemes with different numbers of collocation points $N_c$ in Fig.\ref{wave_error}.  The performance is good for larger values of $N_c$, just like the elliptic PDE system. Compared to standard PINNs and MC-FIPINN, the AFI-PINNs consistently produce errors that are up to an order of magnitude lower.  Additionally, we can see that even with 1000 collocation points, our method can still produce a respectable performance.   In Fig.\ref{wave_predictd_solution}, we report the predicted solution  and the corresponding absolute error  obtained by different methods using $N_{c} = 1000$ collocation points.  The result can be observed more clearly in  that the predicted solution obtained from AFI-PINNs is nearly identical to the true solution.   This phenomenon confirms that by using a small number of collocation points, our AFI-PINNs can handle problems of this type that are time dependent and have a large problem domain. The final distribution of collocation points obtained by AFI-PINNs,  as shown in Fig.\ref{wave_samples}, concentrate on the areas with higher performance function values, which explains why this occurs. As a result, the network can focus more attention on failure regions to improve performance.

\vspace{0.2cm}

\begin{figure}[htbp]
   \begin{center}
      \begin{overpic}[width = 0.3\textwidth]{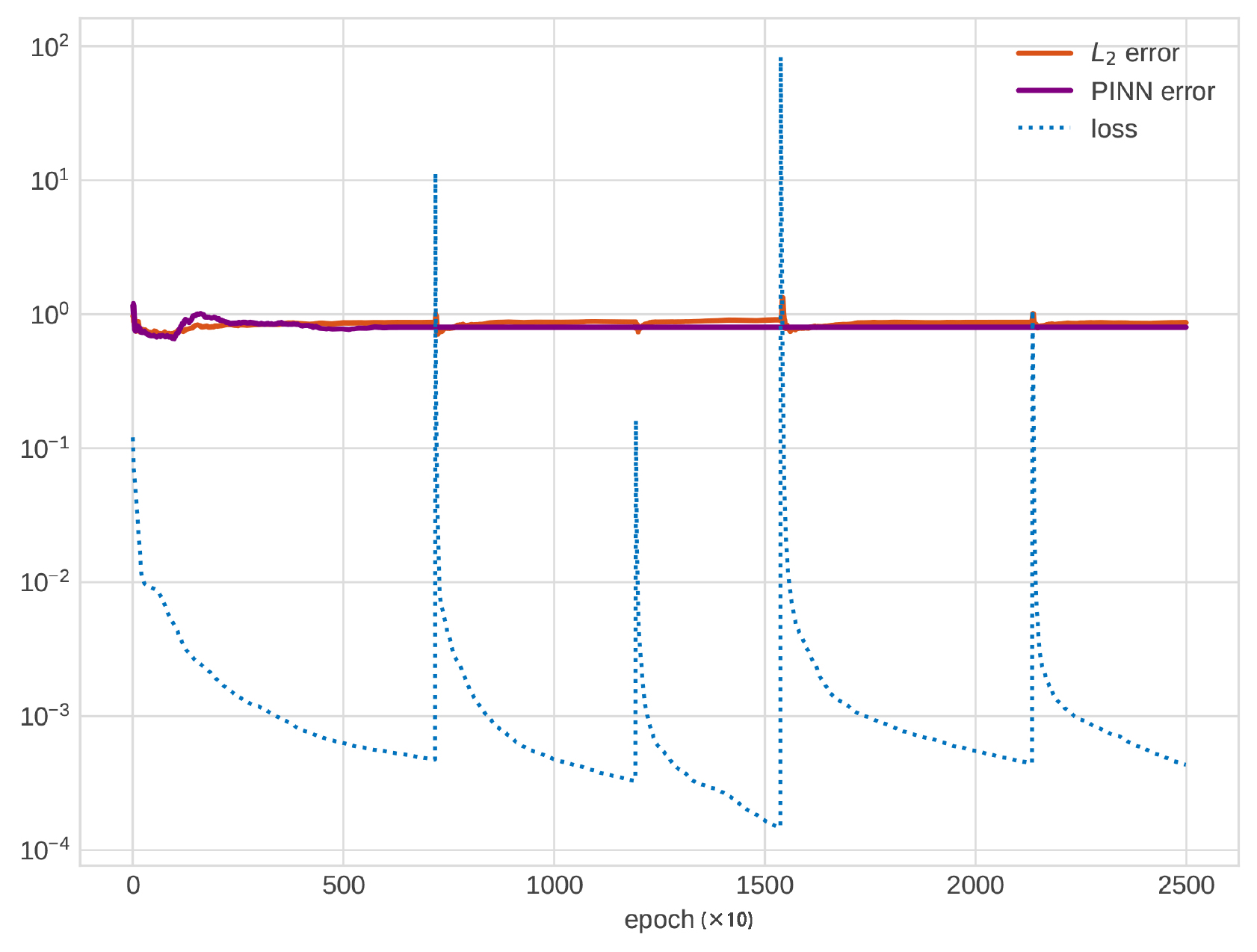}
         \put(35,76){\small MC-FIPINN}
      \end{overpic}
      \begin{overpic}[width = 0.3\textwidth]{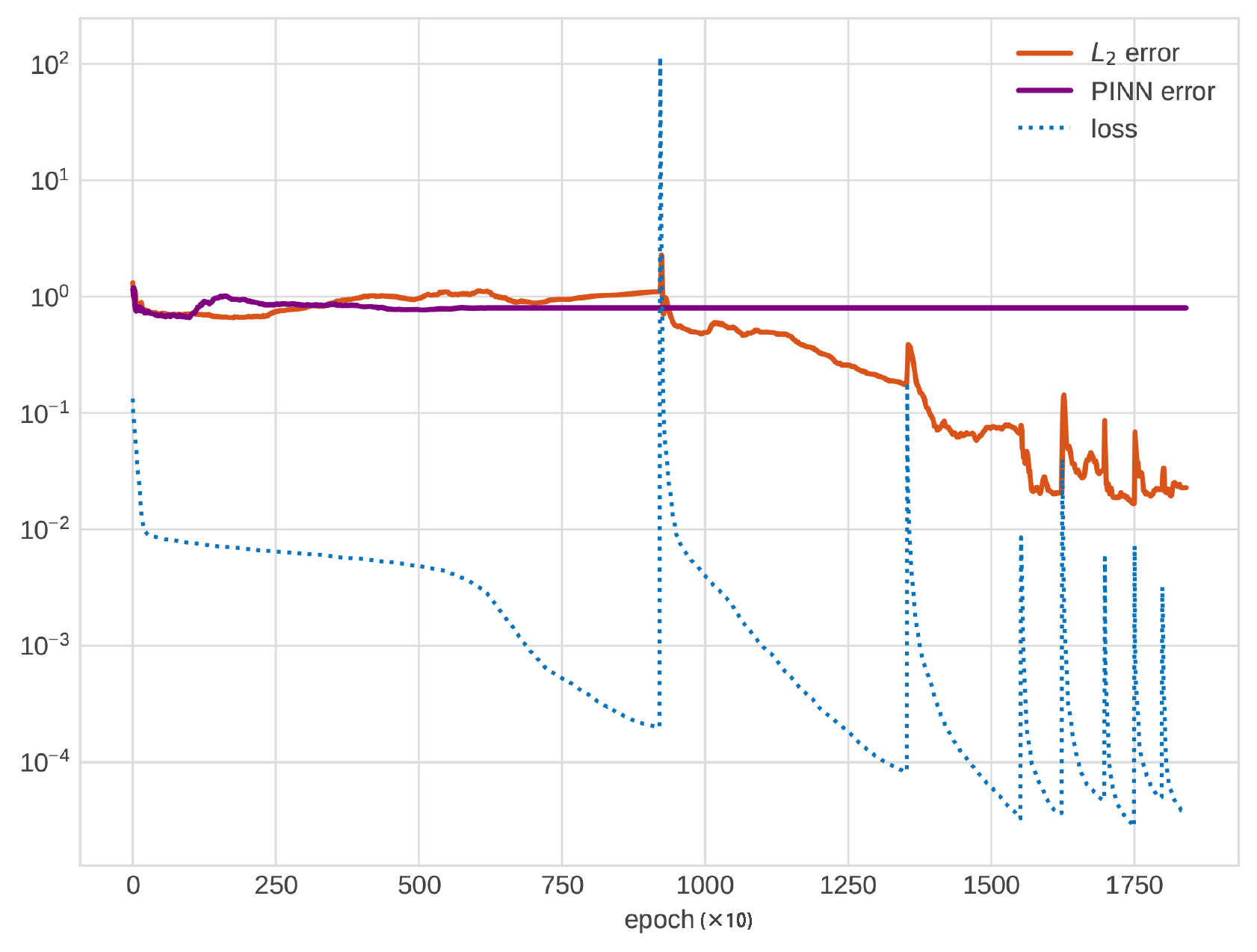}
         \put(35,76){\small R-FIPINN}
      \end{overpic}
      \begin{overpic}[width = 0.3\textwidth]{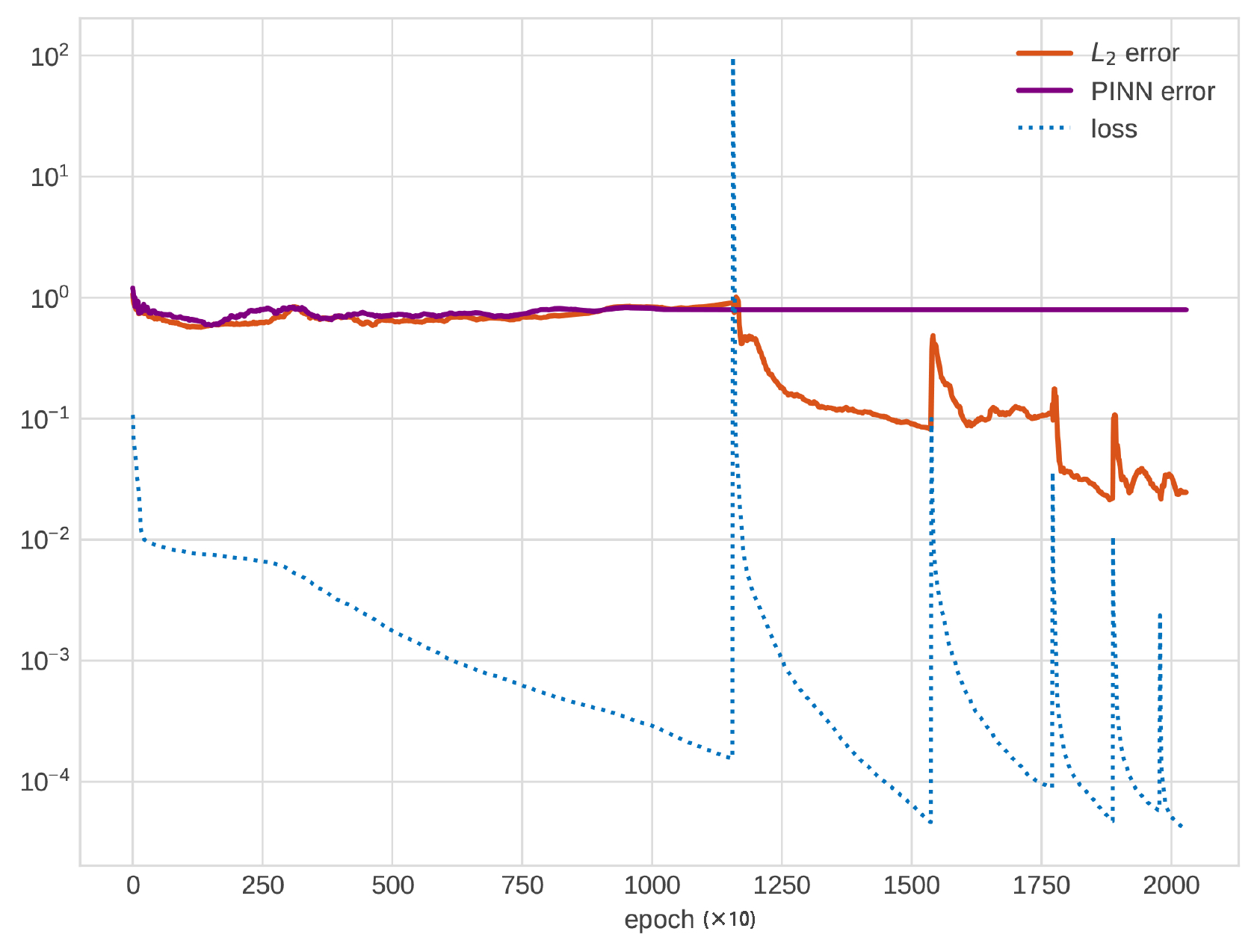}
         \put(35,76){\small G-FIPINN}
      \end{overpic}
   \end{center}
   \vspace{-0.2cm}
   \caption{Relative error and the training loss compared to the vanilla PINN error every 10 epochs during the training process for the wave problem.}
   \label{wave_compared_error}
   \end{figure}

\begin{figure}[htbp]
      \begin{center}
         \begin{overpic}[width = 0.40\textwidth]{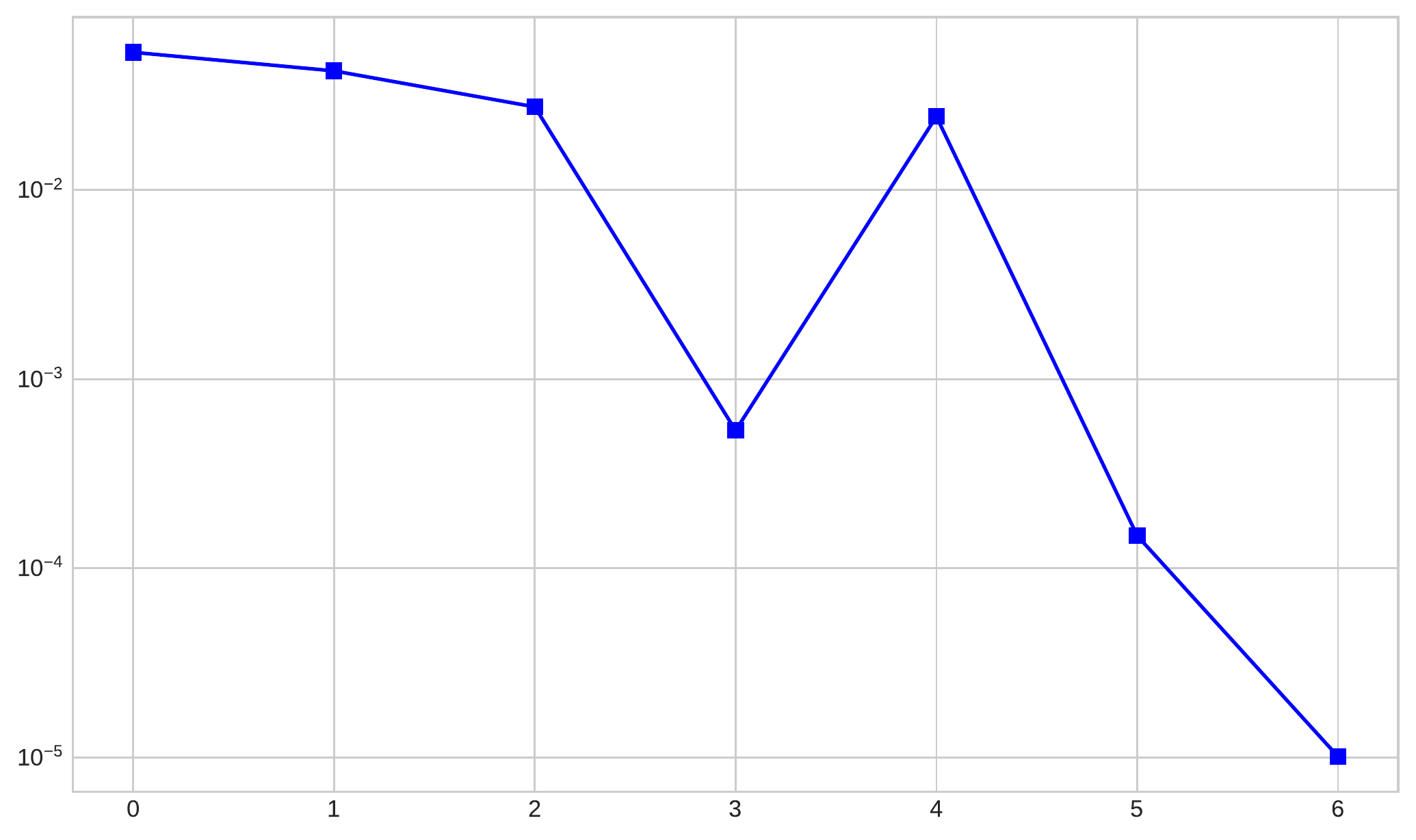}
            \put(35,60){\small R-FIPINN}
         \end{overpic}
         \begin{overpic}[width = 0.40\textwidth]{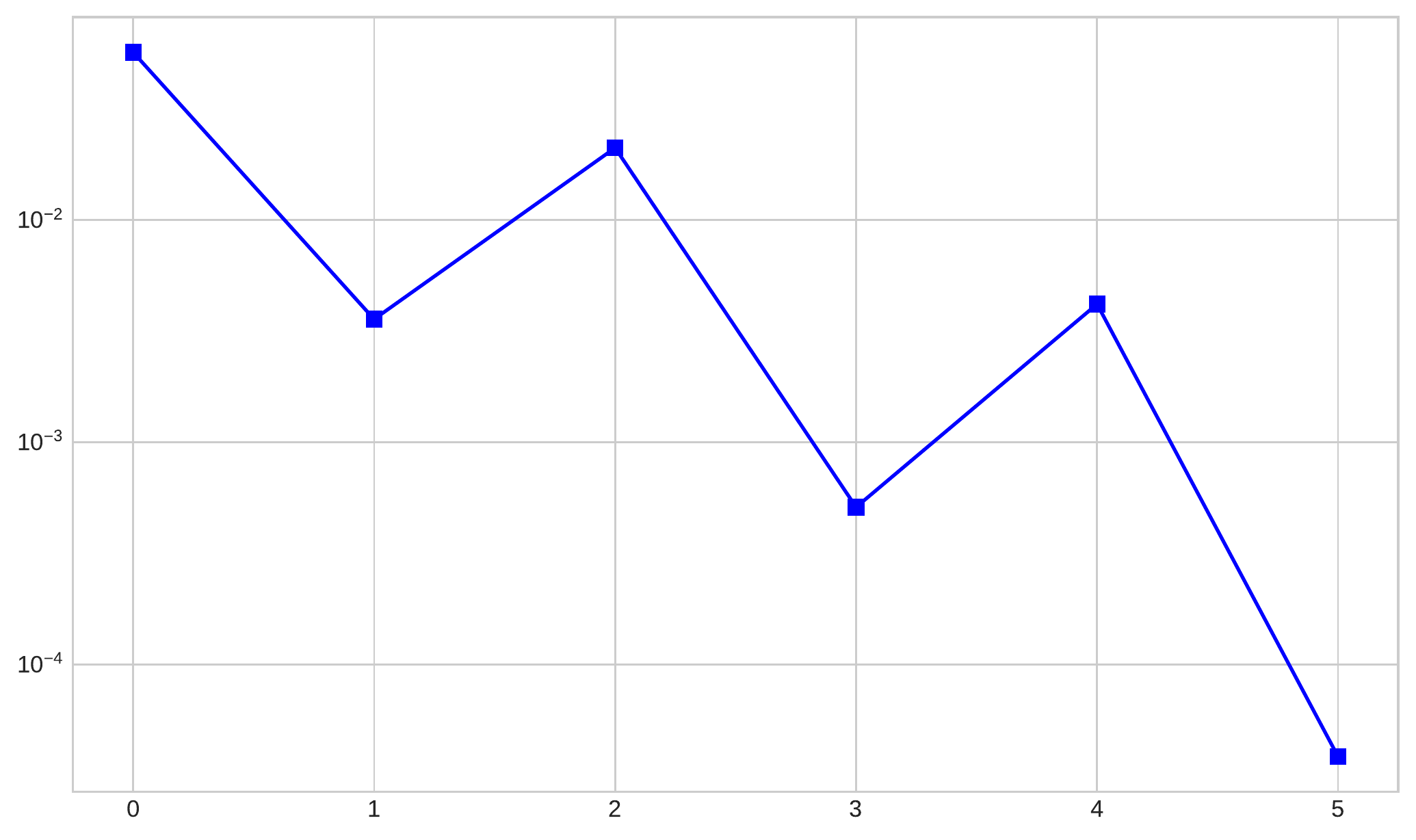}
            \put(35,60){\small G-FIPINN}
         \end{overpic}
      \end{center}
      \vspace{-0.2cm}
      \caption{Estimated failure probabilities over restarts  for the wave problem.}
      \label{wave_failure_probability}
      \end{figure}

We plot the predicted error and training loss during the training process in in Fig.\ref{wave_compared_error} to demonstrate the annealing restart speeds up convergence. It is evident that after each restart, the loss will sharply increase at first before rapidly decreasing. At the same time, as training with our method progresses, the predicted error will gradually decline. Even with numerous restarts, the predicted error obtained by MC-PINNs, however, still follows a straight line. Additionally, the restart can reduce computational costs because it keeps the size of the entire training set constant while dynamically updating the collocation dataset.  Additionally, by restarting, the system will gradually become more reliable, resulting in a decrease in the failure probability as shown in Fig.\ref{wave_failure_probability}, which can also save computational costs by stopping early.

\vspace{0.2cm}
\begin{figure}[htbp]
\centering
\includegraphics[width = 0.45\textwidth]{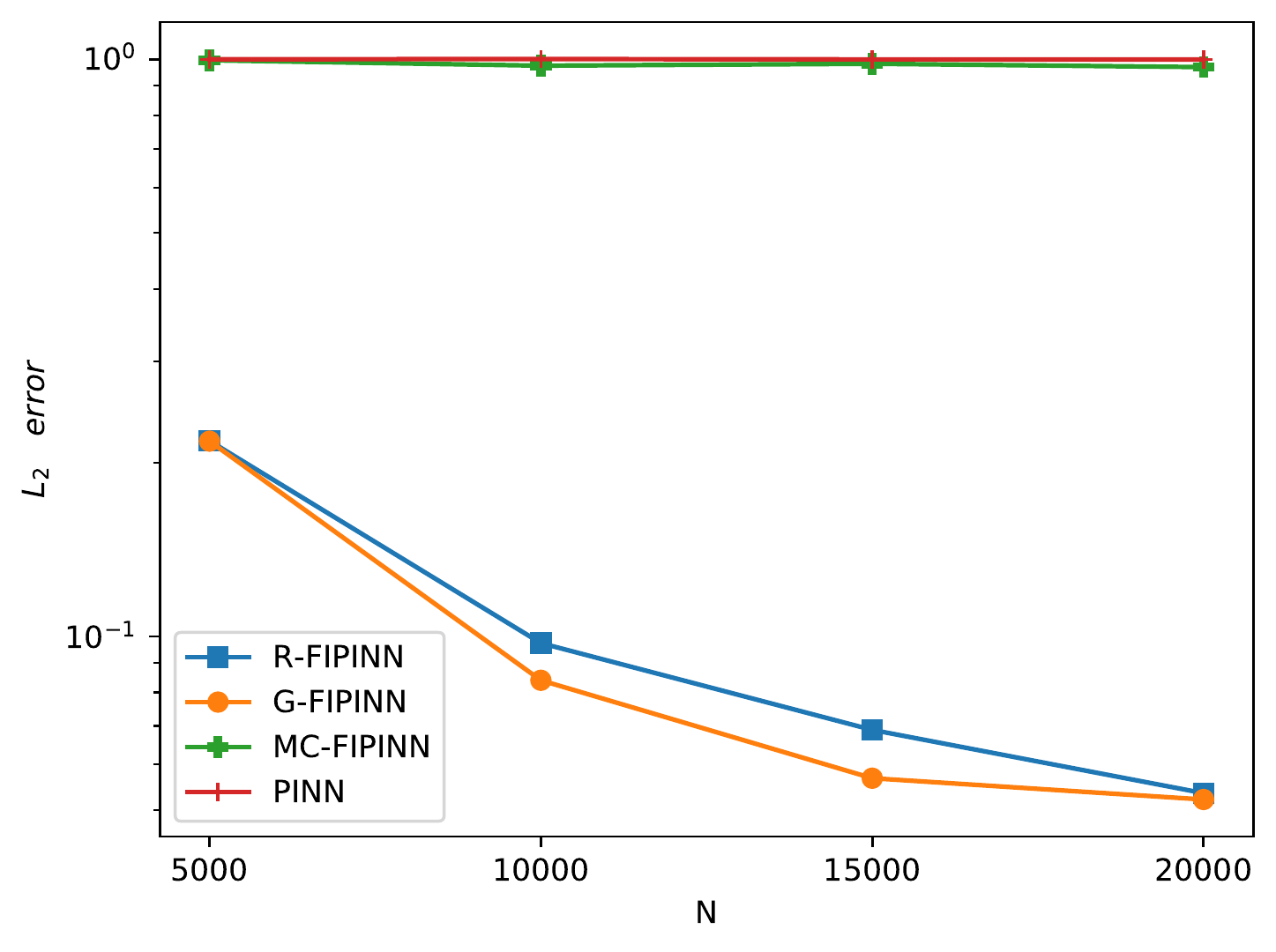}
\caption{L2 error obtained by different methods when $N_{c}$ varies from 5000 to 20000 for the high dimensional problem.}
\label{high_dimensional_error}
\end{figure}

\subsection{High dimensional Poisson equation}

Here we consider a 10-dimensional Poisson equation,
\begin{equation}
-\bigtriangleup u(\newx) = f(\newx),\quad \newx \in [-1,1]^{10},
\end{equation}
where the exact solution is
\begin{equation}
u(\newx) = e^{-10\|\newx\|_{2}^{2}}.
\end{equation}
The boundary conditions can be obtained through the true solution.
This problem is used to test the algorithm under the circumstance of high dimension.

In this example, we use a fully connected network with 7 hidden layers and 32 neurons in each layer. The \textit{Adam} optimizer with maximum 10,0000 iterations and a learning rate 0.005 is used to train the network. We chose to carry out the experiment with a range of collocation dataset sizes, including $\{5000, 10000, 15000,  20000\}$, in order to demonstrate the effectiveness of our algorithm.  The failure probability tolerance  is set to 0.0001 in this case.

Our method significantly outperforms the vanilla PINNs and MC-FIPINN methods with a much smaller predicted error, as shown in Fig.\ref{high_dimensional_error}. Additionally, our method can still perform well with fewer samples. This is crucial for high dimensional PDEs because it can significantly reduce computation costs without increasing the number of collocation samples. Furthermore, our predicted error will continue to decrease as more  collocation data points collected, in contrast to the straight line predicted errors obtained by the vanilla PINNs and MC-FIPINN.

\begin{figure}[htbp]
\begin{center}
   \begin{overpic}[width = 0.3\textwidth]{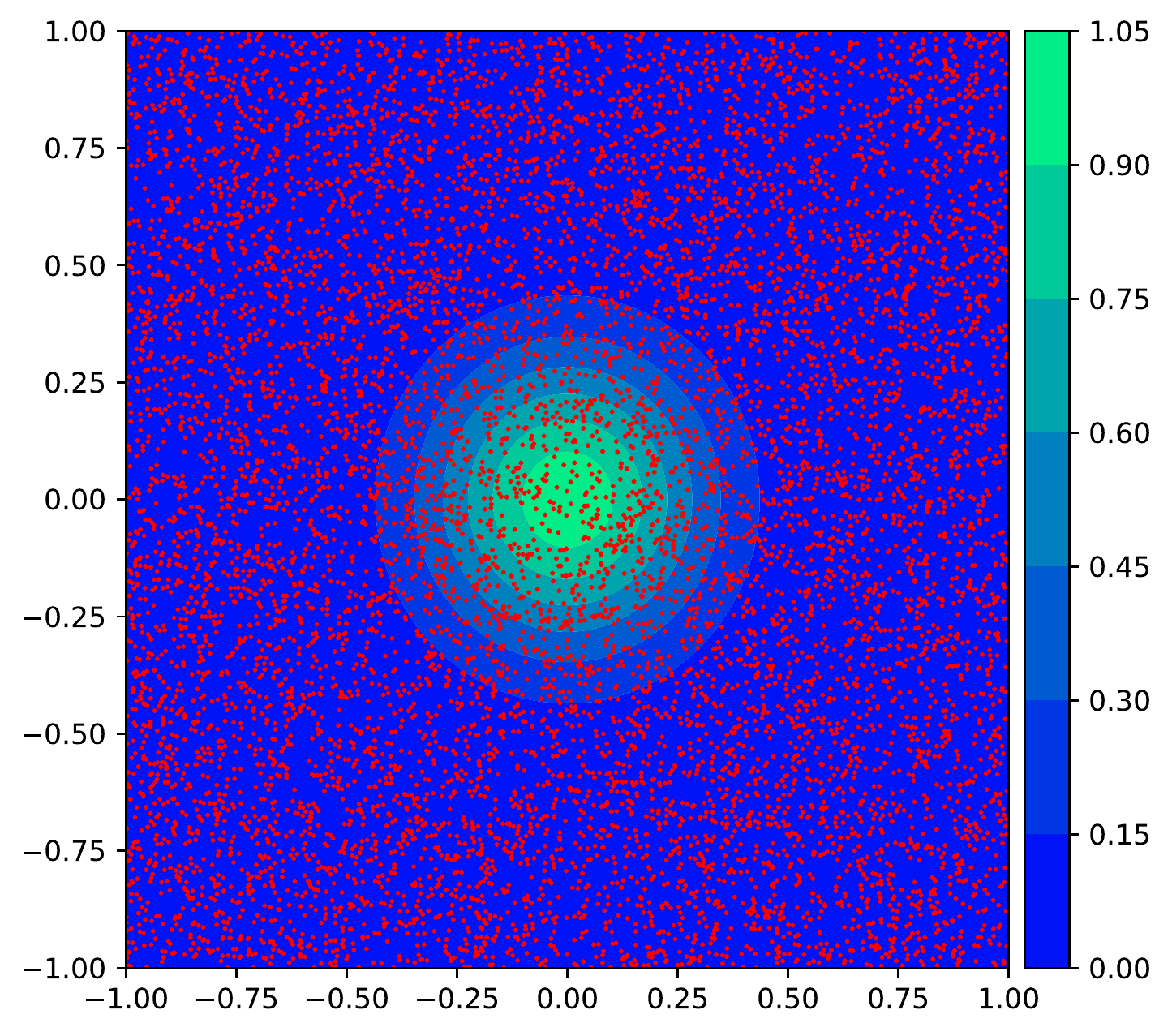}
      \put(35,90){\small MC-FIPINN}
   \end{overpic}
   \begin{overpic}[width = 0.3\textwidth]{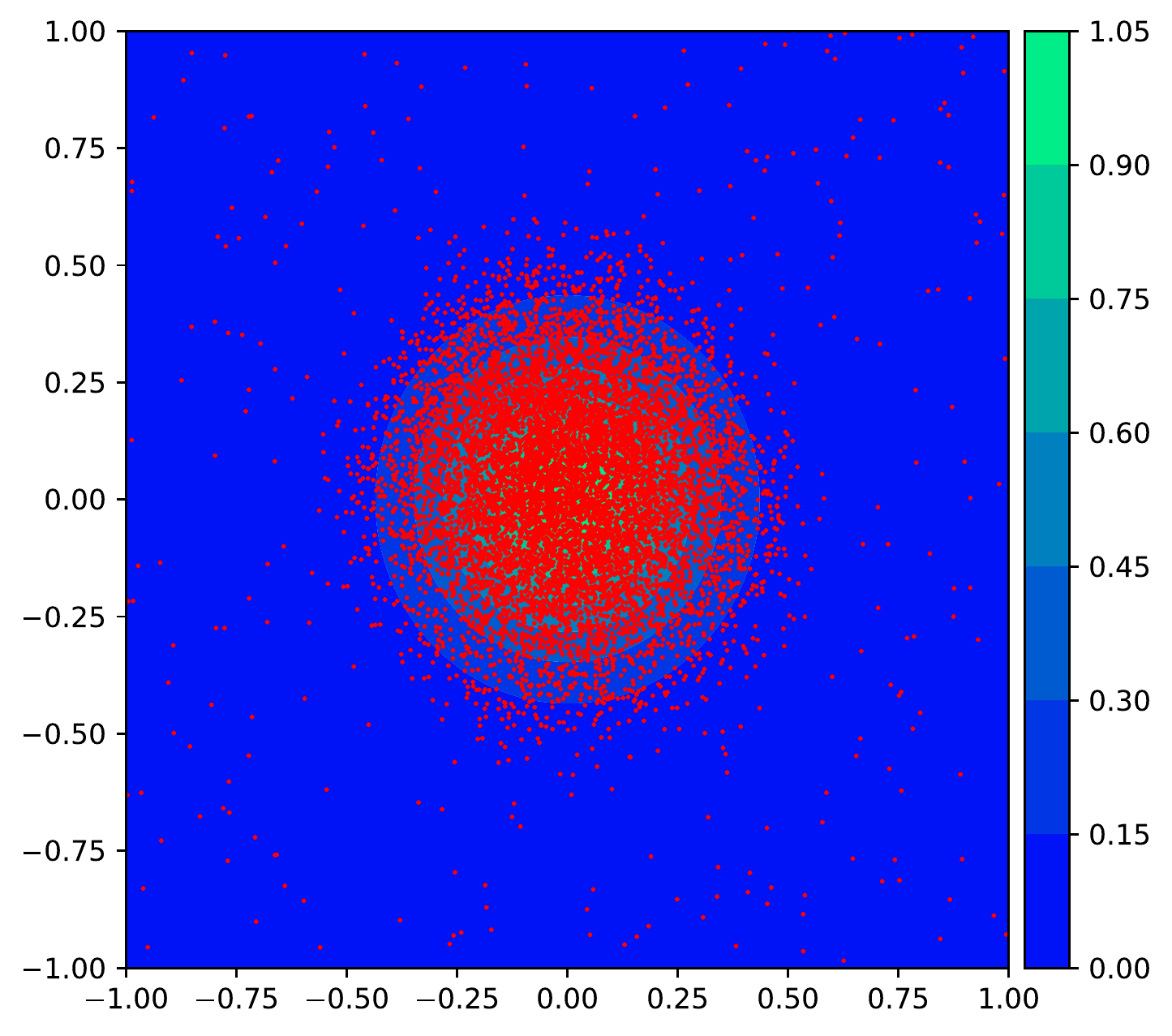}
      \put(35,90){\small R-FIPINN}
   \end{overpic}
   \begin{overpic}[width = 0.3\textwidth]{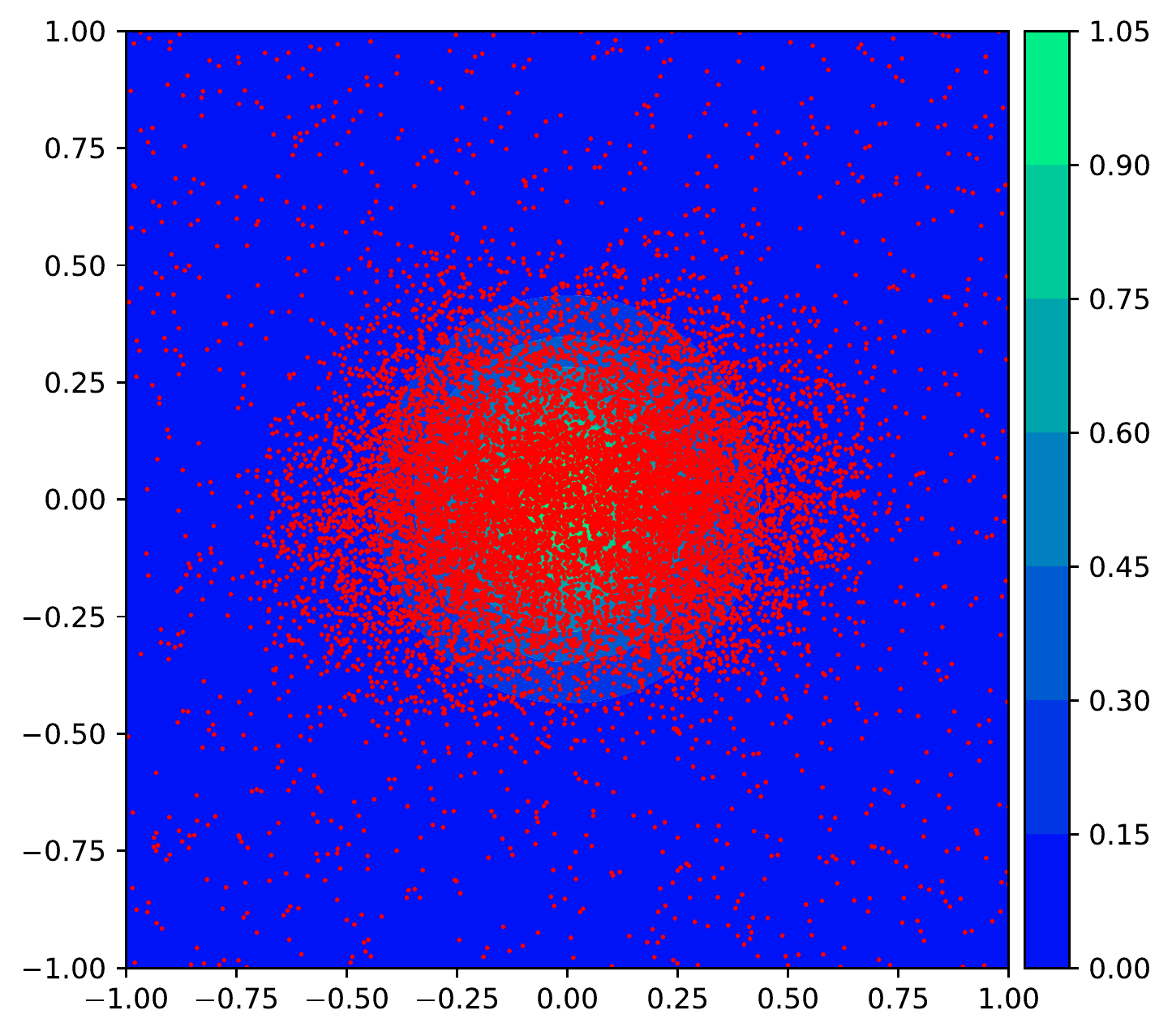}
      \put(35,90){\small G-FIPINN}
   \end{overpic}
\end{center}
\vspace{-0.2cm}
\caption{Final  collocation points  of $x_8, x_9$ for the high dimensional case ($N_{c} = 10000$).}
\label{high_dimensional_samples}
\end{figure}

\vspace{0.3cm}
\begin{figure}[htbp]
   \begin{center}
      \begin{overpic}[width = 0.32\textwidth]{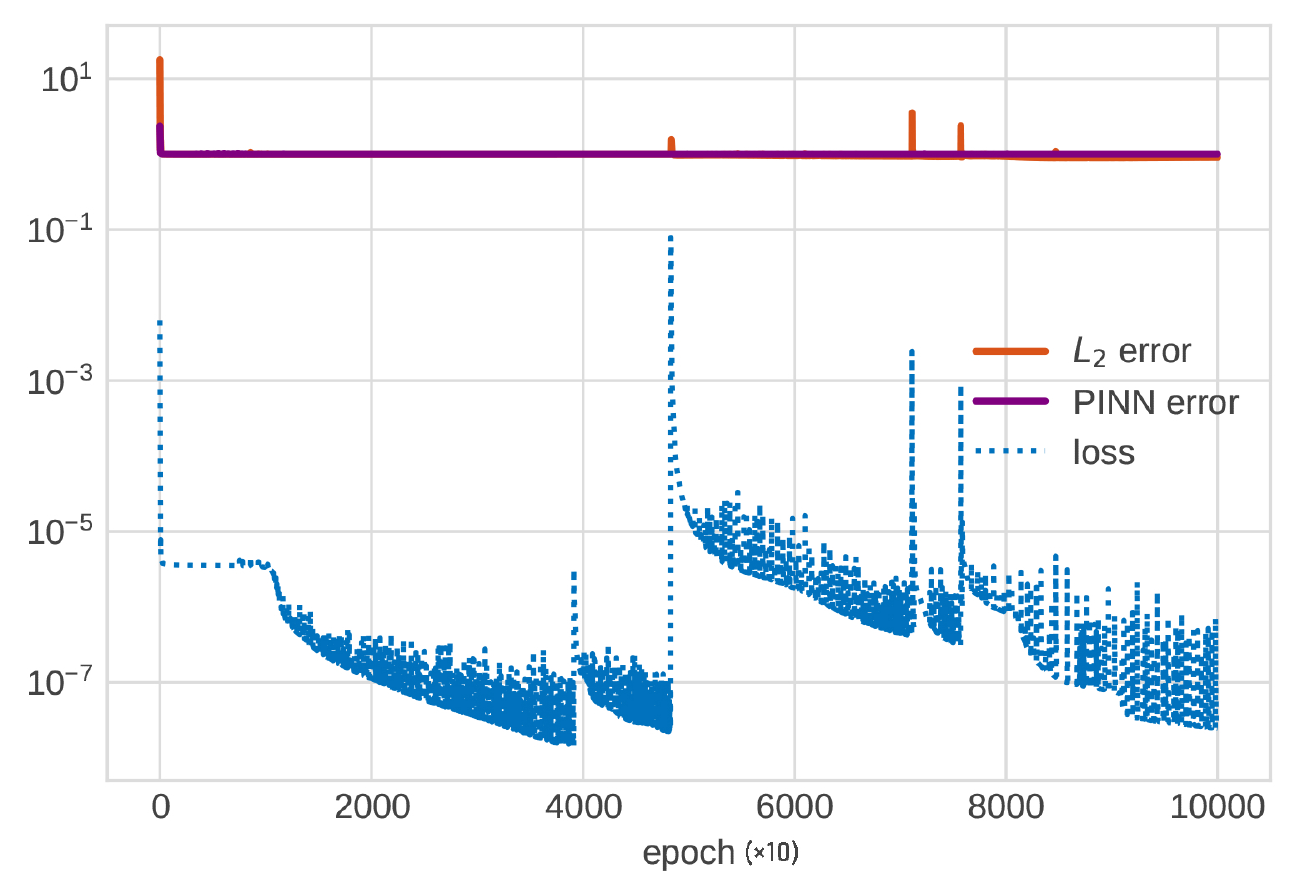}
         \put(35,70){\small MC-FIPINN}
      \end{overpic}
      \begin{overpic}[width = 0.32\textwidth]{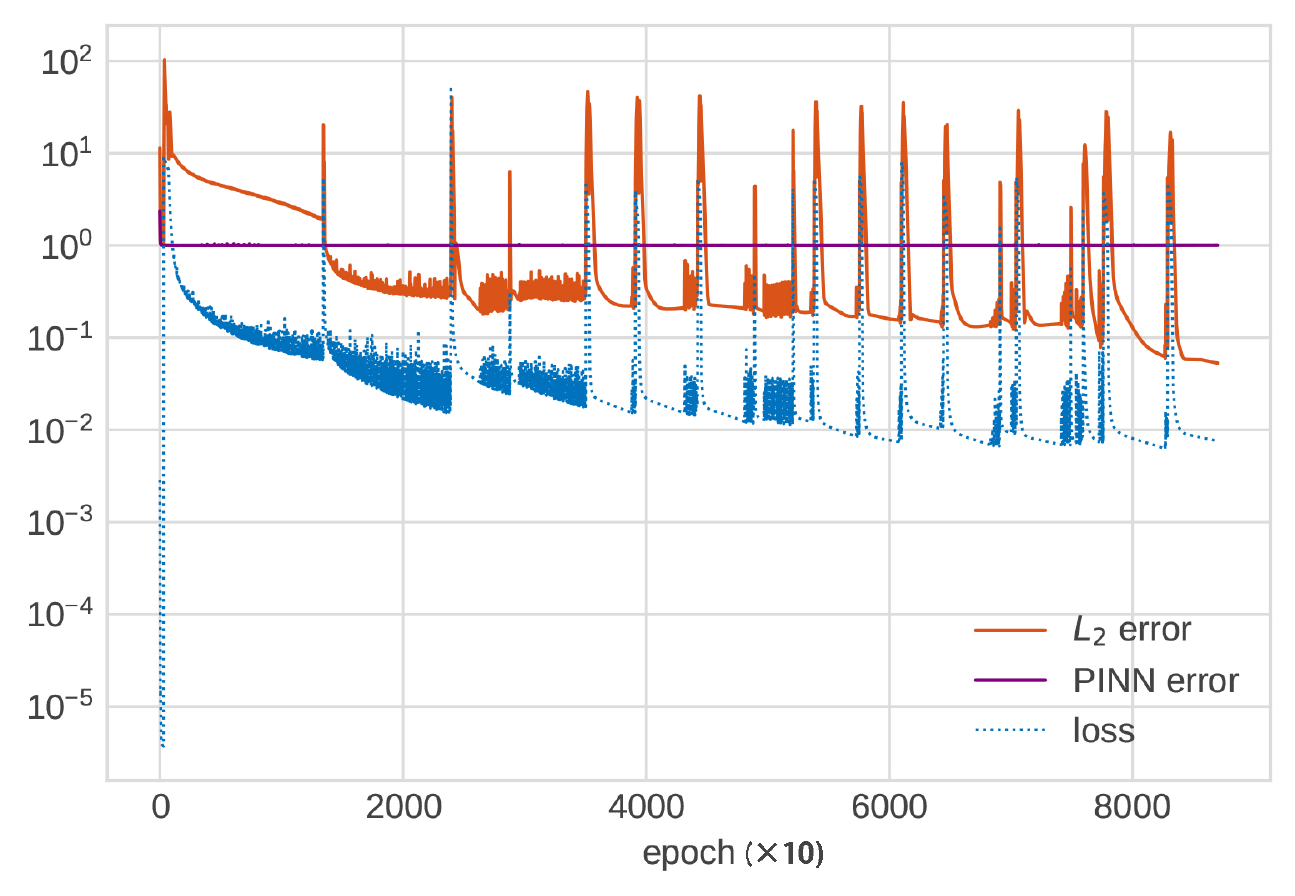}
         \put(35,70){\small R-FIPINN}
      \end{overpic}
      \begin{overpic}[width = 0.32\textwidth]{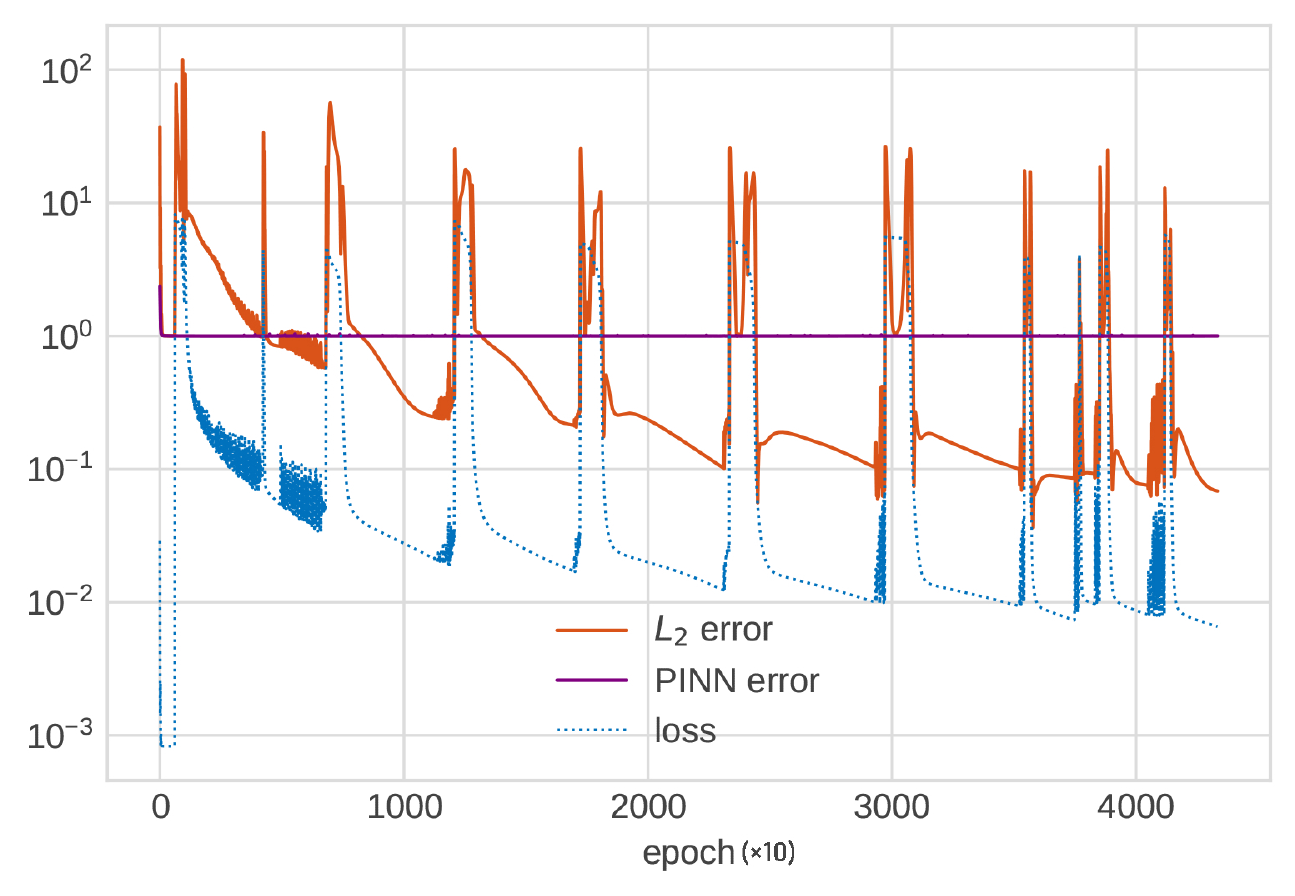}
         \put(35,70){\small G-FIPINN}
      \end{overpic}
   \end{center}
   \vspace{-0.2cm}
   \caption{Relative error and the training loss  during the training process  for the high dimensional case ($N_{c} = 10000$).}
   \label{high_dimensional_compared_error}
   \end{figure}

\vspace{0.3cm}
   \begin{figure}[h]
      \begin{center}
         \begin{overpic}[width = 0.45\textwidth]{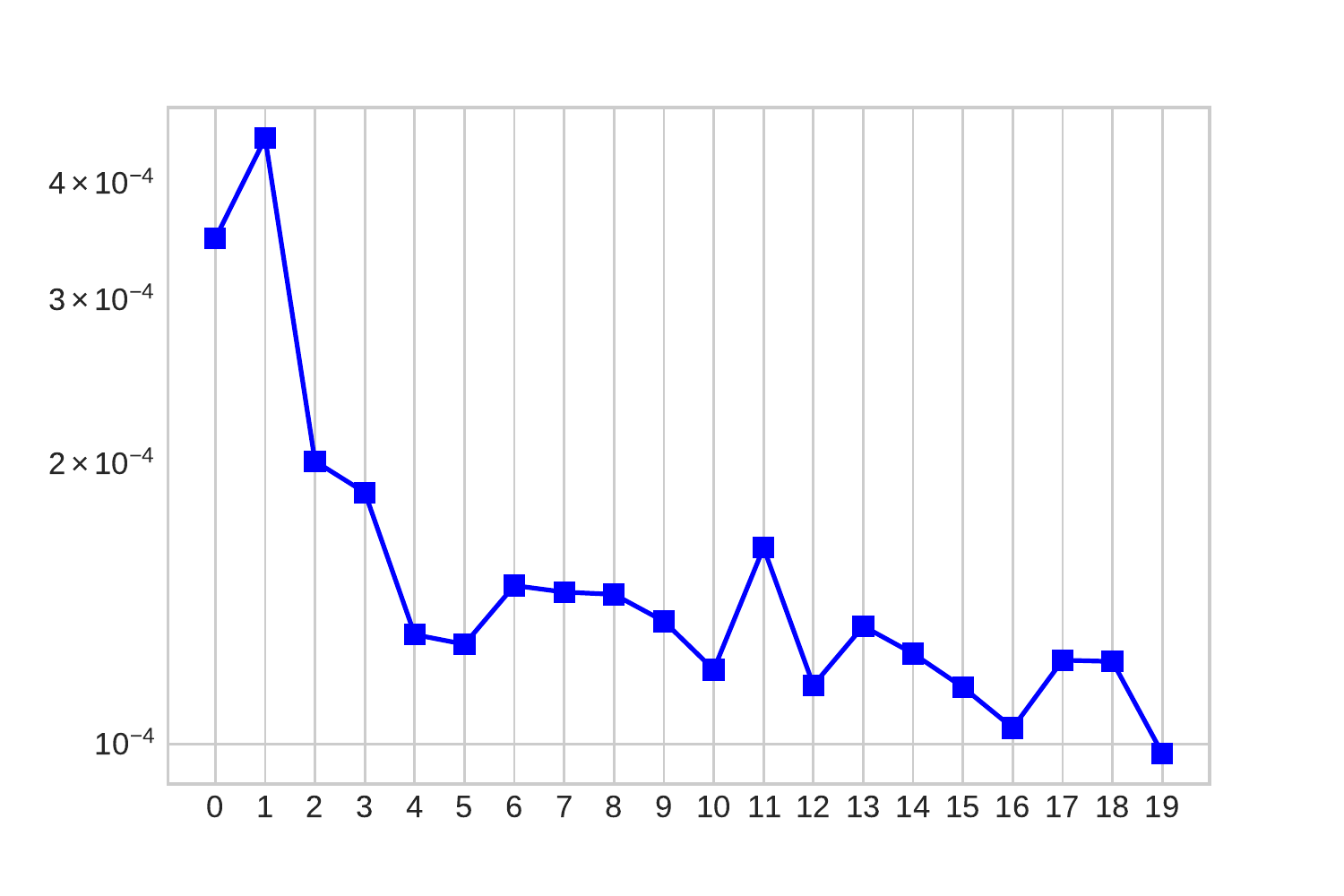}
            \put(35,63){\small R-FIPINN}
         \end{overpic}
         \begin{overpic}[width = 0.45\textwidth]{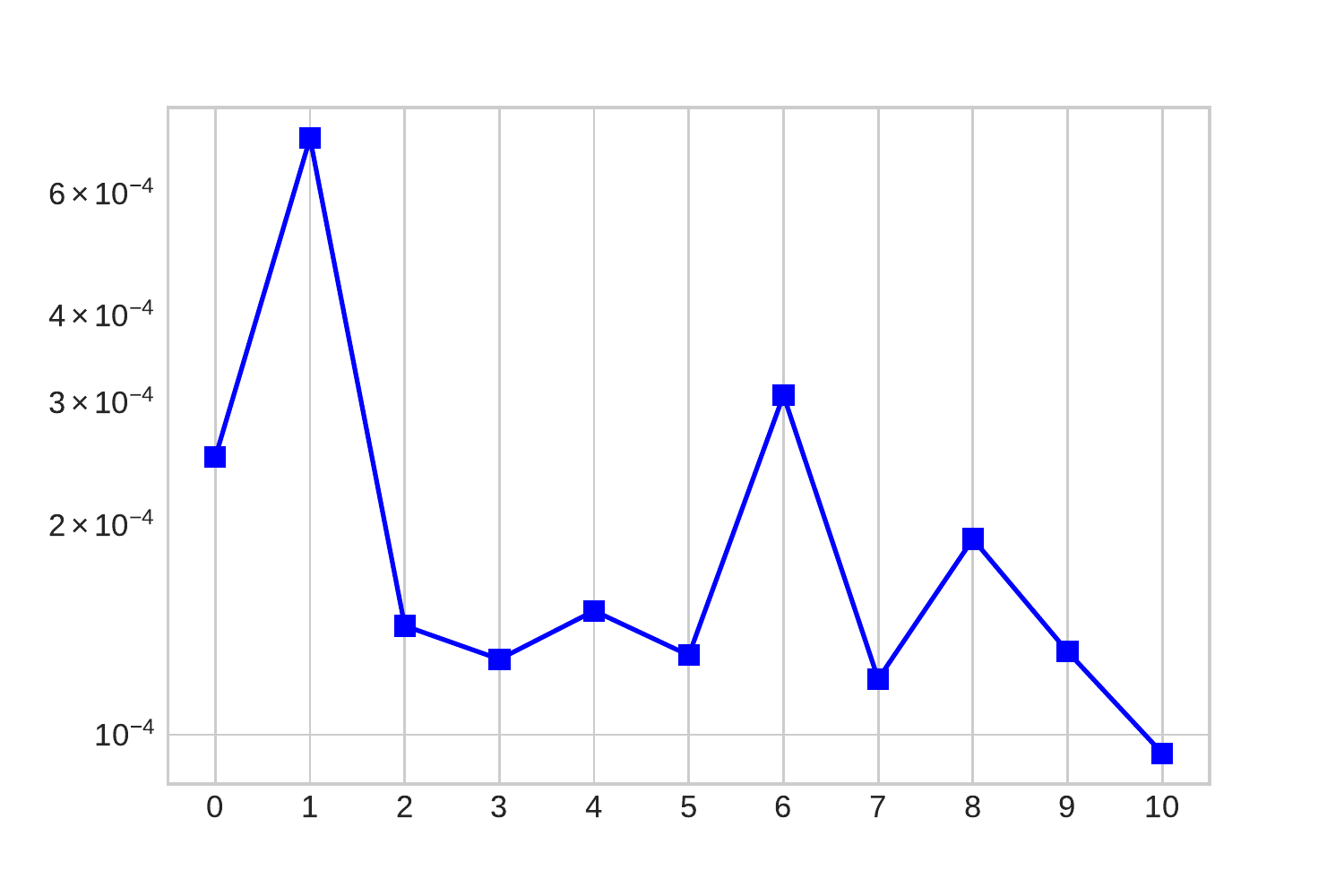}
            \put(35,63){\small G-FIPINN}
         \end{overpic}
      \end{center}
      \vspace{-0.2cm}
      \caption{Estimated failure probabilities over restarts   for the high dimensional case.}
      \label{high_dimensional_failure_probability}
      \end{figure}

In order to make it more clear why our algorithm works, we also plot the final sample distribution in the space of $x_8$ and $x_9$ in Fig.\ref{high_dimensional_samples}.  We can clearly see that samples are concentrated in the failure region, which will give us more useful data for retraining the network. We plot the training loss and the  predicted error in Fig.\ref{high_dimensional_compared_error}  in comparison to the standard PINNs error during training to examine the properties of restart strategy. The restart   strategy will undoubtedly improve the collocation dataset and force the samples to get closer to the failure regions. As a result, after each restart, the loss will suddenly increase before declining.   The system will become more reliable than a vanilla PINN, as shown in Fig.\ref{high_dimensional_failure_probability}, as the predicted error will decline much more quickly.  This is in line with the earlier justification.

\section{Conclusion}
\label{conclusion}
In this work, we have presented two novel extensions to FI-PINNs. The first extension consist in combining with a re-sampling technique, so that the new algorithm can maintain a constant training size. The second extension is to present the subset simulation algorithm as the posterior model (instead of the truncated Gaussian model) for estimating the error indicator, which can more effectively estimate the failure probability and generate new effective training points in the failure region. We investigate the performance of the new approach using several challenging problems, and numerical experiments demonstrate a significant improvement over the original algorithm. In our future work, we shall test the our new algorithm for more challenging problems such as nonlinear time-dependent problems, for which the time domain should be carefully treated \cite{wang2022respecting}. We shall also try to apply the algorithm for operator-learning tasks.

\end{document}